\documentclass[10pt]{article}
\usepackage{amsmath,amsthm,amsfonts,amscd,amssymb,eucal,latexsym}

\renewcommand{\baselinestretch}{1.1}
\begin{document}
\def\pp{\psi\vphantom A}
\def\ff{\varphi\vphantom A}
\def\ov{\overline}
\def\O{{\mathcal O}}
\def\var{\varepsilon}
\def\h#1{\hbox{\rm #1}}
\font\got=eufm10 scaled\magstephalf
\def\Ag{\h{\got A}}
\def\Fg{\h{\got F}}
\def\lam{\lambda}
\def\T{{\cal T}}
\def\C{{\cal C}}
\def\A{{\cal A}}
\def\R{{\cal R}}
\def\B{{\cal B}}
\def\M{{\cal M}}
\def\D{{\cal D}}
\def\N{{\cal N}}
\def\L{{\cal L}}

\def\Bbb{\mathbb}
\def\R{{\cal R}}
\def\A{{\cal A}}
\def\P{{\mathbb P}}
\def\T{{\mathbb T}}
\def\N{{\cal N}}
\def\C{{\cal C}}
\def\d{\hbox{\rm $\,$d}}
\def\End{\hbox{\rm End}}
\def\Aut{\hbox{\rm Aut}}
\def\O{{\cal O}}
\def\H{{\cal H}}
\def\I{{\cal I}}
\def\cC{{\cal C}}
\def\E{{\cal E}}
\def\cT{{\cal T}}
\def\N{{\cal N}}
\def\M{{\cal M}}
\def\D{{\cal D}}
\def\cB{{\cal B}}
\def\K{{\cal K}}
\def\B{{\cal B}}
\def\U{{\cal U}}
\def\V{{\cal V}}
\def\Z{{\mathbb Z}}
\def\ov{\overline}
\def\Om{\Omega}
\def\var{\varepsilon}
\font\got=eufm10 scaled\magstephalf
\def\Ag{\h{\got A}}
\def\Fg{\h{\got F}}
\def\lam{\lambda} 

\allowdisplaybreaks

\font\eightrm=cmr8

{
\renewcommand{\baselinestretch}{1}
\setlength{\textwidth}{15.3cm}
\setlength{\oddsidemargin}{0.6cm}
\setlength{\evensidemargin}{0.6cm}

\title{Noncommutative pressure and the variational principle in
Cuntz--Krieger--type $C^*$--algebras  }
\author{David  Kerr\\ Department of Mathematics\\ University of Toronto\\
Toronto, Ontario, Canada M5S 3G3\\ \\
Claudia Pinzari$^*$\\ Dipartimento di Matematica\\ Universit\`a di Roma Tor
Vergata\\00133 Roma, Italy}
\date{}
\maketitle

\begin{abstract}
Let $a$ be a self-adjoint element of an exact $C^*$--algebra 
$\mathcal{A}$ and $\theta : \mathcal{A} \to \mathcal{A}$ a contractive
completely positive map. We define a notion of dynamical
pressure $P_\theta (a)$ which adopts Voiculescu's approximation
approach to noncommutative entropy and extends the 
Voiculescu--Brown topological entropy and Neshveyev and
St{\o}rmer 
unital-nuclear pressure.
A variational inequality bounding $P_\theta (a)$ below by
the free energies $h_\sigma (\theta ) + \sigma (a)$
with respect to the Sauvageot--Thouvenot entropy $h_\sigma (\theta )$ 
is established in two stages via the introduction of a local state 
approximation entropy, whose associated free energies function as an 
intermediate term. 

Pimsner $C^*$--algebras furnish a framework for investigating
the variational principle, which asserts the equality of 
$P_\theta (a)$ with the supremum of the free energies over all
$\theta$--invariant states. In one direction we extend Brown's result 
on the constancy of the Voiculescu--Brown entropy upon passing to 
the crossed product, and in another we show that the pressure of a
self-adjoint element over the Markov subshift underlying the 
canonical map on the Cuntz--Kreiger algebra $\mathcal{O}_A$ is
equal to its classical pressure. The latter result is extended
to a more general setting comprising an expanded class of 
Cuntz--Krieger-type Pimsner algebras, leading to the variational 
principle for self-adjoint elements in a diagonal subalgebra.
Equilibrium states are constructed from KMS states
under certain conditions in the case of Cuntz--Krieger algebras.
\end{abstract}
}

\vfill

\noindent $^*$Supported by  NATO--CNR.\eject

\section{Introduction}

Let $(X, \theta)$ be a compact topological dynamical system.
The variational principle
in classical ergodic theory, established in full generality by 
Walters \cite{W}, asserts that the topological
pressure of a real-valued continuous function $a$
on $X$ is the supremum of the free energies, i.e., of sums
of the form $h_\sigma(\theta)+\sigma(a)$
where $\sigma$ ranges over all $\theta$--invariant measures of $X$ and
$h_\sigma(\theta)$ stands 
for the (Kolmogorov--Sinai) measure-theoretic entropy of $\theta$. 
Topological pressure was introduced by Walters \cite{W} as a dynamical 
abstraction of the statistical mechanical concept of pressure  
defined as the logarithmic partition function 
density under a thermodynamic limit. Adapting the
approaches of Adler, Konheim, and McAndrew \cite{AKM}, Bowen 
\cite{Bowen2}, and Dinaburg \cite{Din}
to topological entropy, Walter's definition functions not by invoking
a specific sequence of finite subsystems, as in the thermodynamic
notion of entropy or pressure density, but rather samples over 
the dynamical limits of all finite subsystems. Thus, for a lattice
system, the thermodynamic limit is reconceptualized as a dynamical
limit with space translation generating the sequence of subsytems,
and the variational principle for translation-invariant lattice 
systems (see Ruelle \cite{R1}) is subsumed into Walters' 
general result. 

Our ultimate goal
is to investigate the variational principle in a 
noncommutative dynamical setting
which, in analogy to the classical case, captures the 
shift-invariant lattice system model of quantum thermodynamics as
a special instance. Compared to the topological situation, 
noncommutative dynamics presents a much less definitive 
state of affairs for a theory of entropy and pressure. For 
instance, various alternative notions of entropy are available,
from Voiculescu's approximation definition \cite{V} to Connes, Narnhofer, and 
Thirring's \cite{CNT} and Sauvageot and Thouvenot's
\cite{ST} physically motivated approaches
in which the system is observed via Abelian models (see below). 
Recently St\o rmer and Neshveyev, working with a definition of
pressure for unital nuclear $C^*$--algebras and the 
Connes--Narnhofer--Thirring (henceforth abbreviated CNT) entropy, have 
obtained a variational
principle for a class of asymptotically Abelian automorphisms
of AF $C^*$--algebras. 
We will work within the domain of exact $C^*$--algebras, 
replacing the space $X$ by an exact $C^*$--algebra $\A$ and  
$\theta$ by a contractive completely positive self-map of $\A$ and
taking 
the potential to be a fixed self-adjoint element $a$ of $\A$,
and we will establish the variational principle for a class
of $C^*$--dynamical systems which are generally not asymptotically
Abelian.

In Sect. 2 we introduce a notion of pressure for $a$ following
Voiculescu's approach to topological entropy for unital nuclear 
$C^*$--algebras, recently extended
to exact $C^*$--algebras by Brown \cite{B, V}. 
Thus our corresponding partition function is computed by means of 
an optimal approximation, in some sense, of an embedding of 
$\A$ into some $\B(\H)$ via factorizations through finite-dimensional
$C^*$--algebras. 
Our definition reduces to Walters' pressure when the system
$(\A,\theta)$ arises from a topological dynamical system over a compact
space, and also to the pressure introduced by 
Neshveyev and St\o rmer \cite{NS} for unital nuclear $C^*$--algebras.
One advantage of this more general framework
is the yield of an immediate proof of the property
that pressure decreases when taking $C^*$--subalgebras, a fact
that has been already pointed out by Brown for
topological entropy \cite{B}.
 
Among a few other basic properties which easily carry over from 
classical pressure and Voiculescu--Brown topological entropy
or from Neshveyev--St\o rmer pressure, we establish in Sect. 3 
the property of subadditivity in tensor product $C^*$--algebras,
i.e., that the pressure of an element of the form $a\otimes 1+1\otimes b$
with respect to a tensor product map is bounded by the sum of the 
pressures of $a$ and $b$. This fact already implies, in the classical
case, that pressure is a subadditive function. We don't know,
however, whether this still holds for noncommutative pressure.

We next approach the variational principle, first focusing on
a variational inequality which asserts that
that the free energy in a given state is bounded above by the pressure.
In the nuclear case the CNT entropy provides
one natural candidate for defining the free energy,
and indeed the corresponding variational inequality
holds \cite{NS}. In our setting we substitute the Sauvageot--Thouvenot 
entropy, which is defined for unital $C^*$--algebras and
reduces to the CNT entropy in the nuclear case \cite{ST}.
In Sect. 4 we introduce, as an alternative, a measure-theoretic entropy 
for exact $C^*$--algebras which adopts the approximation framework of 
Voiculescu's topological entropy, with the logarithm of the 
rank of the local finite-dimensional
algebra being replaced by the entropy of the induced state on the local
algebra (see Choda \cite{Ch} for the nuclear analogue).
We show that this local state approximation entropy
reduces to the Kolmogorov--Sinai entropy in the classical case,
is a concave function of the invariant state, and majorizes
the Sauvageot--Thouvenot entropy.
The variational inequality is shown to hold if the free energy is
defined via the local state approximation entropy, and as a corollary
we obtain the same inequality using the Sauvageot--Thouvenot entropy. 

In Sect.\ 5 we examine pressure in Cuntz--Krieger algebras
$\O_A$ and crossed product $C^*$--algebras
by a single automorphism $\A\rtimes_\alpha{\Bbb Z}$.  
In the former case we compute the pressure of a self-adjoint 
element $f$ of the canonical Abelian subalgebra of continuous 
functions on the underlying Markov subshift with respect to the natural
unital completely positive map $\theta$ of $\O_A$, with the result
that it equals the classical pressure with respect to the shift epimorphism.
This fact has the consequence that equilibrium states (i.e.,
$\theta$--invariant states
whose free energy reaches the pressure) exist. In particular, 
we recover in the case  
$f=0$ Boca and Goldstein's computation of the topological entropy
of $\theta$ \cite{BG}.
The class of crossed products algebras 
can be regarded, as far as the variational principle is concerned,
as a structurally extreme opposite of that of Cuntz--Krieger algebras.
We generalize Brown's result on the constancy
of topological entropy, so that if $a$ is a self-adjoint element
of $\A$ and $u$ is a unitary in the crossed product implementing $\alpha$,
the pressure of $a$ computed with respect to $\alpha$ in $\A$ 
or $\text{Ad}\,u$ in $\A\rtimes_\alpha{\Bbb Z}$ is the same. 

Regarding $\O_A$ or $\A\rtimes_\alpha{\Bbb Z}$ as a particular case
of the Pimsner $C^*$--algebra $\O_X$ \cite{P} associated to a finitely
generated Hilbert bimodule $X$ over a unital exact $C^*$--algebra
$\A$ 
leads to the problem of investigating the variational principle
in $\O_X$. In Sect. 6 we introduce conditions on 
$X$ which stress
the Cuntz--Krieger behaviour of $\O_X$ rather than the
the crossed product character, and we show that under these conditions
the variational principle holds with the free energy defined
using the Sauvageot--Thouvenot entropy. 
The dynamics here are defined by a unital completely
positive map of $\O_X$ implemented by a basis of the bimodule.
Our main assumptions are the following. First we assume that the left 
action of
$\A$ on $X$ is defined diagonally by a finite set of
endomorphisms
of $\A$.
Then we restrict the space of potentials, selecting self-adjoint elements 
which lie in a ``diagonal subalgebra'' $\D$ of $\O_X$, which is a
noncommutative 
analogue of the canonical maximal Abelian subalgebra of $\O_A$.
Finally we assume that the topological
entropy of the defining set of endomorphisms of $\A$ is zero.
This is the case if, e.g., $\A$ is an inductive limit of finite-dimensional
$C^*$--algebras which are left invariant by the endomorphisms.
This last assumption makes it possible to compute explicitly the pressure
of a potential $a$ in $\D$ which commutes with both the images of $1$ and 
$a$ itself under sufficiently many iterates of the defining endomorphisms.
This is in fact the main step which leads to the proof of the variational
principle.
We also consider a
subclass of potentials of $\D$ for which
 equilibrium
states exist.

In the last section we touch on the problem of the relationship
between the KMS condition and equilibrium, concentrating on the class of
Cuntz--Krieger algebras $\O_A$. To every potential $f\in\C(\Lambda_A)$ 
we associate a one-parameter automorphism group of $\O_A$,
and we show that, if the variation of $f$ is small enough 
and $A$ is aperiodic, the KMS states with respect to this group
are in bijective correspondence with
positive eigenvectors of the Banach space
adjoint $\L_f^*$ of the Ruelle operator $\L_f$ on $\C(\Lambda_A)$. 
A classical
theorem by Ruelle asserts that
if $f$ is H\"older continuous, both $\L_f$ and 
$\L_f^*$
have unique positive eigenvalues, say $h$ and $\mu$, respectively. This
result led Ruelle to a proof of the uniqueness of the equilibrium 
measure for the shift space, which can be identified with the measure
whose Radon--Nykodim derivative with respect to
$\mu$ is $h$ \cite{R68, Bowen, W75}.
We show that, on $\O_A$,  $\mu$ extends naturally  to the unique KMS state
at inverse temperature $1$ and $\nu$ to an equilibrium state of 
$(\O_A,\theta, f)$.

\section{Noncommutative approximation pressure}

\noindent{\it 1.\  Unital exact $C^*$--algebras}\medskip

In this section, unless otherwise stated, $\A$ is a unital exact
$C^*$--algebra, $\theta$ is a unital completely positive map of $\A$
and $a\in \A$ is a self-adjoint element. The collection of finite subsets 
of $\A$ will be denoted by $Pf(\A )$.
We define the
pressure of $a$ by 
approximation through finite-dimensional $C^*$--algebras in the following
way. Let $\pi: \A\to\B(\H)$ be a faithful unital $^*$--representation on a
Hilbert space. 
Since $\A$ is exact and hence nuclearly embeddable \cite{Ki, Was}, 
for any finite subset $\Omega\subset\A$ and for any
$\delta>0$ there is a finite-dimensional $\B$ and unital completely 
positive (henceforth abbreviated as u.c.p.) maps 
$\phi: \A\to\B$ and $\psi: \B\to\B(\H)$
such that $\|(\psi\circ\phi )(x)-\pi(x)\|<\delta$ for all $x\in\Omega$.
We denote by $\text{CPA}(\pi,\Omega, \delta)$ the set of all such
$(\phi,
\psi,B)$. We emphasize that the maps of $\text{CPA}(\pi,\Omega,\delta)$
are unital. 
We set 
\begin{gather*}
\Omega^{(n)}:=\Omega\cup\dots\cup\theta^{n-1}(\Omega),\\
a^{(n)}=\sum_{j=0}^{n-1}\theta^j(a).
\end{gather*} 
We define the {\it partition function}
$$Z_{\theta ,n}(\pi,
a,\Omega,\delta):=
\inf\{\text{Tr }e^{\phi(a^{(n)})} :
(\phi,\psi,\B)\in \text{CPA}(\pi, \Omega^{(n)},\delta)\}$$
where $\text{Tr}$ denotes the trace of $\B$ with the
normalization
$\text{Tr}(e)=1$ for every minimal projection $e\in\B$.
Note that, if $\lambda =\min\text{spec}(a)$, then for any
$(\phi,\psi, \B)\in \text{CPA}(\pi,\Omega^{(n)},\delta)$ we have the
inequality
\begin{gather*}
\text{Tr }e^{\phi(a^{(n)})}\geq e^{\lambda
n}\text{rank}(\B)
\end{gather*}
and so
\begin{gather*}
Z_{\theta ,n}(\pi,a,\Omega,\delta )\geq e^{\lambda
n} rcp(\pi,\Omega^{(n)},\delta),
\end{gather*}
where $rcp$ stands for the Voiculescu--Brown $\delta$--rank
\cite{B, V}. In particular, $$Z_{\theta ,n}(\pi,a,\Omega,\delta)>0.$$
Define 
\begin{gather*}
P_{\theta}(\pi,a,\Omega,\delta):=\limsup_n\frac{1}{n}\log
Z_{\theta ,n}(\pi,a,\Omega,\delta),\\
P_{\theta}(\pi, a, \Omega)=\sup_{\delta>0}P_{\theta}(\pi, a, \Omega,
\delta),\\
P_\theta(\pi, a)=\sup_{\Omega\in Pf(\A)}P_\theta(\pi, a,
\Omega).
\end{gather*}
We will refer to $P_\theta(\pi, a)$ as the {\it approximation pressure}
(or simply {\it pressure}) of
$a$ (with respect to $\theta$).
Note that, referring to the notation of Brown \cite{B} and Voiculescu
\cite{V}, 
\begin{gather*}
Z_{\theta ,n}(\pi,0,\Omega)= rcp(\pi,\Omega^{(n)},\delta),
\end{gather*}
and so
\begin{gather*}
P_{\theta}(\pi,0,\Omega,\delta)= ht(\pi,\theta,\Omega,\delta),\\
P_{\theta}(\pi,0,\Omega)= ht(\pi,\theta,\Omega),\\
P_{\theta}(\pi,0)= ht(\pi,\theta)=\text{the Voiculescu--Brown
entropy of } \theta. 
\end{gather*}

The first fact that we want to establish is that the partition
function, and therefore the pressure, does not 
depend upon the
representation $\pi$. This will be done by generalizing arguments of Brown
\cite{B}
for the  entropy.\medskip

\noindent{\bf Proposition 2.1.} {\it If $\pi_1$ and $\pi_2$ are faithful
and unital $^*$--representations of $\A$,
$$Z_{\theta ,n}(\pi_1, a, \Omega,\delta)=Z_{\theta ,n}(\pi_2, a,\Omega,
\delta ),$$
and so
$$P_\theta(\pi_1, a)=P_\theta(\pi_2, a).$$}\medskip

\noindent{\it Proof.} Given $\epsilon>0$,
choose $(\phi, \psi,\B)\in \text{CPA}(\pi_1,\Omega^{(n)},\delta)$ such 
that $$\text{Tr
}e^{\phi(a^{(n)})}-Z_{\theta ,n}(\pi_1,a,\Omega, \delta )<\epsilon.$$
Consider the map $\pi_2\circ{\pi_1}^{-1}:\pi_1(\A)\to\B(\H_{\pi_2})$.
Apply Arveson's extension theorem \cite{A} to extend this map to a 
u.c.p.\ map $T:\B(\H_{\pi_1})\to\B(\H_{\pi_2})$. Then 
$(\phi, T\circ\psi, \B)\in \text{CPA}(\pi_2,\Omega^{(n)},\delta)$ and so
we
easily
obtain
$Z_{\theta ,n}(\pi_2,a,\Omega,\delta)\leq
Z_{\theta ,n}(\pi_1,a,\Omega,\delta)$. The opposite inequality follows by
exchanging the roles of $\pi_1$ and $\pi_2$.\hfill $\square$\medskip

{\noindent As} a result of this proposition we can avoid specifying the 
representation $\pi$ in the partition function as well as in the 
approximation pressures.
\bigskip

\noindent{\it 2.\ Unital nuclear $C^*$--algebras}\medskip

\noindent{Let }$\Omega$ be a finite subset of $\A$, $\delta>0$, 
and $n\in{\Bbb N}$.
If $\A$ is a nuclear $C^*$--algebra, in the definition of pressure it is 
more natural to replace $\text{CPA}(\pi,\Omega^{(n)},\delta)$ with the 
set $\text{CPA}_{\text{nuc}}(\pi , \Omega^{(n)},\delta)$ of all triples
$(\phi,\psi,\B)$ where $\phi:\A\to\B$ and $\psi:\B\to\A$ are u.c.p.\ maps
and $\B$ is a finite--dimensional $C^*$--algebra such that
$\|(\psi\circ\phi )(x)-x\|<\delta$ for all $x\in\Omega^{(n)}$.
We thus obtain the corresponding nuclear partition function
$Z^{\text{nuc}}_{\theta ,n}(a,\Omega,\delta )$,
nuclear approximation pressures
$P^{\text{nuc}}_\theta(a,\Omega,\delta)$ and
$P^{\text{nuc}}_\theta(a,\Omega)$, and
nuclear pressure $P^{\text{nuc}}_\theta(a)$, as in \cite{NS}.
\medskip

\noindent{\bf Proposition 2.2.} {\it Let $\A$ be a unital nuclear
$C^*$--algebra faithfully and unitally represented on a Hilbert space
$\H$, and let
$a\in\A$ be a self-adjoint element.
Then for any finite subset $\Omega\subset\A$, $\delta>0$, and $n\in{\Bbb
N}$,
$$Z^{\text{\rm nuc}}_{\theta ,n}(a,\Omega,\delta)=Z_{\theta
,n}(a,\Omega,\delta)$$
and so
$$P^{\text{\rm nuc}}_\theta(a)=P_\theta(a).$$}\medskip

\noindent{\it Proof.}
Our arguments generalize the corresponding arguments of Brown 
(Prop.\ 1.4 of \cite{B})
for the  Voiculescu--Brown entropy.
Fix $\Omega\in Pf(\A)$, $\delta>0$, and $n\in{\Bbb N}$. 
We first note that
$Z_{\theta ,}(a,\Omega,\delta)\leq
Z^{\text{nuc}}_{\theta ,n}(a,\Omega,\delta)$ since 
$\text{CPA}_{\text{nuc}}(\Omega^{(n)},\delta)\subset
\text{CPA}(\pi , \Omega^{(n)},\delta)$.
Given $\epsilon>0$
let $(\phi,\psi,\B)\in \text{CPA}(\pi , \Omega^{(n)},\delta)$ be such 
that $$\text{Tr }e^{\phi(a^{(n)})}-Z_{\theta ,n}(a,\Omega,\delta )
<\epsilon.$$
Choose a triple $(\rho,\sigma, \C)\in \text{CPA}_{\text{nuc}}(\pi ,
\Omega^{(n)},
\delta)$ and consider
a (unital) completely positive extension $\Phi:\B(\H)\to\C$ of 
$\rho$, which exists by Arveson's extension theorem \cite{A}. Then
$(\phi, \sigma\circ\Phi\circ\psi, \B)\in
\text{CPA}_{\text{nuc}}(\pi , \Omega^{(n)}, 2\delta)$ and so we easily 
deduce that
${Z^{\text{nuc}}}_{\theta ,n}(a,\Omega,2\delta )\leq
Z_{\theta ,n}(a,\Omega,\delta ).$\hfill $\square$
\bigskip

\noindent{\it 3.\ Not-necessarily-unital exact $C^*$--algebras}\medskip

\noindent{Let }$\A$ be an exact $C^*$--algebra faithfully 
represented on a Hilbert
space $\H$, and $\theta : \A\to\A$ a completely positive contraction.
Following Brown's approach to topological entropy for 
not-necessarily-unital $C^*$--algebras, we introduce a partition 
function $Z^0$ and corresponding approximation
pressures $P^0$ as in the unital case but with respect to the
expanded collection $\text{CPA}_0(\pi , \Omega^{(n)},\delta)$ of triples
$(\phi,\psi,\B)$ where $\B$ is a finite-dimensional $C^*$--algebra and 
$\phi:\A\to\B$ and $\psi:\B\to\B(\H)$ are c.p.\ contractions such that 
$\|\psi\circ\phi(x)-x\| <\delta$ for all $x\in\Omega^{(n)}$.
$P^0_\theta(\pi, a)$ is still independent of the representation $\pi$.
 We would
like to thank G.\ Gong for pointing out why an equality such as the
one in the claim in the proof of the following proposition should 
hold.\medskip

{\noindent\bf Proposition 2.3.} {\it If $\mathcal A$ is unital and exact
and $\theta$ is u.c.p.\ then $P^0_\theta (a)=P_\theta (a)$.}
\medskip

{\noindent\it Proof.} The inequality $P^0_\theta (a) \leq P_\theta (a)$
follows immediately from the definitions.

To establish the reverse inequality, let $\Omega$ be a finite
subset of the unit ball of $\mathcal{A}$ containing $1$, and suppose
$0\leq\delta\leq\frac{1}{4}$. Let $(\phi , \psi , \mathcal{B}) \in
\text{CPA}_0 (\pi , \Omega^{(n)}, \delta )$, and set 
$b=\phi (1)$. Let $p$ be a spectral projection of $b$ 
such that $b_1 :=bp\geq (1-\sqrt{\delta})p$
and $b_2 := b(1-p) < 1-\sqrt{\delta}$. We claim that
$\psi (b_2) < \sqrt{\delta}$. To see this, suppose to the contrary
that $\| \psi (b_2) \| \geq \sqrt{\delta}$. Since $\psi$ is contractive
we have $\left\| \psi \left(b_1 + \frac{1}{1-\sqrt{\delta}}b_2 \right) 
\right\| \leq \left\| b_1 + \frac{1}{1-\sqrt{\delta}}b_2 \right\| 
\leq 1$. On the other hand, since $\psi (b) \geq 1 - \| 
\psi\circ\phi (1) - 1 \| > 1 - \delta$, the positivity of $\psi$ yields 
\begin{align*}
\psi \left(b_1 + \frac{1}{1-\sqrt{\delta}}b_2 \right) 
&= \psi (b_1 + b_2) + \psi \left( \left( \frac{1}{1-\sqrt{\delta}}
-1\right) b_2\right) \\
&> 1 - \delta + \frac{\sqrt{\delta}}{1-\sqrt{\delta}}\psi (b_2), 
\end{align*}
and since $\left\| \frac{\sqrt{\delta}}{1-\sqrt{\delta}} \psi (b_2) 
\right\| > \delta$ this implies $\left\| \psi \left(b_1 + 
\frac{1}{1-\sqrt{\delta}}b_2 \right) \right\| > 1$, producing a
contradiction and thus establishing the claim. 
 
Observe now that 
\begin{align*}
\| \psi (b_1) - 1 \| & = \| \psi (b) - \psi (b_2) - 1 \| \\
&\leq \| \psi (b) - 1 \| + \| \psi (b_2) \| \\
&< \delta + \sqrt{\delta} \\
&< 2\sqrt{\delta}.
\end{align*}
Thus, if $p$ denotes the support projection of $b_1$, we have
$\| \psi (p) - 1 \| \leq \| \psi (p - b_1) \| + \| \psi (b_1) - 1 \| < 
3\sqrt{\delta}$, and so $\| \psi (p)^2 - \psi (p) \| \leq \| \psi (p) 
(\psi (p) - 1) \| < 3\sqrt{\delta}$.
Appealing to Stinespring's theorem \cite{St, Was} we infer that, for 
all $x\in\mathcal{B}$, $\| \psi (pxp) - \psi (p) \psi (x) \psi (p) \| <
4\sqrt{3}\sqrt[4]{\delta}$ and hence
\begin{align*}
\| \psi (pxp) - \psi (x) \| &\leq \| \psi (pxp) - \psi (p) \psi (x) 
\psi (p) \| + \| \psi (p) \psi (x) \psi (p) - \psi (x) \| \\
&< 4\sqrt{3}\sqrt[4]{\delta} \| x \| + 8 \sqrt{\delta} \| x \| \\
&< 16\sqrt[4]{\delta} \| x \|.
\end{align*}

Set $\mathcal{B'} = p\mathcal{B}p$, and define the u.c.p.\ map $\phi' 
: \mathcal{A} \to \mathcal{B'}$ by $\phi' (x) = b_1^{-\frac{1}{2}}
\phi (x) b_1^{-\frac{1}{2}}$, with $b_1$ now being considered as an 
element of $\mathcal{B}'$. Let $\psi' : \mathcal{B'} \to \mathcal{B}
(\mathcal{H})$ be the u.c.p.\ map given by $\psi' (x) = 
\psi (b_1)^{-\frac{1}{2}} \psi \left( b_1^{\frac{1}{2}}xb_1^{\frac{1}{2}} 
\right) \psi (b_1)^{-\frac{1}{2}}$. If $x\in\Omega^{(n)}$ then 
$\| \psi (p\phi (x)p) - x \| \leq \| \psi (p\phi (x) p - \phi (x)) \| + 
\| \psi\circ\phi (x) - x \| < 16\sqrt[4]{\delta} + \delta < 
17\sqrt[4]{\delta}$, and so estimating as does Brown in \cite{B} we obtain 
$\| \psi' \circ \phi' (x) - \psi (p\phi (x) p) \| < 
14(17\sqrt[4]{\delta})$, whence
\begin{align*}
\| \psi' \circ \phi' (x) - x \| &\leq \| \psi' \circ \phi' (x) - \psi 
(p\phi (x) p) \| + {} \\
& \hspace*{1cm} \| \psi (p\phi (x) p - \phi (x)) \| + 
\| \psi\circ\phi (x) - x \| \\
&< 255\sqrt[4]{\delta}.
\end{align*}
We therefore have $(\phi' , \psi' , \mathcal{B'}) \in \text{CPA}(\pi , 
\Omega^{(n)} , 255\sqrt[4]{\delta} ))$. Also note that
\begin{align*}
\big\| b_1^{-\frac{1}{2}} \phi (a^{(n)})b_1^{-\frac{1}{2}} - p\phi 
(a^{(n)})p \big\| &\leq \big\| b_1^{-\frac{1}{2}} \phi (a^{(n)})
b_1^{-\frac{1}{2}} - b_1^{-\frac{1}{2}}\phi (a^{(n)})p \big\| \\
&\hspace*{1.8cm} + \big\| b_1^{-\frac{1}{2}}\phi (a^{(n)})p - 
p\phi (a^{(n)})p \big\| \\
&\leq n \| a \| \big\| b_1^{-\frac{1}{2}} - p \big\| \Big(
\big\| b_1^{-\frac{1}{2}} \big\| + 1 \Big) \\
&\leq 2n \| a \| \frac{\sqrt{\delta}}{(1-\sqrt{\delta})^2}
\end{align*} 
so that $b_1^{-\frac{1}{2}} \phi (a^{(n)})b_1^{-\frac{1}{2}} \leq p\phi 
(a^{(n)})p + 2n\|a\|\frac{\sqrt{\delta}}{(1-\sqrt{\delta})^2}$ and hence
\begin{align*}
\text{log Tr}_\mathcal{B'} e^{\phi' (a^{(n)})} &\leq  
\text{log Tr}_\mathcal{B'} e^{p{\phi (a^{(n)})}p + 
2n\|a\|\frac{\sqrt{\delta}}{(1-\sqrt{\delta})^2}} \\
&= \text{log Tr}_\mathcal{B'} e^{p{\phi (a^{(n)})}p} +
2n\|a\|\frac{\sqrt{\delta}}{(1-\sqrt{\delta})^2} \\
&\leq \text{log Tr}_\mathcal{B'} pe^{\phi (a^{(n)})}p + 
2n\|a\|\frac{\sqrt{\delta}}{(1-\sqrt{\delta})^2} \\
&\leq \text{log Tr}_\mathcal{B} e^{\phi (a^{(n)})} + 
2n\|a\|\frac{\sqrt{\delta}}{(1-\sqrt{\delta})^2}, 
\end{align*}
with the second last inequality following from Proposition 3.17 of 
\cite{OP}. It follows that 
$$P_{\theta}^0 (\mathcal{A} , \Omega , 
255\sqrt[4]{\delta}) \leq P_{\theta}^0 (\mathcal{A} , 
\Omega , \delta ) + 2\|a\|\frac{\sqrt{\delta}}{(1-\sqrt{\delta})^2},$$ 
from which we conclude that $P_{\theta}^0 (a, \Omega ) \leq P_{\theta} 
(a, \Omega )$. Thus $P_{\theta}^0 (a) \leq P_{\theta} 
(a)$.\hfill $\square$

\section{Main Properties}

It is  natural to ask which properties of the Voiculescu--Brown 
topological entropy
or the classical pressure
carry over to the pressure of a self-adjoint element in a
unital exact $C^*$--algebra.
The following result collects some properties inspired from corresponding 
properties of the classical pressure \cite{W}.
\medskip

\noindent{\bf Proposition 3.1.} {\it Let $a$, $b$ be self-adjoint 
elements of $\A$.
\begin{description}
\item{\rm (a)} If $a\leq b$, $P_\theta(a)\leq P_\theta(b)$;
\item{\rm (b)} if $\lambda\in{\Bbb R}$, $P_\theta(a+\lambda
I)=P_\theta(a)+\lambda$. In particular
$P_\theta(\lambda I)=\lambda+ ht(\theta)$; 
\item{\rm (c)} 
$\min\text{spec}(a)+\text{ht}(\theta)\leq
P_\theta(a)\leq\max\text{spec}(a)+ ht(a)$,
so either $P_\theta(a)<\infty$ for all $a$ or
$P_\theta(a)= ht(\theta)=\infty$ for all $a$;
\item{\rm (d)} if $ht(\theta)<\infty$,
$|P_\theta(a)-P_\theta(b)|\leq\|a-b\|;$
\item{\rm (e)} $P_\theta(ca)\leq c P_\theta(a)$ if $c\geq1$ and
$P_\theta(ca)\geq cP_\theta(a)$ if $c\leq1$;
\item{\rm (f)} $|P_\theta(a)|\leq P_\theta(|a|)$.
\end{description}}\medskip

\noindent{\it Proof.} Properties (a) and (d) can be established
using the Peierls--Bogoliubov inequality (cf.\ Cor.\ 3.15 of \cite{OP}) 
as in the proof of Prop.\ 2.4 in \cite{NS} for the nuclear pressure,
while (b) and (e) are immediate from the definition, (c) follows
from (a) and (b), and (f) follows from (e) 
and (a).\hfill $\square$\medskip

The following facts are also very easy to check.\medskip

\noindent{\bf Proposition 3.2.} {\it
\begin{description}
\item{\rm (a)}
$P_\theta(a)=\frac{1}{r}P_{\theta^r}(a^{(r)}),$ if $r\in{\Bbb
N}$;
\item{\rm (b)} 
if  $\theta$ is an automorphism,
$P_\theta(a)=P_{{\theta}^{-1}}(a)$;
\item{\rm (c)} $P_\theta(a+\theta(b)-b) = P_\theta(a)$; 
\item{\rm (d)}
If $\theta$ is an automorphism, $P_\theta(\theta(a))=P_\theta(a)$.
\end{description} }\medskip

\noindent{\it Proof.} The proofs of Propostion 2.4(v)(ii)in \cite{NS} 
for the nuclear approximation pressure can be adapted to establish (a) 
and (c), respectively. To establish (b), we need only note that 
$Z_{\theta ,n}(a,\Omega,\delta)=Z_{\theta^{-1} ,n}(a,\Omega,\delta)$
follows from the observation that $(\phi,\psi,\B)\in
\text{CPA}(\pi,\Omega\cup\dots\cup{(\theta^{-1})}^{n-1}\Omega,
\delta)$
if and only if
$(\phi\circ\theta^{-n+1},
\tilde{\theta}^{n-1}\circ\psi,
\B)\in
\text{CPA}(\pi, \Omega\cup\dots\cup\theta^{n-1}(\Omega),\delta)$,
where $\tilde{\theta} : \mathcal{B}(\mathcal{H}) \to 
\mathcal{B}(\mathcal{H})$ is a u.c.p.\ extension of $\pi\circ\theta\circ
\pi^{-1} : \pi (\mathcal{A}) \to \mathcal{B}(\mathcal{H})$ whose
existence is guaranteed by Arveson's Extension Theorem.
To show (d), we can take $b=a$ in (c). 
\hfill$\square$\medskip

Next we discuss a few properties of the Voiculescu--Brown entropy 
which easily carry over to pressure.\medskip

\noindent{\bf Proposition 3.3.} {\it (Monotonicity)\, Let $\A_0\subset\A$ 
be a $\theta$--invariant $C^*$--subalgebra (i.e., $\theta(\A_0)\subset\A_0$) 
containing $a$.
Then
$$P_{\theta\upharpoonright\A_0}(a)\leq P_{\theta}(a).$$}\medskip

We also have a Kolmogorov--Sinai-type result.\medskip

\noindent{\bf Proposition 3.4.} {\it Let $\{\Omega_\iota : \iota\in I\}$
be a net of finite subsets of $\A$ such that
$\bigcup_{\iota\in I}\bigcup_{j\in{\Bbb N}}\theta^j(\Omega_\iota)$ is total.
Then 
$$P_\theta(a)=\lim_\iota P_\theta(a,\Omega_\iota).$$}\medskip

\noindent{\it Proof.} It is clear that for $\Omega_1\subset\Omega_2$,
$P_\theta(a,\Omega_1)\leq P_\theta(a,\Omega_2)$, and so the limit on the
r.h.s.\ exists and is bounded by the l.h.s.
Let $\Omega\in Pf(\A )$ and $\delta>0$.
Consider $\iota\in I$ and $N\in{\Bbb N}$ such that for any $x\in\Omega$
there is
$x'=\sum_{r\in F, j\leq N}\lambda_{r,j,x} \theta^j({y}_{r,x})$
with $F$ a finite set, $y_{r,x}\in\Omega_\iota$ such that
$\|x-x'\|<\delta$. Set $\delta'=\frac{\delta}{(N+1)\text{Card}
(\Omega_\iota)\max_{r,j,x}|\lambda_{r,j,x}|}$.
For each $n\in{\Bbb N}$ take a triple 
$(\phi,\psi,\B)\in \text{CPA}(\Omega_\iota^{(n+N+1)},
\delta')$
such that
$$\text{Tr
}e^{\phi(a^{(n+N+1)})}<
2Z_{\theta ,n+N+1}(a,\Omega_\iota,\delta')$$
One can easily show that
$(\phi,\psi,\B)\in \text{CPA}(\Omega^{(n)}, 3\delta)$,
and so
\begin{align*}
Z_{\theta ,n}(a,\Omega,3\delta )\leq
\text{Tr }e^{\phi(a^{(n)})}
&\leq \text{Tr
}e^{\phi(a^{(n+N+1)})+(N+1)\|a\|} \\
&\leq 2e^{(N+1)\|a\|}Z_{\theta
,n+N+1}(a,\Omega_\iota,\delta'),
\end{align*}
from which we conclude that $P_\theta(a,\Omega)\leq\lim_\iota 
P_\theta(a,\Omega_\iota)$.\hfill $\square$
\medskip

The next proposition  gives a weak 
version of subadditivity in a tensor product $C^*$--algebra. 
It also extends the entropy tensor product 
inequalities from \cite{V} to pressure. Note that 
the class of exact $C^*$--algebras is closed under taking minimal 
tensor products \cite{Ki}.\medskip

\noindent{\bf Proposition 3.5.} {\it Let $\theta_1 : \A_1 \to \A_1$ and
$\theta_2 : \A_2 \to \A_2$ be u.c.p.\ maps and 
let $a_1$ and $a_2$ be self-adjoint elements of $\A_1$ and $\A_2$,
respectively. Let $\theta : \A_1 \otimes_{\text{\normalfont min}} \A_2 
\to \A_1 \otimes_{\text{\normalfont min}} \A_2$ be the (u.c.p.)\ extension 
of the map $\theta_1 \otimes\theta_2 : \A_1 \otimes \A_2 \to \A_1 
\otimes \A_2$ on the algebraic tensor product. Then
\begin{eqnarray*}
& P_\theta (a_1 \otimes 1 + 1 \otimes a_2 ) \leq P_{\theta_1} (a_1 ) 
+ P_{\theta_2} (a_2 ),
\end{eqnarray*}}\medskip

\noindent{\it Proof.} Let   $\A_1$ and $\A_2$
be faithfully represented on Hilbert spaces
$\H_1$ and $\H_2$ respectively. 
 Then $\A_1 \otimes
\A_2$ is faithfully represented on  
$\mathcal{H}_1 \otimes \mathcal{H}_2$.

Let $\Omega_1 \in Pf(\A_1)$, $\Omega_2 \in Pf(\A_2)$, and $\delta_1,
\delta_2 > 0$. Set
$M = \mbox{max}\{ \| x \| : x \in \Omega_1 \cup \Omega_2 \}$. Suppose
$(\phi_j, \psi_j,\B_j ) \in \text{CPA}(\A_j, \Omega_j^{(n)},
\delta_j )$ for $j=1,2$. Let $\phi : \A_1 \otimes_{\text{min}} \A_2 \to
\B_1
\otimes \B_2$ be the (u.c.p.)\ extension of the map
$\phi_1 \otimes \phi_2 : \A_1 \otimes \A_2 \to \B_1 \otimes \B_2$.
If $x_1 \in \A_1$ and $x_2 \in \A_2$ then
\begin{eqnarray*}
\lefteqn{\| ((\psi_1 \otimes \psi_2 )\circ (\phi_1 \otimes \phi_2)) (x_1 
\otimes x_2 ) - x_1 \otimes x_2 \|} \hspace*{1cm}\\
& = & \| (\psi_1 \circ \phi_1 )(x_1) \otimes (\psi_2 \circ \phi_2 )
(x_2) - x_1 \otimes x_2 \| \\ 
& \leq & \| (\psi_1 \circ \phi_1 ) (x_1) - x_1 
\| \| x_2 \| + \| x_1 \| \| (\psi_2 \circ \phi_2 ) (x_2) - x_2 \|
\end{eqnarray*}
and so  $(\phi, \psi_1 \otimes \psi_2, \B_1 \otimes \B_2) \in
\text{CPA}((\Omega_1 \otimes \Omega_2 )^{(n)}, M(\delta_1 + \delta_2 ))$. 

Let $(e_k^{1})_{k=1}^{\text{rank}(\B_1)}$ and
$(e_l^{2})_{l=1}^{\text{rank}(\B_2)}$ 
be maximal sets of pairwise
orthogonal minimal spectral projections for $\phi ((a_1 \otimes 1
)^{(n)})$
and $\phi ((1 \otimes a_2 )^{(n)})$, respectively. Then 
$(e_k^{1} \otimes e_l^{2})_{1 \leq k \leq \text{rank}(\B_1),
1 \leq l \leq \text{rank}(\B_2)} $ is a maximal set of 
pairwise orthogonal minimal spectral projections for
$\phi ((a_1 \otimes 1)^{(n)} + (1 \otimes a_2 )^{(n)})$, and so
\begin{eqnarray*}
\lefteqn{\text{Tr}_{\B_1 \otimes \B_2}\,
e^{\phi\left(\left(a_1 \otimes 1\right)^{(n)} + \left(1 \otimes 
a_2 \right)^{(n)}\right)}} \hspace*{1.5cm} \\
& = & \sum_{k,l}\, e^{\text{Tr}_{\B_1 \otimes \B_2}\left(\left(
e_k^{1} \otimes e_l^{2}\right) \phi \left(\left(a_1 \otimes 
1\right)^{(n)} + \left(1 \otimes a_2 \right)^{(n)}\right)\right)} \\
& = & \sum_{k,l}\, e^{\text{Tr}_{\B_1}\left(e_k^{1} \phi_1
\left({a_1}^{(n)}\right)\right) +
\text{Tr}_{\B_2}\left(e_l^{2} \phi_2 \left(
{a_2}^{(n)}\right)\right) } \\
& = & \sum_{k} e^{\text{Tr}_{\B_1}\left(e_k^{1} \phi_1
\left({a_1}^{(n)}\right)\right) }
\sum_{l}\, e^{\text{Tr}_{\B_2}\left(e_l^{2} \phi_2
\left({a_2}^{(n)}\right)\right)} \\
& = & \text{Tr}_{\B_1}\, e^{\phi_1\left({a_1}^{(n)}\right)}\,
\text{Tr}_{\B_2}\, e^{\phi_2\left({a_2}^{(n)}\right)}.
\end{eqnarray*}
Therefore
\begin{eqnarray*}
\lefteqn{Z_{\theta ,n}\left(a_1 \otimes 1 + 1 \otimes a_2 ,
\Omega_1 \otimes \Omega_2 , M\left(\delta_1 + \delta_2 \right)\right)} 
\hspace*{2cm} \\
\hspace*{2cm} & \leq & Z_{\theta_1 ,n}(a_1 , \Omega_1 , \delta_1 )
\hspace*{0.5mm} Z_{\theta_2 ,n}(a_2 , \Omega_2 , \delta_2 ), 
\end{eqnarray*}
and since $\A_1 \otimes \A_2$ is dense in $\A_1 \otimes_{\text{min}} \A_2$ 
it follows from Prop.\ 3.4 that
\begin{eqnarray*}
& P_{\theta}(a_1 \otimes 1 + 1 \otimes a_2 ) \leq P_{\theta_1}(a_1 ) 
+ P_{\theta_2}(a_2 ).
\end{eqnarray*}
\hfill $\square$\medskip

If ($X$, $T$) is a topological dynamical system
over compact metric space $X$, and $a$ is a real-valued continuous
function over $X$, 
we denote by $p_T(a)$ the classical topological pressure of $a$ 
considered by Walters  \cite{W}.\medskip

\noindent{\bf Proposition 3.6.} {\it Let $X$ be a compact metric space, 
$T:X\to X$ a continuous function and $a$ a real-valued continuous function on $X$.
Then $$P_{\theta_T}(a)=p_T(a)$$ 
where $\theta_T$ is
the $^*$--homomorphism of $\C(X)$ defined by $\theta_T(f)(x)=f(Tx)$.}\medskip

\noindent{\it Proof.} Since the nuclear and exact approximation
pressures agree for nuclear $C^*$--algebras by Prop.\ 2.2, we can 
appeal to Remark 2.3 of \cite{NS}.\hfill $\square$\medskip

\noindent{\it Remark.}
It is natural to ask whether, as in the
classical situation, the function $a\to P_\theta(a)$ is convex or 
subadditive in the case where $ht(\theta)$ is
finite. 
We just note that in the classical situation tensor
product
subadditivity  combined with monotonicity of the classical pressure when
passing to
closed invariant subspaces implies
subadditivity.

\medskip

\section{Entropy and variational inequalities}

\noindent Our next aim is to establish a variational inequality
bounding the free energy in a given state by the pressure. We first 
introduce a notion of exact--$C^*$--algebraic entropy 
with respect to an invariant state which 
adopts the approximation framework of Voiculescu's 
topological definition \cite{V}, but exercises the entropy 
of the induced local state instead of the logarithm of the rank of 
the local algebra (see \cite{Ch} for the nuclear analogue). 
The local state approximation entropy yields as a 
straightforward consequence of its definition the desired variational
inequality (Prop.\ 4.14), and since it majorizes both
the Sauvageot--Thouvenot and CNT entropies (Prop.\ 4.10) the
inequality will also hold upon substituting either of the latter 
as the entropy term in the free energy. To conclude the
section we collect some facts about the Sauvageot--Thouvenot entropy 
which will be needed in Sect. 6.

Let $\A$ be a unital exact $C^*$--algebra, $\theta : \A \to \A$ a 
u.c.p.\ map, and $\sigma$ a $\theta$--invariant state on $\A$. 
Let $\D$ be an injective $C^*$--algebra and $\iota : \A \to \D$ 
a unital complete order (henceforth abbreviated u.c.o.) embedding. 
For $\Omega\in Pf(\A)$ and $\delta > 0$ we denote 
by $\text{CPA}(\iota , \Omega , \delta )$ the set of all triples 
$(\phi , \psi , \B)$ where $\B$ is a finite-dimensional $C^*$-algebra 
and $\phi : \A \to \B$ and $\psi : \B \to \D$ are u.c.p.\ maps such that
$\| (\psi\circ\phi ) (x) - \iota (x) \| < \delta$ for all $x\in
\Omega$. Since $\A$ is nuclearly embeddable \cite{Ki}, the set 
$\text{CPA} (\iota , \Omega , \delta )$ is non-empty.  
Denote by $\mathfrak{E}(\sigma , \iota )$ the 
set of all states $\omega$ on $\D$ which extend 
the state $\sigma\circ\iota^{-1}$ on $\iota (\A)$.\medskip

\noindent{\bf Definition 4.1.} {\it If $\omega$ is a state on $\D$, 
$\Omega\in Pf(\A)$, and $\delta > 0$, we define the 
completely positive $\delta$--entropy
$$ cpe(\iota , \omega , \Omega , \delta ) = \inf 
\left\{ S(\omega\circ\psi ) : (\phi , \psi , \B)\in \text{CPA}(\iota , 
\Omega , \delta )\right\} $$
of $\Omega$ with respect to $(\iota , \omega )$, and for 
$\omega\in\mathfrak{E}(\sigma , \iota )$ we define the dynamical 
entropies
\begin{align*}
hm_\sigma (\theta , \iota , \omega , \Omega , \delta ) &= 
\limsup_{n\to\infty} \frac{1}{n} cpe(\iota , \omega , \Omega^{(n)}, 
\delta ) \\
hm_\sigma (\theta , \iota , \omega , \Omega ) &= \sup_{\delta >0} 
hm_\sigma (\theta , \iota , \omega , \Omega , \delta ) \\
hm_\sigma (\theta , \iota , \omega ) &= \sup_{\Omega\in Pf(\A)} 
hm_\sigma (\theta , \iota , \omega , \Omega ) \\
hm_\sigma (\theta , \iota ) &= \sup_{\omega\in\mathfrak{E}(\D, 
\sigma , \iota ) } hm_\sigma (\theta , \iota , \omega ).
\end{align*}}\medskip

We will refer to $hm_\sigma(\theta,\iota)$ as the {\it local state
approximation entropy} of $\theta$.\medskip

\noindent{\bf Proposition 4.2. } {\it If $\iota_1 : \A \to 
\D_1$ and $\iota_2 : \A \to \D_2$ are u.c.o.\ embeddings 
into injective $C^*$--algebras $\D_1$ and $\D_2$ then
$$ hm_\sigma (\theta , \iota_1 ) = hm_\sigma (\theta ,
\iota_2 ). $$}\medskip

\noindent{\it Proof.} Since $\D_1$ is injective we can extend
the map $\iota_1 \circ\iota_2^{-1} : \iota_2 (\A)
\to \D_1$ to a u.c.p.\ map $\Upsilon : \D_2 \to \D_1$. Let
$\omega_1 \in \mathfrak{E}(\sigma , \iota_1 )$ and define
$\omega_2 \in \mathfrak{E}(\sigma , \iota_2 )$ by $\omega_2 =
\omega_1 \circ\Upsilon$. Let $\Omega\in Pf(\A)$ and $\delta > 0$,
and suppose $(\phi , \psi , \B) \in \text{CPA}(\iota_2 , 
\Omega^{(n)} , \delta )$. Then for all $x\in\Omega^{(n)}$,
$$ \| ((\Upsilon\circ\psi )\circ\phi )(x) - \iota_1 (x) \| 
= \| \Upsilon ((\psi\circ\phi )(x) - \iota_2 (x)) \| 
\leq \| (\psi\circ\phi )(x) - \iota_2 (x) \| 
< \delta ,$$
so that $(\phi , \Upsilon\circ\psi , \B) \in \text{CPA}(\iota_1 , 
\Omega^{(n)} , \delta )$. Since $S(\omega_1 \circ 
(\Upsilon\circ\psi )) = S(\omega_2 \circ\psi )$, we conclude that
$$ cpe(\iota_1 , \omega_1 , \Omega^{(n)}, \delta ) \leq 
cpe(\iota_2 , \omega_2 , \Omega^{(n)}, \delta ). $$
Thus $hm_\sigma (\theta , \iota_1 , \omega_1) \leq hm_\sigma  
(\theta , \iota_2 , \omega_2 )$ and so, taking the supremum over
$\omega_1\in\mathfrak{E}(\sigma , \iota_1 )$, we obtain 
$hm_\sigma (\theta , \iota_1 ) \leq hm_\sigma (\theta , 
\iota_2 )$. The reverse inequality follows by 
symmetry.\hfill $\square$\medskip

\noindent{\bf Definition 4.3.} {\it In view of the above proposition
and the fact that $\A$ always admits a u.c.o.\ embedding 
into an injective $C^*$--algebra (consider, for instance, its 
universal representation), 
we can define $hm_\sigma (\theta )$ to be $hm_\sigma (\theta ,
\iota )$ for any u.c.o.\ embedding $\iota : \A \to \D$ into 
an injective $C^*$--algebra $\D$.}\medskip

\noindent{\it Remark.} If $\A$ is nuclear, we can dispense 
with state extension and
define $hm_\sigma^{\text{nuc}}(\theta )$, as does Choda in \cite{Ch} 
with different notation, 
by replacing the logarithm of the rank of the local algebra in
Voiculescu's topological definition \cite{V} with the entropy
of the induced local state. We can also adapt Voiculescu's AF 
definition \cite{V} in a similar way to define $hm_\sigma^{\text{AF}} 
(\theta )$ using the local characterization 
for AF algebras. Then
$$ hm_\sigma (\theta ) \leq hm_\sigma^{\text{nuc}} (\theta )
\leq hm_\sigma^{\text{AF}} (\theta ) , $$
with each inequality applying to the appropriate domain of
definition.\medskip

We show that, as for pressure, the local state
approximation entropy can
be computed by means of the larger class of contractive c.p. maps.
 If $\omega$ is a state
on $\D$, $\Omega\in Pf(\A)$, and $\delta>$, we
define
$$cpe_0(\iota,\omega,\Omega,\delta)=
\inf\{S(\omega\circ\psi):(\phi,\psi,\B)\in
\text{CPA}_0(\iota,\Omega,\delta)\},$$ 
and we define $hm_{\sigma}^0
(\theta,\iota,\omega,\Omega,\delta)$, etc, in the 
usual way.
\medskip

\noindent{\bf Lemma 4.4.} {\it Let $\B$ be a finite dimensional
$C^*$--algebra and
$a,b\in\B$ positive elements with $b$ invertible and $b\leq 0$, and 
suppose $a\leq\frac{1}{1+\epsilon}$
and $\|b^{-2}-1\|\leq\epsilon$,
for some $\epsilon<1$.
Let $f:{\Bbb R}^+\to{\Bbb R}$ be a concave function which is 
nonnegative-valued in $[0,1]$
and increasing in some interval $[0,\alpha]$ with 
$\epsilon\leq\frac{\alpha}{2}$. Let
$q$ be a
spectral projection
of $a$ such that $qaq\leq\alpha/2$. Then
$$\text{\rm Tr}_{\mathcal{B}}(f(b^{-1}ab^{-1}))\geq
\frac{1}{1+\epsilon}\text{\rm Tr}_{\mathcal{B}}(f(qaq)).$$}\medskip

\noindent{\it Proof.} Writing the spectral decompositions
$a=\sum_i\nu_jq_j$
and
$b^{-1}ab^{-1}=\sum_i\mu_ip_i$, with $p_i$ and $q_j$ minimal projections
of
$\B$,
we have
$$\mu_i=\text{Tr}_{\mathcal{B}}(b^{-1}ab^{-1}p_i)=
\sum_j\nu_j\text{Tr}_{\mathcal{B}}(b^{-2}p_i)
\frac{\text{Tr}_{\mathcal{B}}(b^{-1}q_jb^{-1}p_i)}{\text{Tr}_{\mathcal{B}}
(b^{-2}p_i)},$$
and so by the concavity of $f$
$$f(\mu_i)
\geq\sum_jf(
\nu_j\text{Tr}_{\mathcal{B}}(b^{-2}p_i)
)
\frac{\text{Tr}_{\mathcal{B}}(b^{-1}q_jb^{-1}p_i)}{\text{Tr}_{\mathcal{B}}
(b^{-2}p_i)}.$$
On the other hand, for all $j$ we have
$$0\leq\nu_j\text{Tr}_{\mathcal{B}}(b^{-2}p_i)\leq(1+\epsilon)\nu_j\leq1 , $$
and thus, since $f$ is nonnegative $[0,1]$,
$$f(\mu_i)
\geq\sum_{\{ j:\nu_j\leq{\alpha/2}\}}f(
\nu_j\text{Tr}_{\mathcal{B}}(b^{-2}p_i))
\frac{\text{Tr}_{\mathcal{B}}(b^{-1}q_jb^{-1}p_i)}{\text{Tr}_{\mathcal{B}}
(b^{-2}p_i)}.$$
If $\nu_j\leq\alpha/2$ then
$$|\nu_j\text{Tr}_{\mathcal{B}}(b^{-2}p_i)-\nu_j|=\nu_j
\text{Tr}_{\mathcal{B}}((b^{-2}-1)p_i)\leq\epsilon\nu_j\leq\alpha/2,$$
and so $\nu_j\leq\nu_j\text{Tr}_{\mathcal{B}}(b^{-2}p_i)\leq\alpha$. Since
$f$ is increasing in $[0,\alpha]$, we have
$$f(\mu_i)\geq\sum_{\{ j:\nu_j\leq\alpha/2\}}f(\nu_j)
\frac{\text{Tr}_{\mathcal{B}}(b^{-1}q_jb^{-1}p_i)}{\text{Tr}_{\mathcal{B}}
(b^{-2}p_i)}.$$
Therefore, summing up over $i$,
\begin{align*}
\text{Tr}_{\mathcal{B}}(f(b^{-1}ab^{-1})) &\geq \sum_i
\sum_{\{ j:\nu_j\leq\alpha/2\}}
f(\nu_j)\frac{\text{Tr}_{\mathcal{B}}(b^{-1}q_jb^{-1}p_i)}
{\text{Tr}_{\mathcal{B}}(b^{-2}p_i)} \\
&\geq \frac{1}{1+\epsilon}\text{Tr}_{\mathcal{B}}(b^{-1}f(qaq)b^{-1}) \\
&\geq \frac{1}{1+\epsilon}\text{Tr}_{\mathcal{B}}(f(qaq)).
\end{align*}$\hfill\square$\medskip

\noindent{\bf Proposition 4.5.} {\it We have
$hm_{\sigma}^0(\theta,\iota)=hm_\sigma(\theta,\iota)$.
}\medskip

\noindent{\it Proof.} Clearly $cpe_0(\iota,\omega,\Omega,\delta)\leq
cpe(\iota,\omega,\Omega,\delta)$. 
To show the reverse inequality, let $\Omega$ be a subset of the unit ball
of $\A$ containing $1$ and $\delta$ a positive number such that
$18\sqrt[4]{\delta}\leq\frac{1}{3}$, and, for
$n\in{\Bbb N}$,
let
$(\phi_n,\psi_n,\B_n)\in\text{CPA}_0(\iota,\Omega^{(n)},\delta)$
be such that
$S(\omega\circ\psi_n)<1+cpe_0(\iota,\omega,\Omega^{(n)},\delta)$.
We start following  the
same
procedure as in the proof of Prop.\ 2.3 to obtain
a corner $\B'_n$ of $\B_n$ obtained by cutting
with a 
(nonzero) spectral projection  $p$ of $\phi_n(1)$ such that
$b_1:=p\phi_n(1)\geq(1-\sqrt{\delta})p$.
We shall need the following estimates proven in
Prop.\ 2.3:
\begin{gather*}
\big\|{b_1}^{-\frac12}-p\big\|\leq\frac{\sqrt{\delta}}
{1-\sqrt{\delta}}, \\
\|\psi_n(pxp)-\psi_n(x)\|<16\sqrt[4]{\delta}\|x\|, \\
\|1-\psi_n(b_1)\|<2\sqrt{\delta}.
\end{gather*}
We first define
$\phi'_n:\A\to\B_n$ and $\psi'_n:\B_n\to\B(\H)$ by
$$\phi'_n(t):={b_1}^{-\frac12}\phi_n(t){b_1}^{-\frac12}+\gamma(t)(1-p),$$
where $\gamma$ is any state of $\A$, and
$$\psi'_n(t):=\psi_n\big(\big({b_1}^{\frac12}+1-p\big)t\big(
{b_1}^{\frac12}+1-p\big)\big).$$
Note that ${\phi'}_n$ is now unital, and
\begin{align*}
\|\psi'_n\phi'_n(t)-\iota(t)\| &\leq \|{\psi}_n(p\phi_n(t)p)-
\iota(t)\|+\|t\|\|{\psi}_n(1-p)\| \\
&\leq 32\sqrt[4]{\delta}\|t\|+\|\psi_n\circ\phi_n(t)-\iota(t)\|
\end{align*}
so that
$(\phi'_n,\psi'_n,\B_n)\in\text{CPA}_0(\Omega^{(n)},33\sqrt[4]{\delta})$.
We next fix $\psi'_n$ in order to obtain a unital map. Define
$\phi''_n:\A\to\B_n\oplus{\Bbb C}$ and 
$\psi''_n:\B_n\oplus{\Bbb C}\to\B(\H)$ by
\begin{gather*}
\phi''_n(t)={\phi'}_n(t)\oplus\gamma(t), \\
\psi''_n(t\oplus\lambda)=
{\psi'}_n(t)+\lambda(1-\psi_n(b_1+1-p)).
\end{gather*}
Note that, for $t\in\A$, 
$$\|\psi''_n\circ\phi''_n(t)-\iota(t)\|\leq
\|\psi'_n\circ{\phi'}_n(t)-\iota(t)\|+\|t\|\|1-\psi_n(b_1)\|+
\|t\|\|\psi_n(1-p)\| , $$
and so 
$(\phi''_n,\psi''_n,\B_n\oplus{\Bbb C})\in
\text{CPA}(\iota,\Omega^{(n)},51\sqrt[4]{\delta})$. 
We next estimate
$\limsup_n\frac{1}{n}S(\omega\circ\psi''_n)$.
If $A\in\B_n$ denotes the density matrix of
$\omega\circ{\psi}_n$,
$B:=\big({b_1}^{\frac12}+1-p\big)A\big({b_1}^{\frac12}+1-p\big)\oplus
\alpha_n$ is
the density matrix
of $\omega\circ\psi''_n$, where
$\alpha_n=\omega(1-\psi_n(b_1+1-p))$.
We claim that 
$$\limsup_n\frac{1}{n}
S(\omega\circ\psi''_n)=
\limsup_n\frac{1}{n}\text{Tr}_{\B_n}
\eta\big(q\big(b_1^{\frac12}+1-p\big)A\big(b_1^{\frac12}+1-p\big)q\big),$$ 
where $\eta(x)=-x\log x$ 
and $q$ is a spectral projection of 
$\big(b_1^{\frac12}+1-p\big)A\big(b_1^{\frac12}+1-p\big)$
 such that $q\big(b_1^{\frac12}+1-p\big)A\big(b_1^{\frac12}+1-p\big)
\leq\frac{1}{3}$
and $(1-q)\big(b_1^{\frac12}+1-p\big)A\big(b_1^{\frac12}+1-p\big)
\geq\frac{1}{3}(1-q)$.
To establish the claim, let
$\lambda_1,\dots,\lambda_N$ be the
eigenvalues of
$\big({b_1}^{\frac12}+1-p\big)A\big({b_1}^{\frac12}+1-p\big)$; then 
$(\lambda_1,\dots,\lambda_N,\alpha_n)$ are the
eigenvalues
of $B$.
Then 
\begin{align*}
\lefteqn{\text{Tr}_{\B_n}\big(\eta\big(q\big(b_1^{\frac12}+1-p\big)A\big(
b_1^{\frac12}+1-p\big)q\big)\big)}\hspace*{2.8cm} \\
&= -\sum_{\{i:\lambda_i \leq \frac{1}{3}\}}\lambda_i\log\lambda_i \\
&\leq -\sum_{\{i:\lambda_i\leq\frac{1}{3}\}}\lambda_i\log\lambda_i
-\sum_{\{i:\lambda_i>\frac{1}{3}\}}\lambda_i\log\lambda_i-\alpha_n
\log\alpha_n \\
&= S(\omega\circ\psi''_n) \\
&\leq\text{Tr}_{\B_n}\big(\eta\big(q\big(b_1^{\frac12}+1\big)A
\big(b_1^{\frac12}+1\big)q\big)\big)+\log 
3-\alpha_n\log(\alpha_n).
\end{align*}
Since $0\leq\alpha_n\leq1$ for all $n$, we have $0\leq-\alpha_n
\log\alpha_n<1$, and
therefore the claim follows by dividing by $n$ and taking the $\limsup_n$.
Applying the previous lemma to the matrices
\begin{gather*}
a=\frac{1}{1+\epsilon}\big(b_1^{\frac12}+1-p\big)A
\big(b_1^{\frac12}+1-p\big), \\
b=b_1^{\frac12}+1-p , 
\end{gather*} 
the function $f=\eta$, and  
$\epsilon=\frac{\sqrt{\delta}}{1-\sqrt{\delta}}$, we see that
$$(1+\epsilon)\text{Tr}_{\B_n}
\Big(\eta\Big(\frac{1}{1+\epsilon}A\Big)\Big)\geq\text{Tr}_{\B_n}
\Big(\eta\Big(\frac{1}{1+\epsilon}q\big(b_1^{\frac12}+1\big)A
\big(b_1^{\frac12}+1\big)q\Big),$$
and so
$\lim_{\delta\to0}\limsup_n\frac{1}{n}S(\omega\circ\psi''_n)\leq
\lim_{\delta\to0}\limsup_n
\frac{1}{n}S(\omega\circ\psi_n)$. We conclude that
$hm_{\sigma}(\theta,\iota,\omega,\Omega)\leq
hm_{\sigma}^0(\theta,\iota,\omega,\Omega)$, completing 
the proof.$\hfill\square$\medskip

We discuss some basic properties of the local state approximation entropy.
\medskip

\noindent{\bf  Proposition 4.6.} {\it
Let $\C$ be a unital $\theta$--invariant $C^*$--subalgebra of $\A$ and
$E:\A\to\C$ a conditional
expectation. If $\sigma$ is a $\theta$--invariant and $E$--invariant state
of $\A$,
$$hm_{\sigma\upharpoonright\C}(\theta\upharpoonright\C)\leq
hm_\sigma(\theta).$$}\medskip

\noindent{\it Proof} Let $\iota:\A\to\B(\H)$ be the universal
representation of $\A$,
$\omega\in\mathfrak{E}(\sigma\upharpoonright\C,\iota\upharpoonright\C)$,
$\Omega\in Pf(\C)$ and $\delta>0$. Extend $\iota\circ
E\circ\iota^{-1}:\iota(\A)\to\B(\H)$ to a
u.c.p. map $\tilde{E}:\B(\H)\to\B(\H)$, so $\iota\circ
E=\tilde{E}\circ\iota$ on $\A$. Thus, by $E$--invariance
of $\sigma$,
$\omega\circ\tilde{E}\in\mathfrak{E}(\sigma,\iota)$.
Now,
if
$(\phi,\psi,\B)\in
\text{CPA}(\iota\upharpoonright\C,\Omega^{(n)},\delta)$
then $(\phi,\tilde{E}\circ\psi,\B)\in
\text{CPA}(\iota\upharpoonright\C,\Omega^{(n)},\delta)$
as well, so
$$cpe(\iota\upharpoonright\C,\omega,\Omega,\delta)\leq
\text{inf}\{S(\omega\circ\tilde{E}\circ\psi): 
(\phi,\psi,\B)\in\text{CPA}(\iota\upharpoonright\C,
\Omega^{(n)},\delta)\}\leq$$
$$\text{inf}\{S(\omega\circ\tilde{E}\circ\psi):
(\phi,\psi,\B)\in\text{CPA}(\iota,
\Omega^{(n)},\delta)\}=cpe(\iota,\omega\circ\tilde{E},\Omega^{(n)},\delta)$$
which implies
$$hm_{\sigma\upharpoonright\C}(\theta\upharpoonright\C,\iota\upharpoonright\C,
\omega,\Omega,\delta)\leq hm_\sigma(\theta,\iota),$$
and the proof is complete.
$\hfill\square$
\medskip

\noindent{\bf Proposition 4.7.} {\it If $k\in\mathbb{N}$ then
$$ hm_\sigma (\theta^k ) = k\,hm_\sigma (\theta ) . $$}\medskip

\noindent{\it Proof.} Since $cpe(\iota , \omega , \Omega , \delta )$ 
is defined by taking an infimum over $\text{CPA}(\iota , \Omega , \delta
)$,
the second half of the proof of Prop.\ 1.3 in \cite{V} can be 
immediately adapted to our situation to establish the equality. 
Explicitly, we have 
$$ hm_\sigma (\theta^k , \iota , \omega , \Omega ,
\delta ) \leq k\,hm_\sigma (\theta , \iota , \omega , \Omega ,
\delta ) $$ 
because 
$$ \text{CPA}\Big(\iota , \bigcup_{j=0}^{n-1}\theta^{jk}(\Omega
), \delta \Big) \supset \text{CPA}\Big(\iota ,
\bigcup_{j=0}^{k(n-1)}\theta^j
(\Omega ), \delta \Big) $$ 
for all $n\in\mathbb{N}$, while the inequality 
$$ hm_\sigma (\theta^k , 
\iota , \omega , \bigcup_{j=0}^{k-1}\theta^j (\Omega ), \delta ) \geq 
k\,hm_\sigma (\theta , \iota , \omega , \Omega , \delta ) $$ 
follows from the observation that
$$ \text{CPA}\Big(\iota , \bigcup_{i=0}^{\lfloor
\frac{n}{k}\rfloor}\theta^{ik}
\Big(\bigcup_{j=0}^{k-1}\theta^j (\Omega )\Big), \delta \Big) 
\subset \text{CPA}\Big(\iota , \bigcup_{j=0}^{n-1}\theta^j (\Omega ),
\delta 
\Big) $$
for all $n\in\mathbb{N}$, whence the proposition follows by taking 
the supremum over all $\Omega\in Pf(A)$, $\delta > 0$, and 
$\omega\in\mathfrak{E}(\sigma , \iota )$ and applying Prop.\ 
4.2.\hfill $\square$\medskip

The proof of Prop.\ 3.4 can be adapted 
to establish the following Kol\-mo\-go\-rov--Sinai-type result.\medskip

\noindent{\bf Proposition 4.8.} {\it If $\iota : \A \to \D$ is a
u.c.o.\ embedding in an injective $C^*$--algebra $\D$, $\{ \Omega_\lambda 
\}_{\lambda\in I}$ is a net of elements of $Pf(\A)$ 
such that $\bigcup_{\lambda\in I}\bigcup_{j\in\mathbb{N}}\theta^j 
(\Omega_\iota )$ is total in $\A$ then
$$ hm_\sigma (\theta ) = \lim_\lambda \sup_{\omega\in\mathfrak{E}
(D, \sigma , \iota ) } hm_\sigma (\theta , \iota , \omega , 
\Omega_\lambda ) . $$}\medskip

We next compare $hm_\sigma (\theta )$ with the Sauvageot--Thouvenot
entropy $h_\sigma (\theta )$.  
We recall from \cite{ST} the notion of a stationary coupling
of $(\A,\theta,\sigma)$ with a unital commutative dynamical
system $(\C,\varsigma,\mu)$. Since
Sauvageot and Thouvenot treat the case in
which $\theta$ is an automorphism, they assume that
$\varsigma$ is an automorphism as well.
Since our $\theta$ is a u.c.p. map, 
we will assume, more naturally, that $\varsigma$ is 
a $^*$--homomorphism.
A {\it stationary coupling} is a 
$\theta\otimes\varsigma$--invariant state
$\lambda$
on $\A\otimes\C$ such that $\lambda(a\otimes 1)=\sigma(a)$,
$\lambda(1\otimes c)=\mu(c)$ for all $a\in\A$ and $c\in\C$.
The Sauvageot--Thouvenot entropy $h_\sigma(\theta)$ is the supremum 
of the quantities
$$h'(\mathcal{P}, \lambda ) := H_\mu\Big(\mathcal{P}\,\,\big|
\bigvee_{k=0}^{\infty}
\varsigma^k (\mathcal{P})\Big)-H_\mu(\mathcal{P})+\sum_{p\in
\mathcal{P}}\mu(p)S(\sigma, \sigma_p ) $$
where $\mathcal{P}$ ranges over all finite partitions of
$\C$ into projections and $\lambda$ over all stationary couplings,
with $\sigma_p$ denoting the state $x\mapsto\frac{1}{\mu(p)}\lambda
(x\otimes p)$ on $\mathcal{A}$ and $S(\cdot ,\cdot )$ denoting Araki's 
quantum relative entropy \cite{Ar}. 
Setting $\mathcal{P}^- = \bigvee_{k=1}^n \varsigma^{-k}\mathcal{P}$,
Sauvageot and Thouvenot show that
$h_\sigma(\theta)$ may be equivalently defined as the supremum
over the same set of $\mathcal{P}$ and $\lambda$ of the quantity
$$ h(\mathcal{P}, \lambda ) :=  H_\mu\Big(\mathcal{P}\,\,\big|
\bigvee_{k=0}^{\infty}\varsigma^k (\mathcal{P})\Big) - H_\lambda
(\mathcal{P}\, |\, \mathcal{A}\otimes\mathcal{P}^- ), $$ 
where $H_\lambda (\cdot \,|\, \cdot )$ denotes the conditional entropy
as defined in Sect.\ 2 of \cite{ST}, here with respect to the 
stationary coupling (also denoted for notational simplicity and
consistency with \cite{ST} by $\lambda$) 
of $(\mathcal{A}\otimes\mathcal{C}, \theta\otimes
\varsigma , \lambda )$ with $(\mathcal{C}, \varsigma , \mu )$
defined by composing $\lambda$ with $id_{\mathcal{A}}\otimes S$, where
$S : \mathcal{C}\otimes\mathcal{C}\to\mathcal{C}$ acts by restricting
functions to the diagonal. We shall find it convenient to use 
the following equivalent expression for $h(\mathcal{P}, \lambda )$
(cf.\ Prop.\ 3.3 of \cite{ST}) which involves the mutual entropy
of $\lambda$ with respect to a partition $\mathcal{Q}$ of $\mathcal{C}$
into projections as defined in \cite{ST} by
$$\varepsilon_\lambda(\A, \mathcal{Q})=
\sum_{q\in \mathcal{Q}}\mu(q)S(\sigma,\sigma_q).$$\medskip   

\noindent{\bf Lemma 4.9.} {\it If $\lambda$ is a stationary coupling
of $(\mathcal{A}, \theta , \sigma )$ with the unital commutative
system $(\mathcal{C}, \varsigma , \mu )$ and $\mathcal{P}$ is a
finite partition of $\mathcal{C}$ into projections then
$$ h(\mathcal{P}, \lambda ) = \lim_{n\to\infty} \frac{1}{n}
\varepsilon_\lambda \Big(\mathcal{A}, \bigvee_{k=1}^{n}\varsigma^k 
(\mathcal{P}) \Big) . $$}\medskip

\noindent{\it Proof.} As described above and in the paragraph
preceding Lemma 2.2 in \cite{ST}, the stationary coupling $\lambda$
defines a stationary coupling (denoted also by
$\lambda$) of $(\mathcal{A}\otimes\mathcal{C}, \theta\otimes
\varsigma , \lambda )$ with $(\mathcal{C}, \varsigma , \mu )$
via the map from $\mathcal{C}\otimes\mathcal{C}$ to $\mathcal{C}$ 
which restricts functions to the diagonal. We may also similarly 
define a stationary coupling (again denoted by $\lambda$) of 
$(\mathcal{A}\otimes\mathcal{C}\otimes
\mathcal{C}, \theta\otimes\varsigma\otimes\varsigma , \lambda )$
with $(\mathcal{C}, \varsigma , \mu )$ using the same map from
$\mathcal{C}\otimes\mathcal{C}$ to $\mathcal{C}$. 
For each $n\geq 2$ we then have by Lemma 2.2 of \cite{ST}
$$ H_\lambda\Big( \bigvee_{k=0}^{n-1} \varsigma^k\mathcal{P} \,\big| \,
\mathcal{A}\otimes\mathcal{P}^- \Big) = 
H_\lambda\Big( \bigvee_{k=0}^{n-2} \varsigma^k\mathcal{P} \,\big| \,
\mathcal{A}\otimes\mathcal{P}^- \Big) + 
H_\lambda\Big( \varsigma^{n-1}\mathcal{P} \,\big|\, \mathcal{A}\otimes
\mathcal{P}^-\otimes\bigvee_{k=0}^{n-2} \varsigma^k\mathcal{P}\Big) , $$
and since
\begin{align*}
H_\lambda\Big( \varsigma^{n-1}\mathcal{P} \,\big|\, \mathcal{A}\otimes
\mathcal{P}^-\otimes\bigvee_{k=0}^{n-2} \varsigma^k\mathcal{P}\Big)
&= H_\lambda\Big( \varsigma^{n-1}\mathcal{P} \,\big|\, \mathcal{A}\otimes
\bigvee_{k=-\infty}^{n-2} \varsigma^k\mathcal{P}\Big) \\
&= H_\lambda (\mathcal{P} \,|\, \mathcal{A}\otimes\mathcal{P}^- )
\end{align*} 
this leads inductively to 
$$ H_\lambda\Big( \bigvee_{k=0}^{n-1} \varsigma^k\mathcal{P} \,\big| \,
\mathcal{A}\otimes\mathcal{P}^- \Big) = 
(n+1) H_\lambda (\mathcal{P} \,|\, \mathcal{A}\otimes\mathcal{P}^- ). $$
Noting now that another application of Lemma 2.2 of \cite{ST} yields 
$$ H_\lambda\Big( \bigvee_{k=0}^{n-1} \varsigma^k\mathcal{P} \,\big| \,
\mathcal{A}\otimes\mathcal{P}^- \Big) =
H_\lambda\Big( \bigvee_{k=1}^{n-1} \varsigma^k\mathcal{P} \,\big| \,
\mathcal{A} \Big) - H_\lambda (\mathcal{P}\, | \, \mathcal{A}) , $$
we obtain
$$ H_\lambda\Big( \bigvee_{k=1}^n \varsigma^k\mathcal{P} \,\big| \,
\mathcal{A} \Big) = H_\lambda (\mathcal{P}\, | \, \mathcal{A}) +
(n+1) H_\lambda (\mathcal{P} \,|\, \mathcal{A}\otimes
\mathcal{P}^- ) . $$
Dividing by $n$ and taking the limit as $n$ tends to infinity yields
$$ H_\lambda (\mathcal{P} \,|\, \mathcal{A}\otimes\mathcal{P}^- ) =
\lim_{n\to\infty}\frac1n H_\lambda\Big( \bigvee_{k=1}^n 
\varsigma^k\mathcal{P} \,\big| \,\mathcal{A} \Big) . $$
Since
$$ H_\mu (\mathcal{P}\, | \, \mathcal{P}^- ) = 
\lim_{n\to\infty}\frac1n H_\mu\Big( \bigvee_{k=1}^n 
\varsigma^k\mathcal{P}\Big) $$
from the classical theory, we conclude that
\begin{align*}
h'(\mathcal{P}, \lambda ) &= H_\mu (\mathcal{P}\, | \,\mathcal{P}^- )
- H_\lambda (\mathcal{P}\, | \, \mathcal{A}\otimes\mathcal{P}^- ) \\
&= \lim_{n\to\infty} \frac1n \bigg[ H_\mu\Big( \bigvee_{k=1}^n 
\varsigma^k\mathcal{P}\Big) - H_\lambda\Big( \bigvee_{k=1}^n 
\varsigma^k\mathcal{P} \,\big| \,\mathcal{A} \Big) \bigg] \\
&= \lim_{n\to\infty} \frac1n \varepsilon_\lambda \Big(\mathcal{A}, 
\bigvee_{k=1}^{n}\varsigma^k (\mathcal{P}) \Big) .
\end{align*}
$\hfill\square$\medskip

\noindent{\bf Proposition 4.10.} {\it If $\theta:\A\to\A$
is a u.c.p. map and $\sigma$ 
is a $\theta$--invariant state on $\A$, then
$$ hm_\sigma (\theta ) \geq h_\sigma (\theta ). $$}\medskip

\noindent{\it Proof.} Let $\iota : \A \to \D$ be a u.c.o.\ embedding 
into an injective $C^*$--algebra $\D$. Suppose 
$\lambda$ is a stationary coupling of $(\A, \theta , \sigma )$ 
with $(\C, \varsigma , \mu )$, with $\mu$ assumed to be 
faithful. Extend 
the state $\lambda\circ(\iota^{-1}\otimes id)$ on $\iota (\A)\otimes\C$ 
to a state $\tilde{\lambda}$ on $\D \otimes \C$. Suppose 
$\mathcal{P}$ is a finite partition of projections in $\C$. 
For each $n\in\mathbb{N}$ and $p\in\bigvee_{k=1}^n\varsigma^k
(\mathcal{P})$, let $\sigma_p$ be the state on $\A$ defined
by $x \mapsto \mu (p)^{-1} \lambda (x\otimes p)$ and $\omega_p$  
the state on $\D$ defined by $y \mapsto \mu (p)^{-1} 
\tilde{\lambda} (y\otimes p)$. Note that $\omega_p$ extends the state
$\sigma_p\circ\iota^{-1}$ on $\iota (\A)$. Let $\omega$
be the state on $\D$ given by the convex combination 
$\sum_{p\in\varsigma (\mathcal{P})}\mu (p)\omega_p$ (which is equal to 
$\sum_{p\in\bigvee_{k=1}^n\varsigma^k (\mathcal{P})}
\mu (p)\omega_p$ for any $n\in\mathbb{N}$).

For every $n\in\mathbb{N}$, 
$\Omega\in Pf(\A)$, and $\delta > 0$ choose
$$ \left(\phi_{(\Omega , \delta ), n}, \psi_{(\Omega ,
\delta ), n}, \B_{(\Omega , \delta ) , n}\right) \in \text{CPA}(\iota , 
\Omega^{(n)}, \delta ) $$
such that 
$$ hm_\sigma (\theta , \iota , \omega , \Omega , \delta ) =
\limsup_{n\to\infty} \frac{1}{n} S\left(\omega\circ\psi_{(\Omega, 
\delta ), n}\right).$$
Set $\Gamma = Pf(\A) \times \mathbb{R}_{>0}$. For each $n\in
\mathbb{N}$, $\{ \psi_{\gamma , n}\circ\phi_{\gamma , n} \}_{\gamma
\in\Gamma }$ is a net converging pointwise in norm to $\iota$, 
so that $\{ \omega\circ\psi_{\gamma , n}\circ\phi_{\gamma , n}
\}_{\gamma\in\Gamma }$ converges weak$^*$ to $\sigma$ and,
for all $p\in\bigvee_{k=1}^n-\varsigma^k
(\mathcal{P})$, $\{ \omega_p\circ\psi_{\gamma , n}\circ
\phi_{\gamma , n} \}_{\gamma\in\Gamma }$ converges weak$^*$ to 
$\sigma_p$. The weak$^*$ lower semicontinuity 
of the relative entropy $S(\cdot , \cdot )$ and the
weak$^*$ compactness of the state space of $\A$ then yields 
a $\gamma_0 = (\Omega_0 , \delta_0 )\in \Gamma$ such 
that, for all $n\in\mathbb{N}$ and $p\in\bigvee_{k=1}^{n}\varsigma^k
(\mathcal{P})$,
$$ S(\sigma , \sigma_p ) < S\left(\omega\circ\psi_{\gamma_0 , n}\circ
\phi_{\gamma_0 , n}, \omega_p\circ\psi_{\gamma_0 , n}\circ
\phi_{\gamma_0 , n}\right) + 1.$$
Since
$$ S\left(\omega\circ\psi_{\gamma_0 , n}\circ\phi_{\gamma_0 , n}, 
\omega_p\circ\psi_{\gamma_0 , n}\circ\phi_{\gamma_0 , n}\right) \leq
S(\omega\circ\psi_{\gamma_0 , n}, \omega_p \circ\psi_{\gamma_0 , 
n}) $$
by the monotonicity of $S(\cdot , \cdot )$, we therefore have
\begin{align*}
h(\mathcal{P}, \lambda ) &= \lim_{n\to\infty} \frac{1}{n}
\varepsilon_\lambda \Big(\mathcal{A}, \bigvee_{k=1}^n\varsigma^k 
(\mathcal{P}) \Big) \\
&= \lim_{n\to\infty} \frac{1}{n} \sum_{p\in \bigvee_{k=1}^n 
\varsigma^k (\mathcal{P})} \mu (p) S(\sigma , \sigma_p ) \\ 
&\leq \limsup_{n\to\infty} \frac{1}{n} \sum_{p\in \bigvee_{k=1}^n
\varsigma^k (\mathcal{P})} \mu (p) S(\omega\circ\psi_{\gamma_0 , n}, 
\omega_p\circ\psi_{\gamma_0 , n}) \\
&\leq \limsup_{n\to\infty} \frac{1}{n} S(\omega\circ\psi_{\gamma_0 , 
n}) \\
&= hm_\sigma (\theta , \iota , \omega , \Omega_0 , \delta_0 ).
\end{align*}

Taking the supremum over all stationary couplings $\lambda$ and 
finite partitions $\mathcal{P}$, we obtain the 
proposition.\hfill $\square$\medskip

Next we show that the local state approximation entropy agrees
with the Kolmogorov--Sinai entropy in the commutative case. 
For an open cover $\mathcal{U}$ of a topological space $X$ we
denote by $\mathfrak{S}(\mathcal{U})$ the set of all $x\in X$ which 
are contained in only one member of $\mathcal{U}$.\medskip

\noindent{\bf Lemma 4.11.} {\it Let $\mu$ be a measure on a compact 
Hausdorff space $X$. If $\mathcal{U}=\{ U_1 , \dots , U_m \}$ is a 
finite open cover of X, then for every $\epsilon > 0$ there
is a open refinement $\mathcal{V}=\{ V_1 , \dots V_m \}$ of 
$\mathcal{U}$ such that there are 
closed sets $H_i \subset V_i$ for $i=1, \dots , m$ such that
$\bigcup_{i=1}^{m} H_i \subset\mathfrak{S}(\mathcal{V})$ and
$\mu (X\setminus\bigcup_{i=1}^{m} H_i ) < \epsilon$.}\medskip

\noindent{\it Proof.} Let $\mathcal{U}=\{ U_1 ,\dots , U_m \}$ be 
a finite open cover of $X$ and $\epsilon > 0$. Set $V_1 = U_1$. 
Let $G_1 \subset U_1$ be a closed set such that $\mu (V_1 \setminus
G_1 ) < \frac{\epsilon}{m^2}$ and set $V_2 = U_2 \cap (X\setminus
G_1 )$. We continue inductively for $k=3,\dots ,m$ so that at the $k$th 
stage we choose a closed set $G_{k-1}\subset V_{k-1}$ such that $\mu 
(V_{k-1}\setminus G_{k-1}) < \frac{\epsilon}{m^2}$ and set
$V_k = U_k \cap (X\setminus \bigcup_{j=1}^{k-1} G_j )$.
 
Put $H_1 = G_1$ and, for $i=2, \dots , m$, $H_i = G_i \setminus
(V_1 \cup \cdots \cup V_{i-1})$. Then $\bigcup_{i=1}^{m} H_i \subset
\mathfrak{S}(\mathcal{V})$, and since $G_1 , \dots , G_m$ are 
pairwise disjoint so are $H_1, \dots , H_m$. Furthermore, for each 
$i=1, \dots , m$ we have
$$ \mu (V_i \setminus H_i ) \leq \sum_{j=1}^{i}\mu (V_j \setminus G_j )
< i\frac{\epsilon}{m^2} \leq \frac{\epsilon}{m}, $$
so that
$$ \mu \bigg( X\setminus\bigcup_{1\leq k \leq m}H_k \bigg) \leq
\mu \bigg( \bigcup_{1\leq k \leq m}(V_k \setminus H_k )
\bigg) \leq \sum_{k=1}^m \mu ( V_k \setminus H_k ) < m
\frac{\epsilon}{m} = \epsilon , $$
as required.\hfill $\square$\medskip

\noindent{\bf Proposition 4.12.} {\it Let $T:X\to X$ be a homeomorphism
of a compact metric space and $\mu$ a $T$--invariant measure on $X$.
If $\theta_T$ is the automorphism of $C(X)$ induced by $T$ and
$\sigma$ denotes the state on $C(X)$ defined by $\mu$, then
$$ h_\mu (T) = hm_\sigma (\theta_T ), $$
where $h_\mu (T)$ is the Kolmogorov-Sinai entropy of $T$.}\medskip

\noindent{\it Proof.} Since the local state approximation entropy
is bounded below by the Sauvageot--Thouvenot entropy (Prop.\ 4.10)
and the latter agrees with the
Kol\-mo\-go\-rov-Sinai entropy in the commutative case, we have
$h_\mu (T) \leq hm_\sigma (\theta_T )$.

To establish the reverse inequality, let $\iota : C(X) \to C(X)^{**}$
be the natural embedding. Note that, since $C(X)$ is nuclear, 
$C(X)^{**}$ is injective \cite{EL}. Suppose $\omega\in\mathfrak{E}(\sigma , 
\iota )$, $\Omega\in Pf(C(X))$ and $\delta > 0$. 
Let $\mathcal{U}$ be an open 
cover of $X$ such that if $U\in\mathcal{U}$ and $x,y\in U$ then 
$|f(x)-f(y)|\leq\delta$ for all $f\in\Omega$. Writing 
$\mathcal{U}=\{ U_1 , \dots , U_r \}$, let $\mathcal{P}$ be the Borel 
partition $\{ U_i \setminus\bigcup_{j=1}^{i-1}U_j : 1\leq i \leq r 
\}$ refining $\mathcal{U}$. Fix $n\in\mathbb{N}$. Note that if $U\in
\bigvee_{j=0}^{n-1}T^j (\mathcal{U})$ and $x,y,\in U$ then  
$|f(x)-f(y)|\leq\delta$
for all $f\in\Omega^{(n)}$. Let $\epsilon > 0$ be small enough so
that if $0\leq a,b \leq 1$ and $|a-b|<\epsilon$ then 
$| a\log a - b\log b | < r^{-n}$. By the lemma there
is a refinement $\mathcal{V}=\{ V_1 , \dots V_m \}$ of 
$\bigvee_{j=0}^{n-1}T^j (\mathcal{U})$ such that there are pairwise 
disjoint closed sets $H_i \subset V_i$ for $i=1, \dots , m$ such that
$\bigcup_{i=1}^{m} H_i \subset\mathfrak{S}(\mathcal{U})$ and
$\mu (X\setminus\bigcup_{i=1}^{m} H_i ) < \epsilon$. Let 
$\Xi = \{ \chi_1 , \dots , \chi_m \}$ be a 
partition of unity subordinate to $\mathcal{V}$ 
and $X_n = \{ x_1 ,\dots , x_m \}$ a finite subset of $X$ such 
that $x_i \in H_i$ for each $i=1, \dots , m$. Defining $\phi_n : 
C(X) \to C(X_n )$ by $f\mapsto f\restriction X_n$ and $\psi_n :
C(X_n ) \to C(X)^{**}$ by $g\mapsto\sum_{1\leq i \leq m} g(x_i ) 
\iota\circ\chi_i $. Since $\Xi$ is subordinate to 
$\bigvee_{j=0}^{n-1}T^j(\mathcal{U})$, for all $f\in\Omega^{(n)}$ 
and $x\in X$ we have
\begin{align*}
| ((\iota^{-1}\circ\psi_n \circ\phi_n )(f))(x)-f(x)| &\leq 
\sum_{1\leq i \leq m}\chi_i (x) | f(x_i) - f(x) | \\
&= \sum_{\{ i\, : \,x\in V_i \}}\chi_i (x) | f(x_i) - f(x) | \\
&< \delta 
\end{align*}
and hence
$$ | (\psi_n \circ \phi_n)(f) - \iota (f) | < \delta , $$
so that $( \phi_n , \psi_n , C(X_n ))\in \text{CPA}(\iota , 
\Omega^{(n)}, \delta )$. 

Now for each $i=1, \dots , m$, we have
$$ \left| \sigma(\chi_i ) - \mu (H_i ) \right| \leq \mu (V_i 
\setminus H_i ) < \epsilon $$
since $\mu (X\setminus\bigcup_{1\leq i\leq m}H_i) < \epsilon$ and $V_i$
does not intersect $H_j$ for $j=1, \dots , m, j\neq i$. Thus,
since $m\leq r^{n}$, our choice of $\epsilon$ yields
\begin{align*}
S(\omega\circ\psi_n ) &= -\sum_{i=1}^{m} \sigma (\chi_i ) \log
\sigma (\chi_i ) \\
&\leq -\sum_{i=1}^{m} \mu (H_i ) \log \mu (H_i ) + 1 . 
\end{align*}
Setting $K = \bigcup_{i=1}^{m} H_i$ we have $| \mu (P\cap K ) - \mu 
(P) | \leq  \mu (X\setminus K )  < \epsilon$ for all $P\in
\bigvee_{i=1}^{n-1}T^i(\mathcal{P})$. Also note that, since each $P\in
\bigvee_{i=1}^{n-1}T^i(\mathcal{P})$ intersects at most
one of $H_1,\dots , H_m$, the partition $\{ P\cap K : P\in
\bigvee_{i=1}^{n-1}T^i(\mathcal{P}) \}$ of $K$ refines 
$\{ H_i : 1\leq i \leq m \}$, and so we infer  
\begin{align*}
-\sum_{i=1}^{m}\mu (H_i ) \log\mu (H_i ) &\leq - \sum_{P\in
\bigvee_{i=1}^{n-1}T^i(\mathcal{P})}\mu (P\cap K ) \log\mu (P\cap K ) \\
&\leq - \sum_{P\in\bigvee_{i=1}^{n-1}T^i(\mathcal{P})}\mu (P) \log\mu
(P) +1 .
\end{align*}
Combining the above two estimates we obtain $S(\omega\circ\psi_n ) 
\leq H_\mu (\bigvee_{j=0}^{n-1}T^j (\mathcal{P}), \theta) + 2$, so that
$$ cpe(\iota , \omega , \Omega^{(n)} , \delta ) \leq H_\mu 
\Big(\bigvee_{j=0}^{n-1}T^j (\mathcal{P}), \theta \Big) + 2. $$
Dividing by $n$ and taking the $\limsup$ yields
$hm_\sigma (\theta , \iota , \omega , \Omega , \delta ) \leq
H_\mu (\mathcal{P}, \theta )$. Taking the supremum over all
$\delta > 0$, $\Omega\in Pf(C(X))$, and $\omega\in\mathfrak{E}(\sigma , 
\iota )$, we conclude that $hm_\sigma (\theta ) = hm_\sigma 
(\theta , \iota ) \leq h_\mu (\theta )$.\hfill $\square$\medskip
 
\noindent{\bf Proposition 4.13.} {\it (Concavity)\, If $\sum_{i=1}^k 
\lambda_i \sigma_i$ is a convex
combination of $\theta$--invariant states $\sigma_i$ on $\A$ then 
$$ \sum \lambda_i hm_{\sigma_i} (\theta ) \leq hm_{\sum \lambda_i 
\sigma_i} (\theta ). $$}\medskip

\noindent{\it Proof.} Let $\iota : \A \to \D$ be an embedding
into an injective $C^*$--algebra. Set $\sigma = \sum_{i=1}^k 
\lambda_i \sigma_i$. For each $i=1,\dots , k$ let $\omega_i \in 
\mathfrak{E}(\sigma_i , \iota )$. 
Then the state $\omega$ defined by $\sum \lambda_i \omega_i$ lies in 
$\mathfrak{E}(\sigma , \iota )$, and
$$ \frac{1}{n} \sum \lambda_i \, \log cpe\left( \iota , \omega_i , 
\Omega^{(n)}, \delta \right) \leq \frac{1}{n}  
\log cpe\left( \iota , \sum\lambda_i \omega_i , \Omega^{(n)}, 
\delta \right) $$
by the concavity of $S(\cdot )$ on state spaces of finite-dimensional
$C^*$--algebras. Therefore 
$$\sum\lambda_i hm_{\sigma_i} (\theta , \iota , \omega_i , 
\Omega , \delta ) \leq hm_{\sum\lambda_i \sigma_i} \left( \theta , 
\iota , \sum\lambda_i \omega_i , \Omega , \delta \right)$$
and hence
$\sum\lambda_i hm_{\sigma_i} (\theta , \iota , \omega_i ) 
\leq hm_{\sum\lambda_i \sigma_i} ( \theta , \iota , \sum\lambda_i 
\omega_i )$. 

Taking the supremum successively for
each $i=1,\dots ,k$ over $\omega_i \in \mathfrak{E}(\sigma_i , 
\iota )$ yields
$$ \sum\lambda_i hm_{\sigma_i}(\theta , \iota ) \leq
hm_{\sum\lambda_i \sigma_i}(\theta , \iota ), $$
establishing the proposition.\hfill $\square$\medskip

The following proposition establishes a variational 
inequality bounding the free energy in a given state by the 
approximation pressure.\medskip

\noindent{\bf Proposition 4.14.} {\it If $a\in \A_{sa}$ and $\sigma$ 
is a $\theta$--invariant state on $\A$, then
$$ P_\theta (a) \geq hm_{\sigma} (\theta ) + \sigma 
(a). $$}\medskip 

\noindent{\it Proof.} Let $\pi : \A \to \mathcal{B}(\mathcal{H})$ 
be a faithful representation of $\A$. Then $P_\theta (a) =
P_\theta (a, \pi )$ and $hm_\sigma (\theta ) = hm_\sigma  
(\theta , \pi )$. 

Suppose $\omega\in\mathfrak{E}(\sigma , 
\pi )$. Let $\Omega$ be a set in $Pf(\A)$ containing $a$, and
suppose $\delta > 0$ and $n\in\mathbb{N}$. 
If $(\phi , \psi , \B) \in \text{CPA}(\pi , 
\Omega^{(n)} , \delta )$ then
\begin{align*}
\text{log Tr}\, e^{\phi (a^{(n)})} &\geq S(\omega\circ\psi ) + 
(\omega\circ\psi ) \big(\phi \big(a^{(n)}\big) \big) \\
&\geq S(\omega\circ\psi ) + n(\omega\circ\pi ) (a) - n\delta \\
&= S(\omega\circ\psi ) + n\sigma (a) - n\delta , 
\end{align*} 
so that
$$ \frac{1}{n}\log Z_{\theta , n} (a , \pi , \Omega , \delta ) \geq  
\frac{1}{n}cpe(\pi , \omega , \Omega , \delta ) + 
\sigma (a) - \delta. $$ 
Hence $P_\theta (a, \pi , \Omega , \delta ) \geq hm_\sigma
(\theta , \pi , \omega , \Omega , \delta ) + \sigma (a) - \delta$ and 
therefore $P_\theta (a, \pi ) \geq hm_\sigma (\theta , \pi , \omega ) + 
\sigma (a)$. 

Taking the supremum over $\omega\in\mathfrak{E}(\sigma , \pi)$ yields 
$P_\theta (a, \pi ) \geq hm_\sigma (\theta , \pi ) + \sigma (a)$, 
thus establishing the proposition.\hfill $\square$\medskip

\noindent In view of Prop.\ 4.14 and the fact that
the Sauvageot--Thouvenot entropy majorizes the CNT entropy
\cite{ST}, we
immediately obtain the following corollary.\medskip

\noindent{\bf Corollary 4.15.} {\it If $\theta\in Aut(\A)$ and $\sigma$ 
is a $\theta$--invariant state on $\A$, then the variational
inequality of the previous proposition also holds when  
$hm_\sigma (\theta )$ is replaced by the Sauvageot--Thouvenot 
or CNT entropy.}\medskip

 Cor. 4.15 leads us to introduce the notion of equilibrium state.
\medskip

\noindent{\bf Definition 4.16.} {\it Let $a$ be a self-adjoint element
of $\A$, and $\theta$ a u.c.p. map of $\A$. An equilibrium state for  $(\A, \theta, a)$
is a $\theta$--invariant state $\sigma$ such that
$h_\sigma(\theta)+\sigma(a)=P_\theta(a)$.}\medskip

To round out this section we discuss some properties of the 
Sauvageot--Thouvenot entropy which we will need in Sect. 6.
The following has been noted by Neshveyev and St\o rmer
(cf.\ Lemma 3.5 in \cite{NS}).
\medskip


\noindent{\bf Proposition 4.17.} \cite{NS} {\it Let $\A$ be a unital
$C^*$--algebra
endowed with a u.c.p. map $\theta$, and let $\B\subset\A$ be a
$\theta$--invariant unital $C^*$--subalgebra. Then, given 
$\epsilon>0$, any invariant state $\sigma$ on $\B$ 
extends to an invariant state $\tilde{\sigma}$ on $\A$
in such a way that
$$h_{\tilde{\sigma}}(\theta)>h_{\sigma}(\theta\upharpoonright\B)-
\epsilon.$$}\medskip

\noindent To establish concavity we need the Donald identity 
(cf.\ Prop.\ 5.23(v) in \cite{OP}), namely, if $\eta$ is a state and 
$\sigma=\sum_i\alpha_i\sigma_i$ is a convex combination of states then
$$\sum_i\alpha_i S(\eta,\sigma_i)=S(\eta,\sigma)+\sum_i\alpha_i
S(\sigma,\sigma_i).$$\medskip

\noindent{\bf Proposition 4.18.} {\it If $\A$ is a unital $C^*$--algebra,
$\theta:\A\to\A$ is a u.c.p. map, and $\alpha\sigma_1 + 
\beta\sigma_2$ is a convex combination of $\theta$--invariant 
states on $\A$, then
$$h_{\alpha\sigma_1+\beta\sigma_2}(\theta)\geq\alpha
h_{\sigma_1}(\theta)+\beta
h_{\sigma_2}(\theta)-(\alpha\log\alpha+\beta\log\beta).$$}\medskip%

\noindent{\it Proof.} Set $\sigma = \alpha\sigma_1 + 
\beta\sigma_2$ and, for $i=1,2$, let $(\lambda_i,{\cal P}_i)$ be a
stationary coupling of $(\A,\theta,\sigma_i)$ with $(\C_i, \varsigma_i
,\mu_i)$ which is optimal within $\epsilon$ with respect to 
$h'(\cdot, \cdot )$. We construct a stationary pairing $\lambda$
of  $(\A,\theta,\sigma)$ with
a commutative system $(\C,\varsigma,\mu)$ as follows.
Set $\C=\C_1\oplus\C_2$, $\varsigma=\varsigma_1\oplus \varsigma_2$, 
$\mu=\alpha\mu_1+\beta\mu_2$,
$\lambda=\alpha\lambda_1+\beta\lambda_2$,
$\mathcal{P}=\mathcal{P}_1\oplus\mathcal{P}_2$.
Then the entropy of this model is
\begin{gather*}
\alpha\Big(H_\mu\Big(\mathcal{P}_1\,\,\big|
\bigvee_{k=0}^{\infty}
\varsigma^k (\mathcal{P}_1)\Big)-H_\mu(\mathcal{P}_1)\Big)+
\beta\Big(H_\mu\Big(\mathcal{P}_2\,\,\big|
\bigvee_{k=0}^{\infty}
\varsigma^k (\mathcal{P}_2)\Big)-H_\mu(\mathcal{P}_2)\Big)\\
-\alpha\log\alpha-\beta\log\beta+
\alpha\sum_{p\in\mathcal{P}_1}\mu_1(p)S(\sigma,
\sigma_p)+\beta\sum_{p\in\mathcal{P}_2}\mu_2(p)S(\sigma, \sigma_p).
\end{gather*}
By the Donald identity, if $\eta$ is any state on $\A$,
\begin{align*}
\lefteqn{\alpha\sum_{p\in\mathcal{P}_1}\mu_1(p)S(\sigma,
\sigma_p)
+\beta\sum_{p\in\mathcal{P}_2}\mu_2(p)S(\sigma,
\sigma_p)}\hspace*{2cm}\\
&= -S(\eta,\sigma)+
\alpha\sum_{p\in\mathcal{P}_1}\mu_1(p)S(\eta,(\sigma_1)_p) \\
&\hspace*{4.2cm}\ +\beta\sum_{p\in\mathcal{P}_2}\mu_2(p)S(\eta,
(\sigma_2)_p) \\
&\geq-\alpha S(\eta,\sigma_1)-\beta S(\eta,\sigma_2)+
\alpha\sum_{p\in\mathcal{P}_1}\mu_1(p)S(\eta,(\sigma_1)_p) \\
&\hspace*{5cm}\ +\beta\sum_{p\in\mathcal{P}_2}\mu_2(p)S(\eta,(\sigma_2)_p)
\end{align*}
by the joint convexity of the relative entropy.
Again applying the Donald identity, the last term coincides with
$$\alpha \sum_{p\in\mathcal{P}_1}\mu_1(p)S(\sigma_1,(\sigma_1)_p)+
\beta \sum_{p\in\mathcal{P}_2}\mu_2(p)S(\sigma_2,(\sigma_2)_p).$$
$\hfill\square$\medskip

\noindent{\bf Lemma 4.19.} {\it Let $\theta:\A\to\A$ be a u.c.p. map,
and assume that $\sigma$ is $\theta^r$--invariant for some
$r\in{\Bbb N}$. Then, for all $j=1,\dots,r-1$,
$$h_{\sigma\circ\theta^j}(\theta^r)\geq h_{\sigma}(\theta^r).$$}\medskip

\noindent{\it Proof.} If $\lambda$ is an optimal stationary pairing of
$(\A,\sigma,\theta^r)$ with $(\C,\mu,\varsigma)$ and $\mathcal{P}$
is a finite partition of projections in $\C$ then
$\lambda_j:=(\lambda\circ\theta^j)\otimes\varsigma$
defines a stationary pairing of $(\A,\sigma\theta^j,\theta^r)$ with
$(\C,\mu,\varsigma)$.
The monotonicity of the quantum relative entropy yields
$$S(\sigma\circ\theta^j,\sigma_p\circ\theta^j)\geq S(\sigma\circ
\theta^j\circ\theta^{r-j},\sigma_p\circ\theta^j\circ\theta^{r-j})=
S(\sigma,\sigma_p),$$
establishing the inequality.\hfill$\square$\medskip

\section{Cuntz--Krieger and Crossed Product Algebras}

In this section we examine pressure in Cuntz--Krieger and crossed
product algebras, exercising in different directions their respective
structures as Pimsner algebras (see \cite{P} and Sect. 6).\bigskip

{\noindent\it 1.\ Cuntz--Krieger algebras}\medskip

{\noindent In} classical ergodic theory the variational principle was 
first proved for lattice systems by Ruelle 
\cite{R1, R2} (see also \cite{Rbook}). If we assume for 
simplicity that the system is one-dimensional,
then it is isomorphic to a  subshift of finite type
\cite{Rbook}. 
The 
partition function corresponding to the classical pressure then takes 
the simple form
\begin{gather*}
Z_n(f)=\sum_C e^{\max\{ \sum_{j=0}^{n-1} f\circ T^j(x): x\in C\}}
\end{gather*}
where the sum is taken over all cylinders $C$ of the subshift obtained
by fixing the first $n$ coordinates.
Inspired by this, we consider the Cuntz--Krieger algebra $\O_A$
(a ``noncommutative subshift of finite type'') associated
to a matrix $A\in M_d(\{0,1\})$ with no row or column identically zero,
introduced by Cuntz and Krieger \cite{CK}.
$\O_A$ is generated by partial isometries $s_1,\dots, s_d$ satisfying
\begin{gather*}
\sum_i s_is_i^*=I\\
s_i^*s_i=\sum_j A_{i j}s_js_j^*.
\end{gather*}
Let $\Lambda_A$ be the one-sided Markov subshift defined by $A$:
\begin{gather*}
\Lambda_A:=\{(x_i)_i\in\{1,\dots,d\}^{\Bbb N}: A_{x_i x_{i+1}}=1,i\in{\Bbb
N}\}.
\end{gather*}
The commutative algebra $\C(\Lambda_A)$ of complex-valued continuous
functions on $\Lambda_A$ sits naturally inside $\O_A$ as the
$C^*$--subalgebra generated by the range projections of the 
iterated products $s_{i_1}\cdots s_{i_r}$. The shift epimorphism
$T: (x_1,x_2,\dots)\in\Lambda_A\to (x_2,x_3,\dots)\in\Lambda_A$
corresponds to the
restriction to $\C(\Lambda_A)$
of the u.c.p. map 
$$\theta: t\in\O_A\to\sum s_its_i^*\in\O_A.$$
(It suffices to check this on the set $\{s_{i_i}\cdots
s_{i_r}(s_{i_1}\cdots s_{i_r})^*\}$, which is total in $\C(\Lambda_A)$).
Let $\alpha=(i_1,\dots,i_r)$ be a finite word of lenth $r=:|\alpha|$
occurring in some element
of $\Lambda_A$; then $s_\alpha:=s_{i_1}\cdots s_{i_r}\neq 0$.
Let $[\alpha]$ denote the cylinder set of
$\Lambda_A$ given by
$$[\alpha]:=\{(x_i)_i\in\Lambda_A: x_1=i_1,\dots,x_r=i_r\}.$$
\medskip

\noindent{\bf Lemma 5.1.} {\it If $f$ is a continuous function on
$\Lambda_A$, then
\begin{gather*}
s_\alpha^*f{s_\beta}=0,\quad |\alpha|=|\beta|,\, \alpha\neq\beta,\\
s_\alpha^*f{s_\alpha}(x)=A_{i_r x_1}f(\alpha x),
\end{gather*}
and thus, for $f\geq0$,
\begin{gather*}
s_\alpha^*fs_\alpha\leq\max\{f(x): \ x\in[\alpha]\} I.
\end{gather*}
}\medskip

These computations turn out to be useful in proving the following
result.
\medskip

\noindent{\bf Theorem 5.2.} {\it For any self-adjoint element  
$f$ of $\O_A$ belonging to the subalgebra $\C(\Lambda_A)$,
the pressure of $f$ with respect to $\theta$ equals
the classical  pressure of $f$ with respect to the shift $T$,
i.e.,
$$P_\theta(f)=p_T(f).$$}\medskip

\noindent{\it Proof.} By the additivity of both $P_\theta$ and
$p_T$ under the addition of scalars (Prop.\ 3.1(b)), we can 
assume $f\geq0$. 
Furthermore, by monotonicity (Prop.\ 3.3)
and Prop.\ 3.6,
\begin{gather*}
p_T(f)=P_{\theta\upharpoonright\C(\Lambda_A)}\leq P_\theta(f).
\end{gather*}
We are thus left to show that $P_\theta(f)\leq p_T(f)$.
By Prop.\ 2.3, we can use
the non-unital
exact definition of pressure. 
We generalize the arguments by Boca and Goldstein \cite{BG} for the 
Voiculescu--Brown entropy.
Consider a finite  set of the form $\Omega(n_0)=\{s_\alpha p_is_\beta^*,
|\beta|\leq|\alpha|\leq n_0\}$
where $p_i:=s_i{s_i}^*$ and $\alpha$ and $\beta$ are finite words
appearing 
in $\Lambda_A$. Consider  the contractive c.p.\ maps
$\rho_m:\O_A\to M_{\vartheta_m}(\O_A)$,
$t\to(s_\mu^*ts_\nu)$, where
$\vartheta_m$ is the number of blocks occurring in elements of $\Lambda_A$
of length $m$. By Lemma 2 in \cite{BG}, for $m\geq n+n_0$,
$j=0,\dots,n-1$, and
$t=s_\alpha p_is_\beta^*\in\Omega(n_0)$, if $|\beta|<|\alpha|$ then
$$\rho_m(\theta^j(t))=\sum_{|\mu|=|\alpha|-|\beta|} x(\mu)\otimes
s_\mu$$
while if $|\alpha|=|\beta|$ then
$$\rho_m(\theta^j(t))=\sum_{j=1}^d x(j)\otimes s_j^*{s_j},$$
where $x(\mu)$ and $x(j)$ are partial isometries of $M_{\vartheta_m}({\Bbb
C})$ depending also on
$i$, $\alpha$, and $\beta$.
Given $\delta>0$ consider $(\phi_0,\psi_0, M_{m_0}({\Bbb C}))$ such that
$$\|(\psi_0\circ\phi_0)(s_r^*s_r)-s_r^*s_r\|+\|(\psi_0\circ\phi_0)
(s_\gamma)-s_\gamma\| <\frac{\delta}{\max(d,\vartheta_{n_0})}$$
for $|\gamma|\leq n_0$ and $r=1,\dots,d$.
Such a triple exists because $\O_A$ is nuclear.
Then by the proof of Prop.\ 3 in \cite{BG}
one can produce an element
$(\phi,\psi,\B)\in\text{CPA}_0(\pi , \Omega(n_0)^{(n)},\delta)$ by
setting 
$\phi:=(\iota\otimes \phi_0)\circ\rho_{n+n_0}$ and
$\psi:=\psi_{n+n_0}\circ(\iota \otimes\psi_0)$. Here
$\psi_m:M_{\vartheta_m}(\O_A)\to\O_A$ takes the matrix $(t_{\alpha
\beta})$
to $\sum s_{\alpha}t_{\alpha\beta}s_\beta^*$. 
We thus compute, by virtue of the previous lemma, 
\begin{align*}
\text{Tr }e^{\phi(f^{(n)})} &=
\sum_{\parbox{3cm}{\centering\scriptsize $\alpha$ a word in 
$\Lambda_A$ of length $n+n_0$}}
\text{Tr }e^{\phi_0(s_\alpha^*f^{(n)}s_\alpha)} \\
&\leq m_0\sum_{|\alpha|=n+n_0} e^{\max\{f^{(n)}(x), x\in[\alpha]\}} \\
&\leq m_0d^{n_0}\sum_{|\alpha|=n}e^{\max\{f^{(n)}(x),x\in[\alpha]\}},
\end{align*}
and therefore by the computation of pressure for (finite type)
subshifts (see, e.g., \cite{DGS})
we obtain
$${P}_\theta(f,\Omega(n_0),\delta)=P^0_\theta(f,\Omega(n_0),\delta)
\leq p_T(f).$$
This inequality implies by the Kolmogorov--Sinai property (Prop.\ 3.4)
that $P_\theta(f)\leq p_T(f)$, 
thus completing the proof.\hfill $\square$
\medskip

\noindent{\it Remark.}
It is not surprising to note that the above theorem produces as a 
special case Boca and Goldstein's result:
$ht(\theta)=h_{\text{top}}(\Lambda_A)=\log r(A)$ \cite{BG}.\medskip

We conclude this subsection with a discussion of the variational
principle in Cuntz--Krieger algebras, comparing
$P_\theta(f)$ with the free energies $h_\sigma(\theta)+\sigma(f)$, where
$h_\sigma(\theta)$ denotes the CNT entropy of $\theta$ computed with
respect to a $\theta$--invariant state $\sigma$.
\medskip

We shall need the following lemma, proven, in a more general form, in 
\cite{PWY}.
\medskip

\noindent{\bf Lemma 5.3.} {\it Any $\theta$--invariant state $\sigma$ on
$\O_A$ containing
$\C(\Lambda_A)$ in its centralizer satisfies
$$h_{\mu}(T)\leq h_\sigma(\theta)$$
where $\mu$ is the $T$--invariant measure on $\Lambda_A$ obtained
restricting $\sigma$. Furthermore any faithful $T$--invariant measure
$\mu$ arises in this way.}\medskip

We have thus obtained the following result.\medskip

\noindent{\bf Theorem 5.4.} {\it Let $f$ be as in Theorem 5.2.
Let $\sigma$ be a $\theta$--invariant state of $\O_A$ centralized
by $\C(\Lambda_A)$, and $\mu$ the shift--invariant measure
on $\Lambda_A$ obtained restricting $\sigma$ to $\C(\Lambda_A)$. Then 
$$h_\mu(T)+\mu(f)\leq h_\sigma(\theta)+\sigma(f)\leq P_\theta(f)=p_T(f).$$
Therefore, by Lemma 5.3, if
($\Lambda_A$, $T$, $f$) admits a faithful equilibrium measure $\mu$,
such a measure extends to
 an equilibrium state $\sigma$ for ($\O_A$, $\theta$, $f$).
}\medskip

\noindent{\it Proof.} Combine the previous Lemma
with Theorem 5.2.\hfill $\square$
\bigskip

{\noindent\it 2.\ Crossed Products}\medskip

{\noindent Now} we turn to crossed products and establish a generalization 
to pressure of a result of Brown \cite{B} which asserts that the 
Voiculescu--Brown entropy of an automorphism of a unital exact 
$C^*$--algebra remains constant under passing to the induced
inner automorphism on the crossed product. Our proof follows
Brown's approach, which in turn is based on a construction
of Sinclair and Smith \cite{SS}.\medskip 

\noindent{\bf Theorem 5.5.} {\it If $\mathcal{A}$ is a unital exact 
$C^*$-algebra, $\theta \in \text{Aut}(\mathcal{A})$, $a$ is a self-adjoint
element in $\mathcal{A}$, and $u$ is the canonical unitary implementing
$\theta$
in $\mathcal{A}\rtimes_\theta\mathbb{Z}$, then
\begin{eqnarray*}
& P_\theta (a) = P_{\text{\normalfont Ad}\, u}(a),
\end{eqnarray*}
where $a$ has been identified on the right with
its image under the natural inclusion $\mathcal{A}\hookrightarrow
\mathcal{A}\rtimes_\theta\mathbb{Z}$.
}\medskip

\noindent{\it Proof.} Without loss of generally we may identify 
$\mathcal{A}$ with its image under a faithful unital representation on a 
Hilbert space $\mathcal{H}$ and, letting $\pi : \mathcal{A}\to\mathcal{B} 
(l^2 (\mathbb{Z},\mathcal{H}))$ be the $*$-monomorphism defined as 
on p.\ 16 of \cite{B}, identify $\mathcal{A}\rtimes_\theta\mathbb{Z}$ with 
the $C^*$-algebra generated by $\pi (\mathcal{A})$ and the image of the 
amplified left regular representation $\lambda$ of $\mathbb{Z}$ in 
$\mathcal{B}(l^2 (\mathbb{Z},\mathcal{H}))$.

The inequality $ P_\theta (a) \leq P_{\text{\normalfont Ad}\, u}
(\pi (a))$ is an immediate consequence of monotonicity (Prop.\ 3.3).

To establish the reverse inequality, we adapt the proof of Brown
\cite{B} for entropy. Let $\Omega\in Pf(\mathcal{A}\rtimes_\theta
\mathbb{Z})$ be of the form $\{\pi (x_1) \lambda_{n_1},\ldots ,
\pi (x_l) \lambda_{n_l} \}$ with $\| x_j \| \leq 1$ for $j=1,\ldots,l$,
 as in the proof of Theorem 3.5 in \cite{B}.
Note that the span of such sets is dense in $\mathcal{A}\rtimes_\theta
\mathbb{Z}$, and so by Prop.\ 3.4 we need only 
show that $P_{\text{Ad}\, u} (id_{\mathcal{A}\rtimes_\theta\mathbb{Z}}, 
a, \Omega ) \leq P_\theta (a)$.
By Lemma 3.4 in \cite{B} there exist a finite set $F\in\mathbb{Z}$
and $\Omega' \in Pf(\mathcal{A})$ such that, if $n\geq 0$, then
\begin{gather*} (\phi ,\psi ,\mathcal{B})\in \text{CPA}(id_{\mathcal{A}} , 
(\Omega')^{(n)}, \delta ) 
\end{gather*}
implies
\begin{gather*}
(\phi' ,\psi' , \mathcal{F} \otimes\mathcal{B} ) \in 
\text{CPA}(id_{\mathcal{A}\rtimes_\theta\mathbb{Z}} ,\Omega^{(n)} ,
2\delta )
\end{gather*}
where $\phi'$ is the u.c.p.\ map $x\mapsto (1\otimes\phi )((p_F
\otimes 1)(x)(p_F\otimes 1))$, $\psi'$ is a u.c.p.\ map, and
$\mathcal{F}$ is the finite-dimensional $C^*$-algebra 
$p_F\mathcal{B} (l^2(\mathbb{Z}))p_F$, with $p_F$ denoting
the projection from $l^2 (\mathbb{Z})$ onto $span\{\xi_t
:t\in F\}$. We can assume that $F$ is of the form $\{-m,-m+1,\ldots
,-1,0,1,\ldots , m-1,m\}$ for some positive integer $m$. By 
Lemma 3.1 in \cite{B}, $\pi (a)= \sum_{t\in\mathbb{Z}} e_{t,t}\otimes
\theta^{-t} (a)$, where the convergence is in the strong operator
topology.

Suppose now that $n\geq 2m+1$ and $(\phi ,\psi ,\mathcal{B})\in 
\text{CPA}(id_{\mathcal{A}} , (\Omega')^{(n)}, \delta )$. Since \\*
$\phi' (\sum_{i=0}^{n-1}
\text{Ad}\, u^i (\pi (a))) = \sum_{i=0}^{n-1} \sum_{t=-m}^{m} p_F 
e_{t,t} p_F \otimes (\phi\circ\theta^{-t-i})(a)$ we have
\begin{align*}
\lefteqn{ \left\| \phi' \left( \sum\nolimits_{i=0}^{n-1}\text{Ad}\, u^i 
(\pi (a))\right) - \sum\nolimits_{k=m}^{n-m-1} \sum\nolimits_{t=-m}^{m} 
p_F e_{t,t} p_F \otimes\left( \phi\circ\theta^{-k}\right) 
(a)\right\| }\hspace*{1cm} \\
&= \left\| \sum\nolimits_{k=-m}^{m-1} \sum\nolimits_{t=-m}^{k} 
p_F e_{t,t} p_F \otimes \left( \phi\circ\theta^{-k}\right) (a) + 
{} \right. \\ 
& \hspace*{1.5cm} \left. \sum\nolimits_{k=n-m}^{n+m-1} 
\sum\nolimits_{t=k-n+1}^{m} 
p_F e_{t,t} p_F \otimes \left( \phi\circ\theta^{-k}\right) (a) 
\right\| \\
&\leq \sum\nolimits_{k=-m}^{m-1} \sum\nolimits_{t=-m}^{k} \left\| p_F 
e_{t,t} p_F \right\| \left\| \left( \phi\circ\theta^{-k}\right) 
(a)\right\| + {} \\
& \hspace*{1.5cm}\sum\nolimits_{k=n-m}^{n+m-1} 
\sum\nolimits_{t=k-n+1}^{m} 
\left\| p_F e_{t,t} p_F \right\| \left\| \left( \phi\circ\theta^{-k}\right) 
(a) \right\| \\
&\leq 4m^2 \| a \| 
\end{align*}
and so
\begin{align*}
\lefteqn{\left| \, \text{log Tr}_{\mathcal{F}\otimes\mathcal{B}} 
e^{\phi' \left( \sum\nolimits_{i=0}^{n-1}\text{Ad}\, u^i (\pi 
(a))\right) } \right. } \hspace*{0.5cm} \\  
& \left. {} - \text{log Tr}_{\mathcal{F}\otimes\mathcal{B}} 
e^{\sum\nolimits_{k=m}^{n-m-1} \sum\nolimits_{t=-m}^{m} p_F 
e_{t,t} p_F \otimes \left( \phi\circ\theta^{-k}\right) (a) } 
\right|  \leq 4m^2 \| a \| 
\end{align*}
by Cor.\ 3.15 in \cite{OP}.
Next observe that if $\beta$ and $\gamma$ are maximal sets of pairwise 
orthogonal spectral
projections for $\sum_{t=-m}^{m} p_F e_{t,t} p_F$ and $\sum_{k=m}^{n-m-1}
(\phi\circ\theta^{-k})(a)$, respectively, then $\beta\otimes\gamma$ is
a maximal set of pairwise orthogonal spectral projections for
$\sum_{k=m}^{n-m-1} \sum_{t=-m}^{m} p_F e_{t,t} p_F \otimes
(\phi\circ\theta^{k})(a)$, and thus
\begin{align*}
\lefteqn{ \text{ Tr}_{\mathcal{F}\otimes\mathcal{B}}\, 
e^{\sum\nolimits_{k=m}^{n-1-m}\sum\nolimits_{t=-m}^{m} p_F e_{t,t} 
p_F \otimes \left( \phi\circ\theta^{-k} \right) (a) }} \hspace*{1cm} \\
&= \sum_{e\in\beta}\sum_{f\in\gamma} e^{\text{Tr}_{\mathcal{F}
\otimes\mathcal{B}} \left[ \left( e\otimes f \right) \left( 
\sum\nolimits_{k=m}^{n-m-1}\sum\nolimits_{t=-m}^{m} p_F e_{t,t} p_F 
\otimes \left( \phi\circ\theta^{-k} \right) (a)\right) \right] } \\
&= \sum_{e\in\beta}\sum_{f\in\gamma}e^{\text{Tr}_{\mathcal{F}
\otimes\mathcal{B}} \left[ \left( e\otimes f \right) \left( \left(
\sum\nolimits_{t=-m}^{m} p_F e_{t,t} p_F \right) \otimes\left( 
\sum\nolimits_{k=m}^{n-m-1}\left( \phi\circ\theta^{-k} \right) (a)\right) 
\right) \right] } \\
&= \sum_{e\in\beta}\sum_{f\in\gamma} e^{\text{Tr}_{\mathcal{F}
\otimes\mathcal{B}} \left[ \left(
e\otimes f \right) \left( 1_\mathcal{F} \otimes 
\sum\nolimits_{k=m}^{n-m-1}\left( \phi\circ\theta^{-k} \right) (a) 
\right) \right] } \\
&= \sum_{e\in\beta}\sum_{f\in\gamma} e^{ \text{Tr}_\mathcal{F} 
\left[ e\cdot 1_\mathcal{F} \right] \,
\text{Tr}_\mathcal{B} \left[ f \left( \sum\nolimits_{k=m}^{n-m-1}\left( 
\phi\circ\theta^{-k} \right) (a)\right) \right] } \\
&= card(\beta ) \sum_{f\in\gamma} e^{\text{Tr}_\mathcal{B}\left[ f 
\left( \sum\nolimits_{k=m}^{n-m-1}\left( \phi\circ\theta^{-k} \right) 
(a) \right) \right] } \\
&= (2m+1) \sum_{f\in\gamma} e^{\text{Tr}_\mathcal{B} \left[ f \left(
\sum\nolimits_{k=m}^{n-m-1}\left( \phi\circ\theta^{-k} \right) (a) \right) 
\right] } \\
&= (2m+1) \text{ Tr}_\mathcal{B}\, e^{\sum\nolimits_{k=m}^{n-m-1}
\left( \phi\circ\theta^{-k} \right) (a) } .\\
\end{align*}
Furthermore, another application of Cor.\ 3.15 in \cite{OP} yields
\begin{align*}
\lefteqn{ \left| \text{log Tr}_\mathcal{B}\, e^{\phi \left( 
\sum\nolimits_{k=0}^{n-1} \theta^{-k}(a)\right)} - \text{log 
Tr}_\mathcal{B} 
\, e^{\phi \left( \sum\nolimits_{k=m}^{n-m-1} 
\theta^{-k}(a)\right) } \right| } \hspace*{1cm} \\
&\leq \left\| \phi \left( \sum\nolimits_{k=0}^{n-1} \theta^{-k}(a) - 
\sum\nolimits_{k=0}^{m-1} \theta^{-k}(a)\right) \right\| \\
&\leq \sum\nolimits_{k=0}^{m-1}  \left\| 
\left( \phi\circ\theta^{-k}\right) (a) \right\| + \sum\nolimits_{k=n-m}^{n-1}  
\left\| \left( \phi\circ\theta^{-k}\right) (a) 
\right\| \\
&\leq 2m. 
\end{align*}
Combining these estimates we obtain
\begin{align*}
\lefteqn{ \left| \,\text{log Tr}_{\mathcal{F}\otimes\mathcal{B}} 
\, e^{\phi' \left( \sum\nolimits_{i=0}^{n-1}\text{Ad}\, u^i 
(\pi (a))\right) } - \text{log Tr}_\mathcal{B}\, e^{\phi 
\left( \sum\nolimits_{k=0}^{n-1} \theta^{-k}(a)\right) }
\right| } \hspace*{1.5cm}\\
&\leq \left| \,\text{log Tr}_{\mathcal{F}\otimes\mathcal{B}} 
\, e^{\phi' \left( \sum\nolimits_{i=0}^{n-1}\text{Ad}\, u^i 
(\pi (a))\right) } \right. \\
& \hspace*{1.7cm}{} - \left. \log \text{ Tr}_{\mathcal{F}\otimes
\mathcal{B}}\, e^{\sum\nolimits_{k=m}^{n-m-1} \sum\nolimits_{t=-m}^{m} 
p_F e_{t,t} p_F \otimes \left( \phi\circ\theta^{-k}\right) (a) } 
\right| + {}\\
& \hspace*{1cm}\left| \,\log \text{ Tr}_{\mathcal{F}\otimes
\mathcal{B}}\, e^{\sum\nolimits_{k=m}^{n-m-1} 
\sum\nolimits_{t=-m}^{m} p_F e_{t,t} p_F \otimes \left( 
\phi\circ\theta^{-k}\right) (a) } \right. \\
& \hspace*{1.7cm}{} - \left. \text{log Tr}_{\mathcal{B}} 
\, e^{\sum\nolimits_{k=m}^{n-m-1} \left( \phi\circ\theta^{-k} \right) 
(a) } \right| +{} \\
& \hspace*{1cm} \left| \text{log Tr}_{\mathcal{B}} \, 
e^{\sum\nolimits_{k=m}^{n-m-1} \left( \phi\circ\theta^{-k} \right)
(a) } - \text{log Tr}_\mathcal{B}\, e^{\phi \left( 
\sum\nolimits_{k=0}^{n-1} \theta^{-k}(a)\right) } \right| \\
&\leq 2m(2m+1) \| a \| +\log (2m+1). 
\end{align*}
Since the above holds for any $(\phi , \psi , B) \in
\text{CPA}(id_{\mathcal{A}} 
,\mathcal{A} , (\Omega')^{(n)}, \delta )$, it follows that
\begin{align*}
\lefteqn{ \log Z_{\text{Ad}\, u, n}(id_{\mathcal{A}\rtimes_\theta\mathbb{Z}}, 
\pi (a), \Omega, 2\delta ) } \hspace*{1cm} \\
&\leq \log Z_{\theta, n}(id_A , a, \Omega' , \delta ) + 2m(2m+1) 
\| a \| +\log (2m+1). 
\end{align*}
Thus, letting $n$ vary while $m$ remains fixed, we infer that
\begin{gather*}
P_{\text{Ad}\, u} (id_{\mathcal{A}\rtimes_\theta\mathbb{Z}}, 
\pi (a), \Omega, 2\delta ) \leq P_\theta (id_\mathcal{A}, a, \Omega' ,
\delta ) 
\end{gather*}
We conclude that $P_{\text{Ad}\, u}(id_{\mathcal{A}\rtimes_\theta
\mathbb{Z})}, \pi (a), \Omega  ) \leq P_\theta (id_\mathcal{A}, a)$, 
completing the proof of the theorem.\hfill $\square$

\section{Variational Principle for bimodule algebras}

\def\U{{\cal U}}
\def\A{{\cal A}}
\def\L{{\cal L}}
\def\K{{\cal K}}
\def\O{{\cal O}}
\def\C{{\cal C}}
\def\B{{\cal B}}
\def\F{{\cal F}}
\def\H{{\cal H}}

This section is divided into three subsections. In the first subsection
we obtain a variational principle for a class
of $C^*$--algebras generated by a Hilbert bimodule \cite{P} which
generalizes the results of subsection 5.1, in the second subsection we
discuss
equilibrium states for the same class, and in the last subsection we
discuss an
application to Matsumoto $C^*$--algebras associated to a subshift.

Let $\A$ be a unital exact $C^*$--algebra faithfully represented on a
Hilbert space, and $X$ a Hilbert
$\A$--bimodule, i.e. 
 $X$ is a  right Hilbert $\A$--module endowed with
a faithful left action of $\A$ given by a 
 unital $^*$--monomorphism
$\A\to\L_\A(X)$
 into
the algebra  of right $\A$--linear endomorphismsms of 
$\A$. We will always assume that $X$ is finitely generated as a right module.
Following \cite{P}, we construct the universal $C^*$--algebra algebra $\O_X$
generated by $X$ and a unital copy of $\A$ satisfying 
\begin{gather*}
{x}^*ax'=<x, ax'>,\quad x,x'\in X, a\in\A, \\
\sum_i x_ix_i^*=I, \quad\{x_i\}_i\text{ a basis of } X,
\end{gather*}
where a basis of $X$ is a finite subset such that
$$x=\sum_i x_i<x_i, x>,\quad x\in X.$$
The Banach subspaces 
$\L_\A(X^{\otimes r}, X^{\otimes s})$ 
are isometrically embedded in $\O_X$ in a manner respecting the
inclusions $\L_\A(X^{\otimes r}, X^{\otimes s})\hookrightarrow 
\L_\A(X^{\otimes r+1}, X^{\otimes s+1})$ under which  
$T\in\L_\A(X^{\otimes r}, X^{\otimes s})$ is identified with
$T\otimes 1_X\in\L_\A(X^{\otimes r+1}, X^{\otimes s+1})$.
The algebra $\O_X$ carries an automorphic action 
$\gamma:{\Bbb T}\to\O_X$ of the circle given by 
\begin{gather*}
\gamma_z(x)=zx,\quad z\in{\Bbb T}, \\
\gamma_z(a)=a,\quad a\in\A
\end{gather*}
and referred
to as the gauge action.
The corresponding fixed-point algebra will be denoted
by $\O_X^{(0)}$.
It is a fact that $\O_X$ is an exact $C^*$--algebra
if $\A$ is exact. This property
 has been proven in \cite{DS} for general Hilbert bimodules.
 However, for the class of special modules we will be considering,
exactness
 of $\O_X$, under the corresponding assumption 
for the coefficient algebra, will result from the proof of Prop.
6.9.
As noticed in \cite{P}, 
if $\alpha$ is an automorphism of $\A$ 
and $X=\A$ with right Hilbert $\A$--module structure
\begin{gather*}
xa=x\alpha(a), \\
<x,y>=\alpha(x^*y),\quad x,y\in X,\, a\in\A
\end{gather*}
and left action given by left multiplication, then
$\O_X=\A\rtimes_\alpha{\Bbb Z}$.
If $\A$ is commutative and finite-dimensional, $\O_X$
is a Cuntz--Krieger algebra. 
Matsumoto algebras associated to 
subshifts \cite{M}  also arise as Pimsner algebras \cite{PWY}.
We further assume that $X$ admits a basis $\{x_i\}_i$
such that
$$<x_i,ax_j>=0,\quad i\neq j,\, a\in\A.$$
In other words,  as a right $\A$--module,
$X=q_1\A\oplus\dots\oplus q_d\A$,
where each $q_i$ is a projection of $\A$,
while the
left $\A$--action is defined by the diagonal action of
$^*$--monomorphisms 
$\rho_i:\A\to \A $ with
$$\rho_i(I)=q_i.\eqno(6.1)$$
Thus $\O_X$ is the universal $C^*$--algebra generated by a unital  copy of
$\A$ and
partial isometries $x_1,\dots, x_d$ satisfying 
$$ax_i=x_i\rho_i(a),\quad i=1,\dots,d,\, a\in\A,\eqno(6.2)$$
$$\sum{x_i}x_i^*=I,\eqno(6.3)$$
$$x_i^*x_i=q_i.\eqno(6.4)$$
All of the bimodules in the above examples admit an
orthogonal basis. 
We recall for convenience Pimsner's construction
of the bimodule generating $\O_A$ \cite{P}. 
Take
 $\A=\oplus_1^d{\Bbb C}$ as the coefficient algebra, and let
 $A=(A_{rs})\in M_d(\{0,1\})$. Then if $\{p_i\}$ stands for the set of
minimal projections of $\A$,  
$X$ is the bimodule $q_1\A\oplus\dots\oplus q_d\A$,
where
 $q_i=\sum_{j}A_{ij} p_j$, and $\rho_i(p_j)=\delta_{i,j} q_i$.
The corresponding basis elements $\{x_i\}$, usually denoted by $\{s_i\}$,
generate $\O_X$ and satisfy 
\begin{gather*}
\sum_i s_is_i^*=I, \\
s_i^*s_i=q_i, \\
s_is_i^*=p_i.
\end{gather*}

\noindent{\it 1.\  The variational principle}\medskip

Our first aim is to give an upper and lower bound for the pressure
of a self-adjoint element $a$ in a suitable $C^*$--subalgebra
of $\O_X$  with respect to  the u.c.p. map 
$$\theta: b\to\sum x_ibx_i^*$$
on $\O_X$.
These bounds will lead, under certain circumstances, to a computation of 
$P_\theta(a)$ and to the variational principle. 
\medskip

\noindent{\it Remark.} In contrast with the class of $C^*$--algebras
considered by Neshveyev and St\o rmer \cite{NS}, $(\O_X,\theta)$ is
usually not asymptotically Abelian. Indeed, if $q_1+\dots+q_d$ is
invertible,
$\theta$ restricts to a unique monomorphism $\sigma$
of $\A'\cap\O_X$
such that $\sigma(t)x=xt$ for $t\in\A'\cap\O_X$ and $x\in X$ 
\cite{PWY}.
If $\O_X$ is simple and there is a nonscalar
element $t\in\A'\cap\O_X$,
one has $x_i\sigma^r(t)=\sigma^{r+1}(t)x_i$, and therefore
 $[x_i,\theta^r(t)]=(\sigma^{r+1}(t)-\sigma^r(t))x_i$. If this 
tended to $0$
for all $i$, then one would have $\sigma(t)=t$, so
that $t$ would be an
element in the
centre
of $\O_X$.
\medskip

If $\alpha=(i_1,\dots,i_r)$ we write $x_\alpha$ for ${x_{i_1}\dots x_{i_r}}$
and denote by $|\alpha|$ the length $r$ of $\alpha$.
We will restrict $a$ to be an element of the following
amplification of the coefficient algebra: 
$$\D:=\left\{b\in\O_X^{(0)}: x_\alpha^*bx_\beta=0, |\alpha|=|\beta|,
\alpha\neq \beta\right\}.$$
Note that $\D$ is a unital  $C^*$--subalgebra
containing $\A$ and elements of the form
$x_\alpha a x_\alpha^*$, $a\in\A$,
and is invariant under $\theta$ and the maps $\text{Ad}\,x_i^*$,
$i=1,\dots, d$.
Furthermore
$\D$  is an exact $C^*$--algebra since it is a subalgebra of
$\O_X^{(0)}$, which is exact.
Notice that the  $\text{Ad}\,x_i^*$'s restrict to endomorphisms of
$\D$.
The closed subspace 
$X_{\D}=X\D$ of $\O_X$ is a Hilbert bimodule
over $\D$ isomorphic to
$q_1\D\oplus\dots \oplus q_d\D$ as a right Hilbert module
with diagonal left action induced by $\text{Ad}\,x_i^*$,
$i=1,\dots,d$. One has
$\O_{X_\D}=\O_X$.
 This construction is familiar in the case of  
Cuntz--Krieger algebras.\medskip

\noindent{\bf Proposition 6.1} {\it If $X$ is the  Hilbert
bimodule defining the Cuntz--Krieger algebra $\O_A$,
$\D$ is the Abelian $C^*$--subalgebra $\C(\Lambda_A)$.}\medskip

\noindent{\it Proof.} The inclusion 
$\C(\Lambda_A)\subset\D$ follows from the fact
that $\D$ is $\theta$--invariant and contains the range projections 
$p_i$, $i=1,\dots,d$. To show the opposite inclusion we consider a
sequence
of conditional expectations $(E_r)_r$ onto the finite-dimensional
$C^*$--subalgebras $\F_r$ generated by $\{{s_\alpha}p_is_\beta^*, 
|\alpha|=|\beta|=r,i=1,\dots, d\}.$ We choose  each $E_r$
to be invariant under a faithful trace of $\O_A^{(0)}$ obtained
restricting
a $\beta$--KMS state $\omega$ of $\O_A$ for the one-parameter group
$t\to\gamma_{e^{2\pi it}}.$
Then the KMS condition 
$$\omega(s_i^*t)=e^\beta\omega(ts_i^*),\quad t\in\O_A,$$
yields
$$s_i^*E_{r+1}(t)s_j=E_{r}(s_i^*ts_j),\quad t\in\O_A^{(0)}.$$
Thus if $t\in\D$ then 
$E_r(t)=\sum_{|\alpha|=r}s_\alpha E_0(s_\alpha^*ts_\alpha)s_\alpha^*$, 
which is contained in $\C(\Lambda_A)$ since the range of $E_0$
is the linear span of the range projections $p_i=s_is_i^*$, $i=1,\dots,
d$.
Since $(E_r)_r$ converges to the identity, we have
$t\in\C(\Lambda_A)$.
$\hfill\square$\medskip

Note that, on the other hand, if $X$ is a Hilbert bimodule defined
by an automorphism $\alpha$ of $\A$, then $\D=\A$.

Recall that  a one-sided subshift 
is a closed subset 
of the compact space $\{1,\dots,d\}^{\Bbb N}$
such that $T(\Lambda)=\Lambda$, where
$T((a_k)_k)=(a_{k+1})_k$
is the left shift epimorphism of the full shift
space.
Let  $\Lambda^{(n)}$ stand for the set
of $n$-tuples $\alpha=(i_1,\dots,i_n)$ for which there is
$(a_k)_k\in\Lambda$
such that $a_1=i_1,\dots, a_n=i_n$, and set $\vartheta_n=\text{Card
}(\Lambda^{(n)})$.

We associate to an orthogonal basis $\{x_i\}_{i=1}^d$ 
of a Hilbert bimodule $X$ the set 
$$\Lambda_{\{x_i\}}=\{(a_k)_k\in\{1,\dots,d\}^{\Bbb N}: x_{a_1}\cdots
x_{a_n}\neq 0 \text{ for all }n\in{\Bbb N}\}.$$  
We will write $\Lambda=\Lambda_{\{x_i\}}$.
It easily checked that $\Lambda$ is a one--sided subshift.
The relation $\sum_i x_ix_i^*=I$
shows that $\Lambda\neq\emptyset$.
For $\alpha=(i_1,\dots,i_n)\in\Lambda^{(n)}$ we set
$x_\alpha=x_{i_1}\dots x_{i_n}$, $q_\alpha=x_\alpha^*x_\alpha$,
$p_\alpha=x_\alpha x_\alpha*$
and $\rho_\alpha=\rho_{i_1}\dots\rho_{i_n}$.
Note that $\C(\Lambda)$ embeds naturally
in $\O_X$ as the $C*$--subalgebra generated by the projections $p_\alpha$,
$\alpha\in\cup\Lambda^{(n)}$.

We define the topological 
entropy of
 the
 the action of 
$X$
on $\A$
 by
$$ht(\A,
X):=\sup_{\Omega\in Pf(\A)}\,\,\sup_{\delta>0}\,\,\limsup_n\frac{1}{n}\log
rcp(\A,\Omega^{(n,X)},\delta),$$
where $$\Omega^{(n,X)}:=\{\rho_\mu(t) : t\in\Omega, |\mu|\leq n-1\}.$$
Our aim is to prove the following.
\medskip

\noindent{\bf Theorem 6.2.} {\it Let $\A$ be a unital exact $C^*$--algebra
and $X$ a Hilbert $\A$--bimodule defined as above by $^*$--monomorphisms 
$\rho_i:\A\to\A$, $i=1,\dots,d$, with the property that
$\sum_{i=1}^d\rho_i(I)$ is invertible. 
Consider
the $C^*$--dynamical system $(\O_X, \theta)$. Suppose that
$ht(\A, X)=0$.
If $a\in\D$ is  positive
and satisfies
\begin{gather*}
[a, x_\alpha^*ax_\alpha]=0 \\
[a, q_\alpha]=0
\end{gather*}
for $|\alpha|$ sufficiently large, then
$$P_\theta(a)=\lim_n\frac{1}{n}\log\sum_{\alpha\in\Lambda^{(n-1)}}
e^{\|x_\alpha^*a^{(n)}x_\alpha\|}.$$
(In particular, if $a\in\C(\Lambda)$ then $P_\theta(a)$ coincides with the
classical pressure of $a$ w.r.t. the shift on $\Lambda$.)
Furthermore one has
$$\sup_{\sigma}h_\sigma(\theta)+\sigma(a)=P_\theta(a),$$
where $h_\sigma(\theta)$ denotes the Sauvageot--Thouvenot entropy
of $\theta$ and the
supremum is taken over all $\theta$--invariant states 
of $\O_X$.}\medskip

Before proving the theorem, we discuss an example where
the condition $ht(\A, X)=0$ is easily checked.\medskip 

\noindent{\it Example.} Consider the case in which
$\A$ is the inductive limit of fi\-ni\-te-di\-mensional
$C^*$--algebras $(\F_r)_{r\in{\Bbb N}}$, each one invariant under
$\rho_1,\dots,\rho_d$, and suppose that
$\A$ admits  a faithful trace
$\tau$.
Let $\Omega$ be a finite set contained in some $\F_{r_0}$.
Then $\Omega^{(n, X)}\subset\F_{r_0}$ for all $n\in{\Bbb N}$.
Consider the $\tau$--preserving conditional
expectation 
$E:\A\to\F_{r_0}$ and let $\iota:\F_{r_0}\to\A$ be 
the inclusion. Clearly $(E, \iota,  \F_{r_0})\in
\text{CPA}(\Omega^{(n,X)},\delta)$
for all $\delta>0$ and $n\in{\Bbb N}$, and so $ht(\A, X)=0$.\medskip

We shall divide the proof of Theorem 6.2 into three parts.
In the first and second part we will give upper and lower
bounds for the pressure of
$a$, and in the 
third part we will prove the variational principle.
We start by showing why we require $\sum_i q_i$ to be 
invertible.\medskip

\noindent{\bf Lemma 6.3.} {\it If $\sum_{i=1}^d q_i$ is invertible,
then $\theta$
is faithful on $\O_X$.}\medskip

\noindent{\it Proof.} The equality $\theta(t^*t)=0$ implies
$q_it^*tq_i=x_i^*\theta(t^*t)x_i=0$, and thus
$tq_i=0$ for all $i$, so that $t=0$.
$\hfill\square$\medskip

Let $a\in\D$ be a self-adjoint element.
We introduce the following notion of pressure for $a$ with respect
to the bimodule $X$. Let $\Omega\in Pf(\A)$, $\delta>0$, and
$n\in{\Bbb N}$.
Setting $a^{(n)}:=\sum_0^{n-1}\theta^j(a)$ as usual,
we define the partition function  
\begin{align*}
\lefteqn{Z_{X,n}(\D, a,
\Omega,\delta)=}\hspace*{1.2cm} \\
&\ \inf\bigg\{\sum_{\alpha\in\Lambda^{(n-1)}}\text{Tr }
e^{\phi\left(x_\alpha^*a^{(n)}x_\alpha\right)} :
(\phi,\psi,\B)\in\text{CPA}(\D,\Omega^{(n, X)},\delta)\bigg\},
\end{align*}
and the corresponding  
$P_{X}(\D, a,\Omega,\delta)$ and 
$P_{X}(\D, a,\Omega)$ are obtained in the usual manner.
We set
$$P_{X}(\D,a)=\sup_{\Omega\in
Pf(\A)}P_X(\D,a,\Omega).$$
We emphasize that  this pressure
is computed
by means of approximations of a faithful representation of
$\D$ via factorizations through finite-dimensional 
$C^*$--algebras. However,
we only let $\Omega$ range over
finite subsets 
of $\A$.
The definition of $Z_{X,n}(\D, a, \Omega, \delta)$ suggests considering
for an element $a\in\O_X$ the sequence
$$\sum_{\alpha\in\Lambda^{(n-1)}}
e^{\max\text{spec }x_\alpha^*a^{(n)}x_\alpha},$$
which resembles the classical partition function defining
the pressure of a subshift (see \cite{DGS}, e.g.).\medskip

\noindent{\bf Lemma 6.4.} {\it If $a$ is a positive element of $\O_X$ then
$$\lim_n\frac{1}{n}\log\sum_{\alpha\in\Lambda^{(n-1)}}
e^{\|x_\alpha^*a^{(n)}x_\alpha\|}=:\ell$$ exists
and 
$$ h_{\text{\normalfont top}}(\Lambda)\leq\ell\leq
\|a\|+h_{\text{\normalfont top}}(\Lambda).$$}\medskip

\noindent{\it Proof.} If $|\alpha|=n$ and $|\beta|=m$ then
\begin{align*}
\lefteqn{x_{\alpha\beta}^*a^{(n+m+1)}x_{\alpha\beta}}\hspace*{0.5cm} \\
&=x_\beta^*\left(x_\alpha^*a^{(n+1)}x_\alpha\right)x_\beta+
x_\beta^*\left(x_\alpha^*\theta^{n}(\theta(a)+
+\dots+\theta^{m}(a))x_\alpha\right)x_\beta\\
&=x_\beta^*\left(x_\alpha^*a^{(n+1)}x_\alpha\right)x_\beta+
x_\beta^*q_\alpha(\theta(a)+\dots+\theta^m(a))q_\alpha x_\beta\\
&=x_\beta^*\left(x_\alpha^*a^{(n+1)}x_\alpha\right)x_\beta+
q_{\alpha\beta}
\left(x_\beta^*(\theta(a)+\dots+\theta^m(a))x_\beta\right)q_{\alpha\beta}
\end{align*}
since 
$x_\beta^*q_\alpha=\sum_{|\gamma|=|\beta|}x_\beta^*q_\alpha x_\gamma
x_\gamma^*=q_{\alpha\beta}x_\beta^*.$
Now the previous term is bounded above by
$$
\|x_\alpha^*a^{(n+1)}x_\alpha\|+\|x_\beta^*a^{(m+1)}x_\beta\|
$$
and so
$s_{n}:= \sum_{\alpha\in\Lambda^{(n)}}
e^{\|x_\alpha^*a^{(n+1)}x_\alpha\|}$
satisfies $s_{n+m}\leq s_ns_m$. It follows that
$\lim\frac{1}{n}\log(s_n)$ exists and equals $\inf\frac{1}{n}\log(s_n)$.
The  upper and lower bounds for $\ell$ follow
from the inequalities $0\leq\|x_\alpha^*a^{(n)}x_\alpha\|\leq n\|a\|$ for
$\alpha\in\Lambda^{(n-1)}$.$\hfill\square$\medskip

\noindent{\bf Proposition 6.5.} {\it If $a\in\D$ is a self-adjoint
element then
$$P_X(\D,a)\leq
ht(\A,X)+\lim_n\frac{1}{n}\log\sum_{\alpha\in\Lambda^{(n-1)}}
e^{\max\text{\normalfont 
spec
}x_\alpha^*a^{(n)}x_\alpha}.$$}\medskip

\noindent{\it Proof.} The proof is 
straightforward once we note that, by Arveson's extension
theorem \cite{A}, every unital c.p.\ map
$\phi:\A\to M_N({\Bbb C})$ extends to a unital
c.p. map $\tilde{\phi}:\D\to M_N({\Bbb C})$.
$\hfill\square$\medskip

Before establishing an upper bound for $P_\theta(a)$ for
certain $a\in\D$,
we shall need a few preliminary results. The first two lemmas are
immediate, and so we omit the proofs. 
\medskip

\noindent{\bf Lemma 6.6.} {\it If $a\in\D$ is a self-adjoint
element, and if
$\lambda\in{\Bbb R}^+$,
\begin{description}
\item{\rm a)} 
$P_X(\D,a+\lambda)\leq P_X(\D, a)+\lambda,$ 
\item{\rm b)} 
$P_X(\D,a-\lambda)\geq P_X(\D, a)-\lambda.$
\end{description}}\medskip

\noindent{\bf Lemma 6.7.} {\it  Set
 $$\phi_m:b\in\O_X\to
(x_\mu^*bx_\nu)_{\mu,\nu\in\Lambda^{(m)}}\in M_{\vartheta_m}(\O_X).$$
Then for $j=0,\dots,n-1$, $|\beta|\leq|\alpha|\leq n_0$, and $t\in\A$ the
$(\mu , \nu)$ entry of
$$\phi_{n+n_0-1}\theta^j({x_\alpha} tx_\beta^*)$$ is 
 nonzero only 
if $\mu$ and $\nu$ are of the form
$\mu=\delta\alpha\mu'$, $\nu=\delta\beta\mu'\nu'$ with
$|\delta|=j$ and $|\nu'|=|\alpha|-|\beta|$. The corresponding entry
is $$\rho_{\mu'}(q_{\delta\alpha} t q_{\delta\beta})x_{\nu'}.$$}\medskip

\noindent{\bf Lemma 6.8.} {\it Let $\A\subset\B(\H)$ be a  unital
$C^*$--algebra and
let 
$\phi:\A\to\B(\H)$ be a unital c.p. map. Let $x$ be an element of the unit
ball of $\A$ and $p$, $q$ projections of $\A$ such that 
$\|\phi(y)-y\|<\delta$ for each $y\in\{x,p,q\}$ and some $\delta<1$. 
Then $\|\phi(pxq)-pxq\|<11\delta^{\frac12}$.}\medskip

\noindent{\it Proof.} By Stinespring's theorem \cite{St}
there is a Hilbert space $\K$, an isometry
$V:\H\to\K$, and a unital $^*$--representation $\pi:\A\to\B(\K)$
such that $\phi(t)=V^*\pi(t)V$, $t\in\A$. In Stinespring's construction 
$\K$ is the tensor product Hilbert bimodule
$\A\otimes_{\Bbb C}\H$, where $\A$ is regarded as a $\A$--${\Bbb C}$
Hilbert bimodule with $\A$--valued
inner product defined by $<a,b>=\phi(a^*b)$. One has $\pi(a)=a\otimes I$
for $a\in\A$
and $V\xi=I\otimes\xi$ for $\xi\in\H$.
One checks that 
\begin{align*}
\|\pi(p)V-Vp\|^2 &\leq \|\phi(p)-\phi(p)p-p\phi(p)+p\| \\
&\leq \|\phi(p)-p\|+\|(p-\phi(p))p\|+\|p(\phi(p)-p)\| \\
&< 3\delta,
\end{align*}
and similarly for $q$.
This implies
\begin{align*}
\|\phi(p)\phi(x)\phi(q)-\phi(pxq)\| &= \|V^*\pi(p)VV^*\pi(x)V
V^*\pi(q)V-V^*\pi(pxq)V\| \\
&\leq 4(3\delta)^{\frac12} \\
&< 8\delta^{\frac12}.
\end{align*}
On the other hand, by our assumption we have 
$$\|\phi(p)\phi(x)\phi(q)-pxq\|<3\delta,$$
which, when combined with the previous estimate, yields the result.
$\hfill\square$\medskip

We are now ready to give an upper bound for $P_\theta(a)$.\medskip

\noindent{\bf Proposition 6.9.} {\it Let $X$ be a Hilbert bimodule
over a unital exact $C^*$--algebra $\A$.  Suppose that 
$X=q_1\A\oplus\dots\oplus q_d\A$ as a right Hilbert module
with left action  defined by unital
$^*$--monomorphisms $\rho_i:\A\to q_i\A q_i$, $i=1,\dots,d$.
Then for any 
self-adjoint element $a\in\D$
commuting asymptotically with the domain projections $\{q_\mu,
\mu\in\cup_n\Lambda^{(n)}\}$, i.e.,
$$\lim_{|\mu|\to\infty}\|[a,\ q_\mu]\|=0, $$
we have
$$P_\theta(a)\leq P_{X}(\D, a)\leq ht(\A,
X)+\limsup_n\frac{1}{n}\log\sum_{\alpha\in\Lambda^{(n-1)}}
e^{\max\text{\normalfont 
spec } x_\alpha^*a^{(n)}x_\alpha}.$$}\medskip

\noindent{\it Proof.} We need only show the first inequality. 
Let $\Omega\subset\A$ be a finite subset of the unit ball containing
$I$. For
$n_0\in{\Bbb N}$ we set $$\Omega(n_0)=\{x_\alpha tx_\beta^*,\quad
|\beta|\leq
|\alpha|\leq n_0,
t\in\Omega\}.$$ Since $\cup_{n_0,\Omega}\Omega(n_0)\cup\Omega(n_0)^*$
is total in $\O_X$, it suffices by monotonicity and the
Kolmogorov--Sinai
property to show that $P_\theta(a,\Omega(n_0))\leq P_{{X}}(\D, a)$ 
for all $\Omega\in Pf(\A)$ and $n_0\in{\Bbb N}$.
Following the proof of Lemma 7.5 in \cite{PWY}, which in turn goes back to
\cite{B}, 
given any subset
 $\Delta\in Pf(\O_X^{(0)})$, $\delta>0$, and $n_0\in{\Bbb N}$,  
we can find a finite subset $F\subset{\Bbb N}$, which depends only
on $\delta$ and $n_0$  and not on
$\Delta$, such that if
$$(\phi,\psi,\B)\in 
\text{CPA}\Big(\O_X^{(0)},\Delta^{(\max F, X)},
\frac{\delta}{2\max_{p\in F}
\vartheta_p}\Big)$$
then there is a triple $(\phi', \psi', \B')\in \text{CPA}_0(\O_X,
\cup_{|\gamma|\leq n_0}\Delta x_\gamma,\delta)$
with $\B'=M_{\vartheta_F}\otimes\B$
and $\phi'=\iota_{M_{\vartheta_F}}\otimes\phi\circ S_F$, where
$\vartheta_F=\sum_{p\in F}\vartheta_p$ and
$$S_F:b\in\O_X\to
(x_\alpha^*m_{|\alpha|-|\beta|}(b)x_{\beta})_{|\alpha|,|\beta|\in F}\in
M_{\vartheta_F}(\O_X^{(0)}).$$
Here $m_{k}$ denotes the natural projection onto the
subspace $\O_X^{(k)}$ of elements which transform like 
$z^k$ under the gauge action.
Let us apply the above construction 
to the parameters $11\delta^{\frac12}$, $n_0$, and any $\Delta\in
Pf(\O_X^{(0)})$, and 
find the corresponding $F$.
Pick $( \tilde{\phi}, \psi, M_N)\in
\text{CPA}(\D,\Omega^{(n+n_0+\max F,
X)},\frac{\delta}{4{\max_{p\in
F}\vartheta_p}^2})$
and extend $\tilde{\phi}$ by Arveson's theorem \cite{A} to a u.c.p. map
$\phi$ on
$\O_X^{(0)}$, so that
$$(\phi,\psi, M_N)\in \text{CPA}\Big(\O_X^{(0)},
\Omega^{(n+n_0+\max F, X)},\frac{\delta}{4{\max_{p\in
F}\vartheta_p}^2}\Big).$$ 
Since $I\in\Omega$, we have
$$\|(\iota-\psi\phi)(\rho_\mu(q_\beta))\|<
\frac{\delta}{4{\max_{p\in F}\vartheta_p}^2},\quad |\mu|+|\beta|\leq
n+n_0+\max
F-1.$$
Also,
$$\|(\iota-\psi\phi)(\rho_\mu(t))\|<
\frac{\delta}{4{\max_{p\in F}
\vartheta_p}^2},\quad |\mu|\leq n+n_0+\max
F-1.$$
By Lemma 6.8,
$$\|(\iota-\psi\phi)(\rho_\mu(q_\beta t q_\alpha))\|<
\frac{11}{2\max_{p\in F}\vartheta_p}\delta^{\frac12} $$
for $|\mu|+|\beta|, |\mu|+|\alpha|\leq
n+n_0+\max F-1$ and
$t\in\Omega$.
We set
$${\Omega'}_n=\{\rho_\mu(q_\alpha t q_\beta), |\mu|+|\alpha|,
|\mu|+|\beta|\leq n-1, t\in\Omega\},$$
and so $(\phi,\psi, M_N)\in \text{CPA}({\Omega'}_{n+n_0}^{(\max F,X)},
\frac{11}{2\max_{p\in F}\vartheta_p}\delta^{\frac12}).$
By the construction at the beginning we can find
  $$(\phi',\psi',\B')\in \text{CPA}_0(\cup_{|\gamma|\leq
n_0}{\Omega'}_{n+n_0}x_\gamma, 11\delta^{\frac12}).$$
Consider the u.c.p. map $\psi_m: (t_{\mu \nu})\to
M_{\vartheta_m}(\B(\H))\to\sum_{|\mu|, |\nu|=m}{x_\mu}t_{\mu
\nu}x_\nu^*\in\B(\H)$ and the contractive c.p. map $\phi_m$ defined
in Lemma 6.7.
We claim that $$(\phi'',\psi'',M_{\vartheta_{n+n_0-1}}\otimes\B')\in
\text{CPA}_0(\O_X, \Omega(n_0)^{(n)},11\vartheta_{n_0}\delta^{\frac12}),$$
where
$\psi'':=\psi_{n+n_0-1}\circ(\iota_{M_{\vartheta_{n+n_0-1}}}\otimes\psi')$
and 
$\phi'':=(\iota_{M_{\vartheta_{n+n_0-1}}}\otimes\phi')\circ\phi_{n+n_0-1}$.
To establish the claim, note first that, for $j=0,\dots,n-1$, $t\in\Omega$,
and $|\beta|\leq|\alpha|\leq n_0$,
$$\|(\psi''\phi''-\iota)(\theta^j(x_\alpha t x_\beta^*))\|\leq 
\|\iota_{M_{\vartheta_{n+n_0-1}}}\otimes(\psi'\phi'-\iota)\circ
\phi_{n+n_0-1}\circ\theta^j(x_\alpha t x_\beta^*)\|$$
Using notation from Lemma 6.7, the term on the right is bounded by
$$\max_{|\delta|=j}\sum_{|\nu'|=|\alpha|-|\beta|}\max_{|\mu'|}
\|(\psi'\phi'-\iota)(\rho_{\mu'}(q_{\delta\alpha} t
q_{\delta\beta})x_\nu')\|\leq11\vartheta_{n_0}\delta^{\frac12}.$$
Notice that the range of $\phi''$ is a matrix algebra
of rank $\vartheta_{n+n_0-1}\vartheta_FN$. However, we can reduce
this rank by taking into account degeneracies, and so
we will consider $\phi''$ as a
maps with range 
in matrices of rank $\sum_{p\in F}\vartheta_{n+n_0-1+p}N$.
Assume for the moment that $a\geq0$. Given  
 $\epsilon>0$ let $q\in{\Bbb N}$ be such that
$\|q_\gamma a q_\gamma-a^{\frac12} q_\gamma a^{\frac12}\|<\epsilon$, $|\gamma|\geq
q$.
Then for all $j$,
\begin{align*}
q_\gamma\theta^{j}(a)q_\gamma &= \sum_{|\beta|=j}
x_{\beta}
\rho_\beta(q_\gamma)a\rho_\beta(q_\gamma)
x_\beta^* \\
&\leq \sum_{|\beta|=j}x_{\beta}
a^{\frac12}\rho_\beta(q_\gamma)a^{\frac12}
x_\beta^* \\
&\hspace*{1.5cm} \ +\Big\|\sum_{|\beta|=j}x_\beta
(\rho_\beta(q_\gamma)a\rho_\beta(q_\gamma)-
a^{\frac12}\rho_\beta(q_\gamma)a^{\frac12})x_\beta^*\Big\| \\
&\leq \theta^j(a)+\epsilon.
\end{align*}
We compute, for $n\geq \max F+n_0+q$, 
\begin{align*}
\lefteqn{\text{Tr }\exp({\phi''(a^{(n)})})}\hspace*{1cm} \\
&= \sum_{p\in
F}\,\sum_{\alpha\in
\Lambda^{(n+n_0-1+p)}}\text{Tr }\exp({\tilde{\phi}
(x_\alpha^*a^{(n)}x_{\alpha})}) \\
&\leq
\sum_{p\in
F}\,\sum_{\gamma\in\Lambda^{(p+n_0+q)}}\,\sum_{\alpha\in\Lambda^{(n-q-1)}}
\text{Tr }\exp({\tilde{\phi}
(x_\alpha^*x_\gamma^*a^{(n)}x_{\gamma}x_\alpha)}) \\
&\leq \sum_{p\in F}\,\,\bigg[
\sum_{\gamma\in\Lambda^{(p+n_0+q)}}\exp({(p+n_0+q)\|a\|}) \\
&\hspace{1cm} \times\sum_{\alpha\in\Lambda^{(n-q-1)}}\text{Tr }
\exp({\tilde{\phi}(x_\alpha^*q_\gamma
a^{(n-p-n_0-q)}q_\gamma x_{\alpha})})\bigg],\tag{6.5}
\end{align*}
where we have used the fact that, by the Peierls--Bogoliubov 
inequality (cf.\ Prop.\ 31.5 in \cite{OP}),
for $p\in F$ and
$|\gamma|=p+n_0+q$,
\begin{align*}
\lefteqn{\sum_{\alpha\in\Lambda^{(n-q-1)}}
\text{Tr }\exp({\tilde{\phi}(x_\alpha^*x_\gamma^*a^{(n)}x_\gamma 
x_\alpha)})}\hspace*{1cm} \\
&\leq \exp{(n_0+p+q)\|a\|}
\sum_{\alpha\in\Lambda^{(n-q-1)}}\text{Tr
}\exp({\tilde{\phi}(x_\alpha^*q_\gamma
a^{(n-p-n_0-q)}q_\gamma
x_\alpha)}).
\end{align*}
Now $(6.5)$ is bounded by
\begin{gather*}
\sum_{p\in
F}\vartheta_{p+n_0+q}\exp({(p+n_0+q)\|a\|+(n-p-n_0-q)\epsilon)}\hspace*{1.5cm}
\\
\hspace*{1cm}\times\sum_{\alpha\in\Lambda^{(n-q-1)}}
\text{Tr }\exp({\tilde{\phi}
(x_\alpha^*a^{(n-p-n_0-q)}x_{\alpha})}) \\
\hspace*{1cm}\leq\sum_{p\in
F}\vartheta_{p+n_0+q}\exp({2(p+n_0+q)\|a\|+(n-p-n_0-q)\epsilon}) \\
\hspace*{1.5cm}\times\sum_{\alpha\in\Lambda^{(n-q-1)}}\text{Tr }
\exp({\tilde{\phi}
(x_\alpha^*a^{(n-q)}x_{\alpha})}).
\end{gather*}
This shows that
$$P_\theta(\O_X,a,\Omega(n_0),11\vartheta_{n_0}\delta^{\frac12})\leq
P_{X}(\D,a)+\epsilon,$$
and so, by the arbitrarity of $\epsilon$, $P_\theta(\O_X,a)\leq
P_X(\D,a)$.
For general $a$ we write $a=a_+-\lambda I$ with
$a_+$ positive and $\lambda\in{\Bbb R}^+$.
By Lemma 6.6 we have
$$P_\theta(a)=P_\theta(a_+)-\lambda\leq
P_X(\D, a_+)-\lambda\leq P_X(\D, a),$$
and the proof is complete.
$\hfill\square$
\medskip

We next give a lower bound for $P_\theta(a)$. 
\medskip

\noindent{\bf Proposition 6.10.}
{\it Let $a$ be a positive element of $\D$ such that 
there is $p\in{\Bbb N}$ for which
$$ [a,\ x_\alpha^*a{x_\alpha}]=0,\eqno(6.6)$$
and
$$[a,\ q_\alpha]=0,\eqno(6.7)$$
for $|\alpha|\geq p$. If $\sum_i^d q_i$ is invertible then
$$P_\theta(a)\geq
\lim_n\frac{1}{n}
\log\sum_{\alpha\in\Lambda^{(n-1)}}e^{\|x_\alpha^*a{x_\alpha}\|}.$$
}\medskip

\noindent{\it Proof.} 
Suppose $b\in\D$ satisfies $(6.6)$ and $(6.7)$ for $|\alpha|\geq r$.
Consider the $C^*$--subalgebra $\C(b,r)$ generated by
$$\{x_{\alpha}C^*(b,I)x_\alpha^*,\quad 
|\alpha|=nr,
n=0,1,2,\dots\}.$$ 
Notice that $\C(b,r)$ is $\theta^r$--invariant.
For fixed $\alpha$ with $|\alpha|=nr$, 
 $x_{\alpha}C^*(b,I)x_\alpha^*$ 
is a commutative $C^*$--algebra. Furthermore, if $|\alpha|=hr$ and
 $|\beta|=kr$ with $h<k$, and $s,t\in\C^*(b,I)$, then
$x_\alpha sx_\alpha^*{x_\beta}tx_\beta^*$
is nonzero only if  $\beta=\alpha\beta'$ for some $\beta'$ of length
$(k-h)r\geq r$, and in this case
\begin{gather*}
x_\alpha sx_\alpha^*{x_\beta}tx_\beta^*=
{x_\alpha} s q_\alpha
x_{\beta'}tx_\beta^*={x_\beta}x_{\beta'}^*s{x_{\beta'}}q_\beta
tx_\beta^* \hspace*{1.5cm}\\
\hspace*{3cm} =x_\beta t q_\beta x_{\beta'}^*sx_{\beta'}x_\beta^*=
x_\beta tx_\beta^*x_\alpha s x_\alpha^*,
\end{gather*}
so that $\C(b,r)$ is commutative.
 Let $T_r$ denote the epimorphism
of the spectrum of $\C(b,r)$ obtained transposing
$\theta^r$.
Consider the open (and closed) cover $\U$ of the spectrum
of $\C(b,r)$ defined by the characteristic functions
$\{x_\alpha x_\alpha^*, |\alpha|=r\}.$
Then by the monotonicity of pressure (Prop.\ 3.3) and the fact
that the noncommutative pressure reduces to the classical
pressure on commutative $C^*$--algebras, we obtain
$$P_{\theta^r}(b)\geq p_{T_r}(b)\geq
\lim_n\frac{1}{n}\sum_{\alpha\in\Lambda^{(rn-r)}}
e^{\|x_\alpha^*(b+\theta^r(b)+\dots+\theta^{r(n-1)}(b))x_\alpha\|}.$$
Suppose $a\in\D$ satisfies $(6.6)$, $(6.7)$ for $|\alpha|\geq p$
and set $a_r=a+\theta(a)+\dots+\theta^{r-p}(a)$ for $r\geq p$.
Then $a_r$ satisfies the corresponding relations for $|\alpha|\geq r$. 
Therefore, for $r\geq p$,
$$P_{\theta^r}(a+\dots+\theta^{r-p}(a))\geq
\lim_n\frac{1}{n}\log\sum_{\alpha\in\Lambda^{(rn-r)}}
e^{\|x_\alpha^*a_r^{(n)}x_\alpha\|}$$
where
\begin{gather*}
a_r^{(n)}=(a+\dots+\theta^{r-p}(a))+
(\theta^r(a)+\dots+\theta^{2r-p}(a))\hspace*{2cm} \\
\hspace*{3cm}\ +\dots+(\theta^{r(n-1)}(a)+\dots+
\theta^{rn-p}(a)).
\end{gather*}
Now
$$\|x_\alpha^*a_r^{(n)}x_\alpha\|\geq
\|x_\alpha^*(a+\dots+\theta^{rn-1}(a))x_\alpha\|-n(p-1)\|a\|$$
and, by the monotonicity of pressure with respect to the 
self-adjoint element and scalar
additivity (Prop.\ 3.1),
\begin{align*}
P_{\theta^r}(a+\dots+\theta^{r-p}(a))&=
P_{\theta^r}(a^{(r)}-
(\theta^{r-p+1}(a)+\dots+\theta
^{r-1}(a))) \\
&\leq rP_\theta(a)+(p-1)\|a\| 
\end{align*}
and so
$$P_\theta(a)+\frac{p-1}{r}\|a\|\geq
\lim_n\frac{1}{n}\log\sum_{\alpha\in\Lambda^{(n-r)}}
e^{\|x_\alpha^*a^{(n)}x_\alpha\|}-\frac{p-1}{r}\|a\|.$$
Since 
$$\lim_n\frac{1}{n}\log\sum_{\alpha\in\Lambda^{(n-r)}}
e^{\|x_\alpha^*a^{(n)}x_\alpha\|}=
\lim_n\frac{1}{n}\log\sum_{\alpha\in\Lambda^{(n-1)}}
e^{\|x_\alpha^*a^{(n)}x_\alpha\|},$$
for all $r$, we obtain the result letting $r\to\infty$.
$\hfill\square$\medskip

\noindent{\it Proof of Theorem 6.2.}
Combining Prop.\ 6.9 and Prop.\ 6.10, we obtain the proof of
the
first
part of the Theorem. 
Assume that $a\in\C(\Lambda)$. Then, if 
$a\geq0$, $\|x_\alpha^*a^{(n)}x_\alpha\|$ is the supremum of
$a^{(n)}$ in the cylinder set
$\{(i_j)\in\Lambda: (i_1,\dots,i_{|\alpha|})=\alpha\}$, thus the
pressure formula reduces to
the 
classical pressure formula for a positive continuous function on $\Lambda$
(see, e.g., \cite{DGS}).
For the last part we will adapt  from \cite{NS} 
a proof of the variational principle for a class of asymptotically Abelian
$C^*$--algebras. By the additivity of pressure under addition of
scalars, we may assume $a\geq0$. Consider the unital, commutative
$\theta^r$--invariant $C^*$--algebra
$\C(a_r, r)$ introduced in the proof of Prop.\ 6.10.
By the classical variational principle, given $\epsilon>0$ there exists
a $\theta^r$--invariant state $\mu_r$ on $\C(a_r, r)$ such that
$$h_{\mu_r}(\theta^r\upharpoonright_{\C(a_r,r)})+\mu_r(a_r)>
P_{\theta^r\upharpoonright_{\C(a_r,r)}}-\epsilon.$$
By Prop.\ 4.17,   
$\mu_r$ extends to a $\theta^r$--invariant state $\tilde{\sigma}_r$
on $\O_X$ in such a way that
$$h_{\tilde{\sigma}_r}(\theta^r)>
h_{\mu_r}(\theta^r\upharpoonright_{\C(a_r,r)})-1.$$
Then $\sigma_r:=\frac{1}{r}\sum_0^{r-1}\tilde{\sigma}_r\theta^j$ is 
$\theta$--invariant. By Prop.\ 3.3 in \cite{ST} 
$h_{\sigma_r}(\theta)=\frac{1}{r}h_{\sigma_r}(\theta^r)$. By 
concavity
of the Sauvageot--Thouvenot entropy
(Prop.\ 4.18) and
Lemma 4.19, one has
\begin{align*}
h_{\sigma_r}(\theta)=\frac{1}{r}h_{\sigma_r}(\theta^r)
&\geq \frac{1}{r^2}\sum_0^{r-1}h_{\tilde{\sigma}_r\theta^j}
(\theta^r)-\frac{\log r}{r} \\
&\geq\frac{1}{r}h_{\tilde{\sigma}_r}(\theta^r)-\frac{\log r}{r} \\
&> \frac{1}{r}
h_{\mu_r}(\theta^r\upharpoonright_{\C(a_r,r)})-\frac{1}{r}-\frac{\log
r}{r}.
\end{align*}
Since
$$\sigma_r(a)\geq\frac{1}{r}\mu_r(a_r)-\frac{p-1}{r}\|a\|,$$
we infer that
$$\limsup_r
h_{\sigma_r}(\theta)+\sigma_r(a)\geq
\limsup_r\frac{1}{r}(P_{\theta^r\upharpoonright_{\C(a_r,r)}}(a_r)-\epsilon)
=
P_\theta(a).$$
The last equality has been proven in Prop.\ 6.10, taking into
account Prop.\ 6.9.
$\hspace*{\fill}\square$\medskip

\noindent{\it 2.  Equilibrium states}\medskip

We start this subsection proving that, at least if we restrict further
the space of potentials, equilibrium states exist for $\O_X$.
This generalizes Theorem 5.4 to the algebras $\O_X$.
\medskip

\noindent{\bf Proposition 6.11.} {\it Assume that $ht(\A,X)=0$ and that
$a$ is a self-adjoint element of the commutative $C^*$--subalgebra 
$\C(\Lambda)\subset\O_X$, so, by Theorem 6.2,
 $P_\theta(a)=p_T(a)$. Then any faithful equilibrium 
measure for $(\Lambda,T,a)$ extends to an equilibrium state
for $(\O_X,\theta,a)$ for the CNT entropy, and thus also for the
Sauvageot--Thouvenot entropy and  
the local state approximation entropy.}\medskip

\noindent{\it Proof.} By Prop.\ 8.4 and Lemma 8.3 in \cite{PWY} any
faithful shift-invariant
measure $\mu$ on $\Lambda$ extends to a $\sigma$--invariant state
of $\O_X$ such that ${h^{\text{CNT}}}_\sigma\geq h_\mu(T).$
Therefore if we start with an equilibrium measure for
the Kolmogorov--Sinai entropy, by Prop.\ 4.10
and the fact that the Sauvageot--Thouvenot entropy majorizes the CNT
entropy we obtain
\begin{gather*}
{hm}_\sigma(\theta)+\sigma(a)\geq
h^{\text{ST}}_\sigma(\theta)+\sigma(a)\geq\hspace*{5.5cm} \\
\hspace*{3.5cm}h^{\text{CNT}}_\sigma(\theta)+\sigma(a)\geq h_\mu(T)+\mu(a)=
p_T(a)=P_\theta(a).
\end{gather*}
$\hfill\square$\medskip

We next give an upper
bound for the local state approximation entropy of $\theta$
which is similar to the corresponding bound for the pressure
(Prop.\ 6.9). This bound, together with Prop.\ 4.6, will lead, in a similar
way,
to a computation of $hm_\sigma(\theta)$, and will also be
useful
when discussing equilibrium states.\medskip

\noindent{\bf Proposition 6.12.} {\it Let $X$ be a Hilbert bimodule over
a unital exact $C^*$--algebra $\A$ satisfying the same assumptions as
in Prop.\ 6.9,  let $\sigma$  be a $\theta$--invariant
state of $\O_X$ and $m$ the
probability measure on $\Lambda$
obtained restricting $\sigma$ to $\C(\Lambda)$.
Then
$$hm_\sigma(\theta)\leq h_m(T)+ht(\A,X)$$
where $h_m(T)$ denotes the Kolmogorov--Sinai entropy
of the shift $T$ of $\Lambda$.}\medskip

\noindent{\it Proof.} The proof parallels that of Prop.\ 6.9,
with the same local approximations being employed to obtain an upper
bound for $hm_\sigma(\theta)$. 
Thus, as in Prop.\ 6.9, given a finite subset $\Delta\subset\O_X^{(0)}$,
$\delta>0$, and $n_0\in{\Bbb N}$, let $F=F(\delta,n_0)\subset{\Bbb N}$ be
a finite subset independent of $\Delta$ such that
for any $(\phi,\psi,\B)\in\text{CPA}(\O_X^{(0)}, \Delta^{(\max F,X)},
\frac{11\delta^{\frac12}}{2\max_F\theta})$ there is a triple
$$(\phi',\psi',\B')
\in\text{CPA}_0(\O_X,\cup_{|\gamma|\leq n_0}\Delta
x_\gamma,
11\delta^{\frac12})$$ where
$\B'=M_{\sum_F\vartheta_p}\otimes \B$. 
Here we shall need recall,
from Lemma 7.5 in \cite{PWY}, that $\psi'$ is of the form
$\psi'=\tilde{S}_{F,f}\circ(\iota\otimes\psi):\B'\to\B(\H)$
where $f\in\ell^2({\Bbb N})$ has support in $F$,
$\|f\|_2\leq1$, and 
$$\tilde{S}_{F,f}:t=(t_{\alpha,\beta})\in
 M_{\sum_F\vartheta_p}(\B(\H))\to
\sum_{|\alpha|,|\beta|\in F}f(|\alpha|)\ov{f(|\beta|)}
x_\alpha t_{\alpha,\beta}x_\beta^*\in\B(\H).$$

Let $\Omega\subset\A$ be a finite subset containing $I$.
Pick
$$(\tilde{\phi},\psi,M_N({\Bbb C}))\in 
\text{CPA}\Big(\A,\Omega^{(n+n_0+\max
F,X)},\frac{\delta}{4{\max^2}_{p\in F}\vartheta_p}\Big)$$ of minimal rank
and 
follow the
same procedure 
as in the proof of Prop.\ 6.9 
to obtain a triple
$$(\phi'',\psi''
,M_{\vartheta_{n+n_0-1}}\otimes M_{\sum_{p\in
F}\vartheta_p}\otimes M_N)\in
\text{CPA}_0(\O_X,\Omega(n_0)^{(n)},
11\vartheta_{n_0}\delta^{\frac12})$$
where
$\psi''=\psi_{n+n_0-1}\circ(\iota_{M_{n+n_0-1}}\otimes\psi')$.
Pick
$\omega\in\mathfrak{E}(\sigma,\iota)$. 
The positive linear functional $\omega\circ\psi''$ is determined by
it values on each matrix unit, which are given by
 $$\omega\circ\psi''(e_{\alpha,\beta}\otimes
e_{\mu,\nu}\otimes
e_{i,j})=f(|\mu|)\ov{f(|\nu|)}\omega(x_{\alpha\mu}\psi(e_{i,
j}){x_{\beta\nu}^*}).$$
We have that
$$ S(\omega\circ\psi'')
\leq
S(\text{diag}(|f(|\mu|)|^2\omega(x_{\alpha\mu}\psi(e_{i,i})
x_{\alpha\mu}^*)_{|\alpha|=n+n_0-1,|\mu|\in F, i=1,\dots,N}))$$
by the estimate on page 60 in \cite{OP}, e.g., and this last expression 
is bounded above by
$$S\Big(\text{diag}\Big(\frac{1}{N}|f(|\mu|)|^2m(p_{\alpha\mu})
\Big)_{\alpha,\mu,i}\,\Big),$$
which in turn is bounded above by 
\begin{align*} 
\lefteqn{-\sum_{\alpha,\mu}
|f(|\mu|)|^2m(p_{\alpha\mu})
\log\Big(\frac{1}{N}|f(|\mu|)|^2m(p_{\alpha\mu})\Big)}\hspace*{1.4cm} \\
&= \sum_{\alpha,\mu}-|f(|\mu|)|^2m(p_{\alpha\mu})
\log(m(p_{\alpha\mu}))+\log N 
\sum_\mu|f(|\mu|)|^2m(p_\mu) \\
&\hspace*{4cm} \ -\sum_\mu|f(|\mu|)|^2m(p_\mu)\log m(p_\mu),
\end{align*}
with the equality following from the $T$--invariance of $m$. Finally,
using the equality
$$\sum_\mu |f(|\mu|)|^2m(p_{\alpha\mu})={\|f\|_2}^2m(p_\alpha)$$
and the concavity of $x\mapsto-x\log x$, we see that the last 
displayed expression is bounded by
\begin{align*}
\lefteqn{-\sum_{|\alpha|=n+n_0-1}{\|f\|_2}^2m(p_\alpha)
\log({\|f\|_2}^2m(p_\alpha))} \\
&\hspace*{1.5cm}\ + \log 
N\sum_{|\mu|\in F}|f(|\mu|)|^2m(p_\mu)
-\sum_{|\mu|\in F}|f(|\mu|)|^2m(p_\mu)\log m(p_\mu).
\end{align*}
Therefore, since $\|f\|_2\leq1$, Prop.\ 4.5 yields
$$hm_\sigma(\theta,\iota,\omega,\Omega(n_0))\leq
h_\mu(T)+ht(\A,X).$$ 
$\hfill\square$\medskip

We next derive a few consequences on equilibrium states from the previous
proposition.
There is a natural conditional expectation $E:\O_X\to\D_0$
where $\D_0$ is the $C^*$--subalgebra of $\D$ generated
by elements of the form $x_\alpha ax_\alpha^*$, $a\in\A$,
$\alpha\in\cup_n\Lambda^{(n)}$, defined in the 
following way. Compose the average over the gauge action
$\O_X\to\O_X^{(0)}$ with the pointwise norm limit $P:\O_X\to\O_X^{(0)}$ of
the maps
$$t\to P_n(t)=\sum_{|\alpha|=n}
x_\alpha x_\alpha^*tx_\alpha x_\alpha^*.$$
One has
$$E\circ\theta=\theta\circ E.$$
\medskip

\noindent{\bf Corollary 6.13.} 
{\it Let $(\O_X,\theta)$ be the
$C^*$--dynamical
system constructed as in Theorem 6.2. Assume that
$ht(\A,X)=0$.
 If $\sigma$ is a $\theta$--invariant
state and the restriction $m$ of $\sigma$ to $\C(\Lambda)$ is a 
faithful measure, then $\sigma\circ E$ is  a $\theta$--invariant
state centralized by $\C(\Lambda)$ for which
$$h_m(T)=hm_{\sigma\circ E}(\theta)=h_{\sigma\circ E}^{\text{ST}}(\theta)=
h_{\sigma\circ E}^{\text{CNT}}(\theta),$$
where $h^{\text{ST}}$ and $h^{\text{CNT}}$ denote respectively
the Sauvageot--Thouvenot and CNT entropy.
If moreover $a\in\D_0$ is a self--adjoint element and $\sigma$ is an
equilibrium
state for $(\O_X,\theta,a)$, then 
$\sigma\circ E$ is an equilibrium state for the same system
(both with respect to $hm$).
}\medskip

\noindent{\it Proof.} Note that 
$\sigma\circ E$ is a $\theta$--invariant state since $E$
commutes with $\theta$. Furthermore $\sigma\circ E$ is centralized
by $\C(\Lambda)$ since $E$ is a conditional expectation onto $\D_0$,
which contains $\C(\Lambda)$, and $\C(\Lambda)$ commutes with $\D_0$.
Therefore Props.\ 8.2 and 8.3 in \cite{PWY} can be applied to
$\sigma\circ E$.
Using, Prop.\ 6.12,
Prop.\ 4.10, the fact that Sauvageot--Thouvenot entropy majorizes
the CNT entropy \cite{ST}, and Props.\ 8.2 and 8.3 of 
\cite{PWY}, respectively, we infer that
$$h_m(T)\geq hm_{\sigma\circ E}(\theta)\geq
h_{\sigma\circ E}^{\text{ST}}(\theta)\geq
h_{\sigma\circ E}^{\text{CNT}}(\theta)\geq
h_m(T).$$
Assume now that $\sigma$ is an equilibrium state for $(\O_X,\theta,a)$.
Then
$$hm_{\sigma\circ E}(\theta)+\sigma(a)=h_m(T)+\sigma(a)\geq
hm_\sigma(\theta)+\sigma(a)=P_\theta(a),$$
and so $\sigma\circ E$ must be an equilibrium state for the same system 
as well. $\hfill\square$\medskip

The following is a converse of Prop.\ 6.11.\medskip

\noindent{\bf Corollary 6.14.} {\it Let $(\O_X,\theta)$ be the
$C^*$--dynamical
system constructed as in Theorem 6.2. Assume that $ht(\A,X)=0$ and let
$a$ be a selfadjoint element of the canonical Abelian subalgebra
$\C(\Lambda)$. Let $H_\sigma(\theta)$ be 
either the local state approximation entropy,
the Sauvageot--Thouvenot entropy, or the CNT entropy.
If $\sigma$ is a
$\theta$--invariant equilibrium state
for $(\O_X,\theta,a)$ w.r.t. $H_\sigma(\theta)$,
then the measure $m$ obtained restricting
$\sigma$ to $\C(\Lambda)$ is an equilibrium measure for
$(\Lambda,T,a)$. Furthermore one has
$$H_\sigma(\theta)=h_m(T)$$
where $h_m$ is the Kolmogorov--Sinai entropy.
}\medskip

\noindent{\it Proof.} By the comparison between the various state--based
entropies (Prop.\ 4.10),
 the fact that $P_\theta(a)$ coincides with the classical
pressure of $a$ (Prop.\ 6.2),
Prop.\ 3.6,
and Prop.\ 6.12 under the
assumption $ht(\A,X)=0$, we have
$$p_{T}(a)=
P_\theta(a)=H_\sigma(\theta)\leq hm_\sigma(\theta)+\sigma(a)\leq
h_m(T)+m(a),$$
so that $m$ is an equilibrium state for $(\Lambda,T,a)$. Since all the
inequalities become equalities, we conclude that
$H_\sigma(\theta)=h_m(T)$.$\hfill\square$
\medskip

\noindent{\it 3.\  An application to 
Matsumoto algebras associated to subshifts}\medskip

We conclude this section with  an application  to Cuntz--Krieger
algebras,
or, more generally, to Matsumoto $C^*$--algebras.
\medskip

\noindent{\bf Corollary 6.15.} {\it Let $\O_\Lambda$ denote the Matsumoto
algebra associated to a subshift of one of the following kinds:
\begin{description}
\item{\rm a)} finite type subshifts,
\item{\rm b)} sofic subshifts,
\item{\rm c)} $\beta$--shifts.
\end{description}
Then for any real-valued $f\in\C(\Lambda)\subset\O_\Lambda$,
$P_\theta(f)$ equals the classical pressure of $f$ with respect to
the shift $T$:
 $$P_\theta(f)=p_T(f).$$
Furthermore any shift-invariant measure $\mu$ on $\C(\Lambda)$
extends to a $\theta$--invariant state  $\sigma$ on $\O_\Lambda$
with the propery $h_\sigma(\theta)\geq h_\mu(T)$. In particular,
if $\mu$ is an equilibrium measure for $(\Lambda, T,f)$, the corresponding
extension is an equilibrium state for $(\O_\Lambda, \theta,f)$.}\medskip

\noindent{\it Proof.} For Matsumoto $C^*$--algebras 
the coefficient algebra $\A$ is commutative and commutes with
$\C(\Lambda)$. Furthermore,
the growth of the local
completely positive $\delta$--ranks $rcp(\A, \Omega^{(n,X)},\delta)$ is
polynomial (see\cite{PWY}), and
so $ht(\A, X)=0$, implying the first part of the assertion.
Let $\mu$ be a $T$--invariant measure on $\Lambda$.
That $\mu$ extends to a $\theta$--invariant state $\sigma$ on $\O_\Lambda$
with entropy as least as large has been proven in Theorem 8.6 of \cite{PWY}.
The rest is now clear.
$\hfill\square$\medskip

\section{The KMS condition and equilibrium in $\O_A$}

We will show how certain
equilibrium states for Cuntz--Krieger algebras can be constructed
from KMS states
with respect to a suitable one-parameter
automorphism group in the case where the self-adjoint element has small
variation on the underlying subshift of finite type and is a H\"older
continuous function. 
The key idea is to establish a connection between
KMS states with respect to this group and the 
Perron--Frobenius--Ruelle theorem for subshifts
of finite type \cite{R68, Bowen, W75}.

Let $A$ be  $\{0,1\}$-matrix with no row or column identically zero,
and let
$a\in\C(\Lambda_A)$  be a self-adjoint element. Consider for
$\beta\in{\Bbb R}$ the
unitary group of $\C(\Lambda_A)$, $U_{\beta,a}(t)=\exp(it(\beta-a))$, 
and define the one-parameter automorphism group of $\O_A$ 
$${\alpha^{\beta,a}}_t(s_i)=U_{\beta,a}(t)s_i,\quad i=1,\dots,d,$$
where $s_1,\ldots,s_d$ are the generating partial isometries of $\O_A$
\cite{Watatani}.
We shall also need a positive operator of $\C(\Lambda_A)$
whose spectral properties and their relation with 
equilibrium states
were first studied by Ruelle \cite{R68}
in the case of the full $2$--shift and
in a more general setting by
Bowen \cite{Bowen} and Walters \cite{W75}.
Set
$$\L_a(f)(x)=\sum_{i:A_{i x_1}=1}\exp(a(ix))f(ix),\quad
f\in\C(\Lambda_A),$$
where $x=(x_k)_k\in\Lambda_A$.
Notice that we can write, in $\O_A$,
$$\L_a(f)=\sum_is_i^*e^afs_i.$$
Thus $\L_a$ extends in an obvious way to an operator on
$\O_A$, which we will denote by $\ov{\L}_a$.
We begin by establishing some partial results.
\medskip

\noindent{\bf Lemma 7.1.} {\it $\C(\Lambda_A)$ is contained in the 
algebra of fixed points under $\alpha^{\beta,a}$ for all $\beta\in{\Bbb
R}$ and all self-adjoint $a\in\C(\Lambda_A)$.}\medskip

\noindent{\it Proof.}  For all $j\in{\Bbb N}$ 
and $t\in{\Bbb R}$ we have $\theta^j(U_{\beta,a}(t))\in\C(\Lambda_A)$, and
therefore $[\theta^j(U_{\beta,a}(t)),\C(\Lambda_A)]=0$.
On the other hand, if $f\in\C(\Lambda_A)$ is of the 
form $f=s_{i_1}\dots s_{i_r}(s_{i_1}\dots s_{i_r})^*$, then
for all $t\in{\Bbb R}$ 
$${\alpha^{\beta,a}}_t(f)=U_{\beta,a}(t)\dots\theta^{r-1}(U_{\beta,a}(t))f
\theta^{r-1}(U_{\beta,a}(-t))\dots
U_{\beta,a}(-t)=f$$
since $\theta$ is multiplicative on $\C(\Lambda_A)$ (in fact, 
one can easily check that
$\theta$ is
multiplicative on the relative commutant of $\{s_1^*s_1,\dots,
s_d^*s_d\}$ in
$\O_A$), completing the
the proof.\hfill$\square$\medskip

The following result is well-known. We supply 
a proof for convenience.\medskip

\noindent{\bf Lemma 7.2.} {\it If $\mu$ is a state on
$\C(\Lambda_A)$ such that, for some $\beta\in{\Bbb R}$,
$$e^\beta\mu(f)=\mu(\L_a(f)), \quad f\in\C(\Lambda_A),$$
then $$\min(a)+\log r(A)\leq\beta\leq\max(a)+\log r(A).$$}\medskip

\noindent{\it Proof.} For $n\in{\Bbb N}$ and
$x=(x_k)_k\in\Lambda_A$,
$$\L_a^n(1)(x)=\sum_{A_{i_1 i_2}=\dots=
A_{i_n x_1}=1}e^{a^{(n)}(i_1,\dots,i_n,x)}.$$
Therefore
$$e^{n\min a}\vartheta_n\leq\L_a^n(1)\leq e^{n\max a}\vartheta_n,$$
where, as usual, $\vartheta_n$ denotes the cardinality of the set of words
of length $n$ appearing in $\Lambda_A$.
Applying $\mu$ yields
$$e^{n\min a}\vartheta_n\leq e^{n\beta}\leq e^{n\max a}\vartheta_n,$$
and so computing $\lim_n\frac{1}{n}\log (\cdot )$
and using the fact that $\lim_n\frac{1}{n}\log\vartheta_n=\log r(A)$ (see,
e.g., \cite{DGS}) we
obtain the desired estimate.\hfill$\square$\medskip

We next describe a bijective correspondence between KMS states for
$\alpha^{\beta,a}$ and positive eigenvectors of the Banach space adjoint
of the Ruelle operator
$$\L_a^*:\C(\Lambda_A)^*\to\C(\Lambda_A)^*.$$
 We start showing that $(\alpha^{\beta,a},1)$--KMS 
states restrict to positive eigenvectors of $\L_a^*$.
\medskip

\noindent{\bf Lemma 7.3.} {\it If $\omega$ is a 
$(\alpha^{\beta,a},1)$--KMS state on $\O_A$ then
$$\sum_i\omega(s_i^*e^abs_i)=e^\beta\omega(b),\quad b\in\O_A.$$
In particular, if $\mu:=\omega\upharpoonright\C(\Lambda_A)$ then
$$\L_a^*(
\mu)=
e^\beta\mu.$$}\medskip

\noindent{\it Proof.} 
We first note that, for $i=1,\dots,d$,
${\alpha^{\beta,a}}_{-i}(s_j^*)=s_j^*e^{-\beta+a}$. Thus,
by the KMS property,
if $b\in\O_A$ then
$$\omega(b)=\sum_j\omega(bs_js_j^*)=\sum_j\omega(\alpha^{\beta,a}_{-i}
(s_j^*)bs_j)=e^{-\beta}\sum_j\omega(s_j^*e^abs_j).$$
\hfill$\square$\medskip

\noindent{\bf Lemma 7.4.} {\it If $A$ is aperiodic then $r(A)>1$.}\medskip

\noindent{\it Proof.} Let $N$ be a positive integer such that
all of the entries of $A^N$ are positive. Since $A\in M_d(\{0,1\})$,
$A^N_{ij}\geq 1$ for all $i,j$, and therefore $A^{pN}_{ij}\geq d^{p-1}$
for $p\in{\Bbb N}$, so that
$$r(A)=\lim_n\|A^n\|^{1/n}=\lim_p\|A^{Np}\|^{1/Np}\geq d^{1/N}>1.$$
\hfill$\square$\medskip

We define a metric on $\Lambda_A$ by $d(x,y)=\frac{1}{k}$ where $k$
is the least integer
for which $x_k\neq y_k$. For
$f\in\C(\Lambda_A)$, we set
\begin{align*}
\text{var}_0(f) &= \max f-\min f, \\
\text{var}_n(f) &= \sup\Big\{|f(x)-f(y)|, d(x,y)\leq \frac{1}{n+1}\Big\},
\quad n\in{\Bbb N}.
\end{align*}
Note that
contunuity implies that $\text{var}_n(f)\to0$ as $n\to\infty$.
\medskip

\noindent{\bf Lemma 7.5.} {\it Let $A$ be an aperiodic $\{0,1\}$--matrix
and $a\in\C(\Lambda_A)$ a self-adjoint element such that
$\text{\rm var}_0(a)<\log r(A)$.
Then any
 $(\alpha^{\beta,a},1)$--KMS
state
on $\O_A$  is
gauge-invariant 
and faithful.}\medskip

\noindent{\bf Proof.} We first establish gauge invariance. 
Let $\omega$ be a $(\alpha^{\beta,a},1)$--KMS state.
We
need to show that if $b\in\O_A^{(k)}$ for some 
$k>0$, then $\omega(b)=0$. It suffices to pick $b$ nonzero and of the form
$b=s_{i_1}\dots s_{i_s}(s_{j_1}\dots s_{i_r})^*$ with $s-r=k$.
Since $\omega$ is an  $\alpha^{\beta,a}$--KMS state, it is
$\alpha^{\beta,a}$--invariant, and therefore
$$\omega({\alpha^{\beta,a}}_{it}(b))=\omega(b),\quad t\in{\Bbb R},$$
that is,
$$e^{-k\beta t}\omega(e^{ta}\dots\theta^{s-1}(e^{ta})b
\theta^{r-1}(e^{-ta})\dots e^{-ta})=\omega(b).$$
Now by Lemma 7.1 $\theta^j(e^{ta})$ lies in the centralizer of $\omega$,
and thus the above equality reduces to
$$e^{-k\beta t}\omega(\theta^r(e^{ta})\dots\theta^{s-1}(e^{ta})b)=
\omega(b),$$
and so
$|\omega(b)|\leq e^{-k\beta t}\|e^{ta}\|^k\|b\|$.
Hence, for $t\geq0$,
$$\frac{|\omega(b)|}{\|b\|}\leq e^{kt(-\beta+\max a)},$$
which implies, assuming $\omega(b)\neq0$,
that $\beta\leq\max a$.
On the other hand, by Lemma 7.3 the restriction of $\omega$
to $\C(\Lambda_a)$ is a positive eigenvector of the transposed Ruelle 
operator with eigenvalue $e^\beta$,
and so $\beta\geq\min a+\log r(A)$ by Lemma 7.2. Therefore we 
must have $\log r(A)\leq\text{var}_0(a)$, which contradicts our assumption.

We next show that $\omega$ is faithful. By the previous part,
it suffices to show that the restriction of $\omega$ to $\O_A^{(0)}$
is faithful. Set $\I=\{b\in\O_A^{(0)}: \omega(b^*b)=0\}$.
Clearly $\I$ is a closed left ideal of $\O_A^{(0)}$. If $c$ ranges over 
a dense set of analytic vectors for
$\alpha^{\beta,a}$ and $b\in\I$, then by the KMS property
$$\omega(c^*b^*bc)=\omega(\alpha^{\beta,a}_{-i}(c)c^*b^*b)=0,$$
so that $\I$ is a two-sided closed ideal of $\O_A^{(0)}$. Since
$A$ is aperiodic, $\O_A^{(0)}$ is a simple AF $C^*$--algebra, and so
$\I=0$.
\hfill$\square$\medskip

\noindent{\bf Proposition 7.6.} {\it Let $A$ be an aperiodic 
$\{0,1\}$--matrix.
If $a\in\C(\Lambda_A)$ is a self-adjoint element 
the map
$$\omega\to\mu:=\omega\upharpoonright\C(\Lambda_A)$$
sets up a surjective correspondence between the set of
$(\alpha^{\beta,a}, 1)$--KMS states of $\O_A$ and the set of
probability measures  on $\Lambda_A$ for which 
$$\L_a^*\mu=e^\beta\mu.$$
If in addition $\text{\rm var}_0(a)<\log r(A)$, this map is a 
bijection.
}\medskip

\noindent{\it Proof.} Let $\omega$ be
a $(\alpha^{\beta,a}, 1)$--KMS state. By Lemma 7.3 the measure $\mu$ 
corresponding
to
$\omega\upharpoonright\C(\Lambda_A)$ is an eigenvector of 
the transposed Ruelle operator
with eigenvalue $e^\beta$.

We show that any state $\mu$ on $\C(\Lambda_A)$ arising as
an eigenvector of the transposed Ruelle operator with eigenvalue $e^\beta$
 is the restriction of
a (gauge--invariant) $(\alpha^{\beta,a},1)$--KMS state.
Consider the  $C^*$--subalgebras $F_n$ of $\O_A$
  linearly spanned by elements of the form 
$s_{i_1}\ldots s_{i_n}\C(\Lambda_A)(s_{j_1}\ldots s_{j_n})^*$.
Note that $F_n\subset F_{n-1}$ and that $\cup_n F_n$ is dense in
$\O_A^{(0)}$.
Recursively define for each $n=0,1,2,\dots$ a state $\omega_n$ on $F_n$ 
by $\omega_0=\mu$ and
$$\omega_n(b)=e^{-\beta}\sum_i\omega_{n-1}(s_i^*e^{a/2}be^{a/2}s_i),
\quad b\in F_n.$$
Then $\omega_1$ extends $\omega_0$, and one can check that
$\omega_n$ extends
$\omega_{n-1}$ for all $n$.
Consider the gauge-invariant state $\omega$ of $\O_A$ that extends 
$\omega_n$ on $F_n$. By construction, $\omega$ satisfies the scaling property
$$\sum_i\omega(s_i^*e^abs_i)=e^\beta\omega(b),\quad b\in\O_A.$$
We show that $\omega$ is a $(\alpha^{\beta,a}, 1)$--KMS state.
Consider elements of the form 
$b=s_{i_1}\dots s_{i_s}(s_{j_1}\dots s_{j_r})^*$,
$c=s_{h_1}\dots s_{h_r}(s_{k_1}\dots s_{k_s})^*$.
We need to show that
$\omega(bc)=\omega(\alpha^{\beta,a}_{-i}(c)b)$.
If $(j_1,\dots,j_r)\neq(h_1,\dots,h_r)$ then
$\omega(bc)=0$. We also have 
$\omega(\alpha^{\beta,a}_{-i}(c)b)=0$
since the l.h.s.\ is
$$\omega(e^{\beta-a}s_{h_1}\dots e^{\beta-a}s_{h_r}
s_{k_s}^*e^{-\beta+a}\dots s_{k_1}^*e^{-\beta+a}
s_{i_1}\dots s_{i_s}(s_{j_1}\dots s_{j_r})^*),$$
which, by an iteration of the scaling property, is seen to be zero.
Assume now that
$(j_1,\dots,j_r)=(h_1,\dots,h_r)$. The scaling property 
tells us again that if $\omega(bc)$ is 
nonzero, we must have $(i_1,\dots,i_s)=(k_1,\dots,k_s)$. Clearly,
the above computation shows that if 
$\omega(\alpha^{\beta,a}_{-i}(c)b)\neq0$
the same condition holds.
Again, by the scaling property,
the latter is 
\begin{align*}
\lefteqn{\omega(
s_{i_s}^*e^{-\beta+a}\cdots s_{i_1}^*e^{-\beta+a}s_{i_1}\cdots 
s_{i_s}(s_{j_1}\cdots s_{j_r})^*s_{j_1}\cdots s_{j_r})} \hspace*{2.5cm} \\ 
&= \omega (s_{i_1}\cdots s_{i_s}(s_{j_1}\cdots s_{j_r})^*s_{j_1}\cdots 
s_{j_r}(s_{i_1}\cdots s_{i_s})^*) \\ 
&= \omega (bc).
\end{align*}

We next show that the map $\omega\to\mu$ is one-to-one if
$\text{var}_0(a)<\log r(A)$.
Since $\omega$ is gauge-invariant by Lemma 7.5, it is determined by its
restriction to
$\O_A^{(0)}$.
Applying ${\ov \L}_a^n$ on words of the form
$s_{i_1}\dots s_{i_n}(s_{j_1}\dots s_{j_n})^*$ 
yields an element of
$\C(\Lambda_A)$, and so again by Lemma 7.3
the restriction of $\omega$ to $\O_A^{(0)}$
is determined uniquely by its values on $\C(\Lambda_A)$,
and the proof is complete.\hfill$\square$\medskip

We recall the Perron--Frobenius--Ruelle theorem for subshifts of finite
type.
\medskip

\noindent{\bf Theorem 7.7.} { \cite{R68, Bowen, W75}
 \it 
Let $A$ be an aperiodic $\{0,1\}$--matrix and $a\in\C(\Lambda_A)$ 
a self-adjoint
element satisfying 
$$\sum_n\text{\rm var}_n(a)<\infty.$$ Then 
\begin{description}
\item{\rm (a)} $\L_a$ admits a strictly positive eigenvector
$h$ which is unique up to a scalar factor,
\item{\rm (b)} $\L_a^*$  admits a unique probability measure eigenvector
$\mu$,
\item{\rm (c)}
one has
${\L_a}h=\lambda h$ and $\L_a^*\mu=\lambda\mu$, where
$\log\lambda=p_T(a)=\log r(\L_a)$, with $r(\L_a)$  the
spectral radius of $\L_a$,
\item{\rm (d)} 
$$\frac{\L_a^n(f)}{\lambda^n}\to\frac{\mu(f)}{\mu(h)}h$$
uniformly for all $f\in\C(\Lambda_A)$, 
\item{\rm (e)} $\nu(f):=\mu(hf)$, $f\in\C(\Lambda_A)$ is the unique 
equilibrium measure for $(\Lambda_A, T, a)$, 
\item{\rm (f)} $\mu$ an $\nu$ are faithful.
\end{description}}

We are now in a position to establish a connection between 
$(\alpha^{\beta,a}, 1)$--KMS states and equilibrium measures
for $\O_A$.
\medskip

\noindent{\bf Theorem 7.8.} {\it Assume that $A$ is aperiodic
and that the self-adjoint element $a$ satisfies
\begin{gather*}
\sum_n\text{\rm var}_n(a)<\infty.
\end{gather*}
Then 
$\O_A$ admits a $(\alpha^{\beta,a},1)$--KMS state
if and only if $\beta=P_{\theta}(a)$. 
If $\text{\rm var}_0(a)<\log r(A)$ then
there is exactly one such state $\omega$.
If $h$ is the unique strictly positive eigenvector of $\L_a$ with
$\omega(h)=1$,
$\sigma(b):=\omega(hb)$ is a faithful equilibrium state for 
$(\O_A,\theta, a)$.}\medskip

\noindent{\it Proof.} 
By Prop.\ 7.6 the set of
$(\alpha^{\beta, a}, 1)$--KMS states 
corresponds surjectively to
the set of probability measures eigenvectors of $\L_a^*$ with eigenvalue
$e^\beta$,
and therefore by Theorem 7.7, there is a 
$(\alpha^{\beta, a}, 1)$--KMS state if and only if $\beta=P_\theta(a)$.
If $\text{var}_0(a)<\log r(A)$, there is exactly one such state, again by
Prop.\ 7.6 and Theorem 7.7.
Furthermore, by (e) of Theorem 7.7, the restriction 
of  $\sigma$  to
$\C(\Lambda_A)$ is the unique equilibrium measure $\nu$ for 
$(\Lambda_A, T, a)$. We note that $\sigma$ is $\theta$--invariant,
for if $b\in\O_A$ then
$$\sigma(\theta(b))=\sum_i\omega(hs_i bs_i^*)=e^{-P_\theta(a)}
\sum_i\omega(s_i^*e^a hs_ib)=\omega(hb)=\sigma(b).$$
Since
$\omega$ contains $\C(\Lambda_A)$ in its centralizer,
the same holds
for $\sigma$. Thus,
by Lemma 5.3,
$h_\mu(T)\leq h_\sigma(\theta)$, and hence
$$p_T(a)=h_\nu(T)+\nu(a)\leq h_\sigma(\theta)+\sigma(a)\leq
P_\theta(a)=p_T(a),$$
which establishes that $\sigma$ is an equilibrium state for
$(\O_A,\theta,a)$.
\hfill$\square$\bigskip

\noindent{\bf Acknowledgments.} Part of this 
paper was written while C.P. was visiting MIT, Cambridge, MA
on a leave of absence from the University of Rome. She
gratefully acknowledges the
hospitality extended to her by Professor I. M. Singer.
She also thanks the organizing committe, and in particular
G. Elliott, for inviting her to participate in part of the special
year dedicated to Operator Algebras held at MSRI during the year
2000-01.

\end{document}